\documentclass[11pt,reqno]{amsart}
\usepackage[utf8x]{inputenc}
\usepackage{amssymb,amsmath,latexsym}
\usepackage{epsfig}

\usepackage{color}

\usepackage{epic,eepic,wrapfig,color, ifthen}
\usepackage{url}
\newtheorem{thm}{Theorem}[section]
\newtheorem{cor}[thm]{Corollary}
\newtheorem{lemma}[thm]{Lemma}
\newtheorem{prop}[thm]{Proposition}
\newtheorem{defn}[thm]{Definition}

\newtheorem{remark}[thm]{Remark}
\newtheorem{example}[thm]{Example}
\numberwithin{equation}{section}

\newcommand{\formula}[2][nolabel]
{\ifthenelse{\equal{#1}{nolabel}}
 {\begin{align*} #2 \end{align*}}
 {\ifthenelse{\equal{#1}{}}
  {\begin{align} #2 \end{align}}
  {\begin{align} \label{#1} #2 \end{align}}
 }
}

\def\pf{{\medskip\noindent {\bf Proof. }}}
\def\qed{{\hfill $\Box$ \bigskip}}
\def\R{{\mathbb R}}
\def\P{{\mathbb P}}
\def\E{{\mathbb E}}
\def\1{{\bf 1}}

\DeclareMathOperator{\dist}{dist}

\renewcommand{\bar}{\overline}

\newcommand{\cal}[1]{\mathcal{#1}}
  
\def\sD {{\cal D}}  
  \def\sI {{\cal I}}
\def\sJ {{\cal J}}  \def\sL {{\cal L}}
\def\sM {{\cal M}} \def\sN {{\cal N}}

 \def\bB {{\mathbb B}}

\def\bbeta {{\boldsymbol \beta}}
\def\bth {{\boldsymbol \theta}}

\def\R {{\mathbb R}}
\def\N {{\mathbb N}}

\def\nn{\nonumber}

\def\E{{\mathbb E}}
\def\P{{\mathbb P}}

\def\bea{\begin{align*}}
\def\eea{\end{align*}}
\def\bee{\begin{equation}}
\def\eee{\end{equation}}

\newcommand{\p}{\rho}

\makeatletter
\@addtoreset{equation}{section}

\makeatother

\begin{document}
\bibliographystyle{plain}

\title[Estimates on the tail probabilities of subordinators]
{ \bf Estimates on the tail probabilities of subordinators and applications to general time fractional equations}

\author
{Soobin Cho \quad and \quad Panki Kim}

\address[Cho]{Department of Mathematical Sciences,
Seoul National University,
Building 27, 1 Gwanak-ro, Gwanak-gu
Seoul 08826, Republic of Korea}
\curraddr{}
\email{soobin15@snu.ac.kr}

\address[Kim]{Department of Mathematical Sciences,
Seoul National University,
Building 27, 1 Gwanak-ro, Gwanak-gu
Seoul 08826, Republic of Korea}\thanks{This research is  supported by the National Research Foundation of Korea(NRF) grant funded by the Korea government(MSIP) (No. 2016R1E1A1A01941893).
}
\curraddr{}
\email{pkim@snu.ac.kr}

 \date{}

\maketitle

\begin{abstract}
In this paper, we study estimates on tail probabilities $\P(S_r \ge t)$ of several classes of subordinators under mild assumptions on the tail of its L\'evy measure. As an application of that result, we obtain two-sided estimates for fundamental solutions of general homogeneous time fraction equations including those with Dirichlet boundary conditions.

\end{abstract}

\bigskip


\bigskip

 
 \bigskip

\section{Introduction}

\subsection{Motivation}

\smallskip

The fractional-time diffusion equation $\partial_t^\beta u = \Delta u$ $(0<\beta<1)$ has been used in various fields to model the diffusions on sticky and trapping environment. Here, $\partial^\beta_t$ is the Caputo derivative of order $\beta$ which is defined as
\begin{align*}
\partial^{\beta}_t u(t) := \frac{1}{\Gamma(1-\beta)} \frac{d}{dt}\int_0^t (t-s)^{-\beta}(u(s)-u(0)) ds,
\end{align*}
where $\Gamma$ is the gamma function defined as $\Gamma(z):=\int_0^\infty x^{z-1}e^{-x}dx$. Motivated by this equation, following \cite{Ch}, we consider the following generalized fractional-time derivatives. Let $w:(0, \infty) \to [0, \infty)$ be a function which satisfies the following condition.

\smallskip

\noindent {\bf (Ker.)} $w$ is a right continuous non-increasing function satisfying $\lim_{s \to 0+} w(s) = \infty$, $\lim_{s \to \infty} w(s) = 0$ and $\int_0^{\infty} \min\{1, s\} (-dw(s)) < \infty$.

\begin{defn}
For a function $u: [0, \infty) \to \R$, the generalized fractional-time derivative $\partial^{w}_t$ with respect to the kernel $w$ is given by
\begin{align*}
\partial^{w}_t u(t) := \frac{d}{dt} \int_0^t w(t-s)(u(s)-u(0))ds,
\end{align*}
whenever the above integral makes sense.
\end{defn}

For example, if $w(t)= \frac{1}{\Gamma(1-\beta)}t^{-\beta}$ for some $0<\beta<1$, then the fractional-time derivative $\partial^w_t$ is nothing but the Caputo derivative of order $\beta$.

In \cite{Ch}, Zhen-Qing Chen established the probabilistic representation for the fundamental solution of generalized fractional-time equation $\partial^{w}_t u(t) = \sL u$ where $\sL$ is the infinitesimal generator of some uniformly bounded strongly continuous semigroup in a Banach space. This procedure can be described as follows: For a given function $w$ satisfying condtion {\bf (Ker.)}, we define a Bernstein function $\phi$ by

\begin{align}\label{e:defphi}
\phi(\lambda) := \int_0^{\infty} (1 - e^{-\lambda s}) (-dw(s)) \qquad \text{for all} \;\; \lambda \ge 0.
\end{align}

\noindent Since $|1-e^{-\lambda s}| \le (1+\lambda)\min\{1, s\}$, we see from {\bf (Ker.)} that $\phi$ is well-defined. Let $\{S_r, r \ge 0\}$ be a subordinator (non-negative valued L\'evy process with $S_0=0$) whose Laplace exponent is given by \eqref{e:defphi}, that is, $\phi(\lambda)= -\log \E\big[\exp(-\lambda S_1)\big]$ for all $\lambda \ge 0$. Then, define its inverse as $E_t:= \inf\{r>0 : S_r >t \}$ for $t>0$. Since condition {\bf (Ker.)} holds, we have $\lim_{s \to 0+} w(s) = \infty$ so that $S_r$ is not a compounded Poisson process. Therefore, almost surely, $r \mapsto S_r$ is strictly increasing and hence $t \mapsto E_t$ is continuous. Denote by $T_t$ the semigroup corresponding to the generator $\sL$ in a Banach space. Then, for every $f \in \sD(\sL)$, the unique solution (in some suitable sense) to the following general homogeneous fractional time equation
\begin{align}\label{e:tfe}
\partial^w_t u(t,x) = \sL u(t,x) \quad \text{with} \;\; u(0,x) = f(x)
\end{align}
is given by 
\begin{align}\label{e:soln}
u(t,x) = \E_x\big[T_{E_t}f(x)\big].
\end{align}

In \cite{CKKW}, the second named author, jointly with Zhen-Qing Chen, Takashi Kumagai and Jian Wang, proved that when $T_t$ is the transition semigroup of a symmetric strong Markov process, \eqref{e:soln} is the unique weak solution to equation \eqref{e:tfe} (see \cite[Theorem 2.4]{CKKW} for a precise statement).  Moreover, they obtained two-sided estimates for the fundamental solution under the condition that $\phi$ satisfies ${\bf WS}(\alpha_1, \alpha_2)$ for some $0<\alpha_1 \le \alpha_2 <1$ (see Definition \ref{d:wsp} for the definition of ${\bf WS}(\alpha_1, \alpha_2)$). The key ingredients to obtain those estimates were the estimates on tail probabilities $\P(S_r \ge t)$ and $\P(S_r \le t)$ established in \cite{JP, Mi}. Particularly, the weak scaling conditions for $\phi$ were needed to get sharp estimates on $\P(S_r \ge t)$. 

In this paper, we study estimates on upper tail probabilites $\P(S_r \ge t)$ of a general class of subordinators. Our results cover some cases when the lower scaling index $\alpha_1$ of $\phi$ is $0$ and the upper scaling index $\alpha_2$ of $\phi$ is $1$. Indeed, we will see that the lower scaling index has no role in tail probability estimates. On the other hand, when the upper scaling index is $1$, various phenomena can arise in the asymptotic behaviors of $\P(S_r \ge t)$ as $t \to \infty$. To assort those phenomena, we impose conditions on the tail measure $w$ instead of the Laplace exponent $\phi$ and then obtain estimates on $\P(S_r \ge t)$ under each condition. More precisely, we will consider the three cases: (i) $w$ is a polynomial decaying function; (ii) $w$ decreases subexponentially or exponentially; (iii) $w$ is finitely supported. (See, Section \ref{s:subordinator} for details.)

As applications to these tail probability estimates, we then establish two-sided estimates for fundamental solution of a general time fractional equation including the ones with the Dirichlet boundary condition, given by \eqref{e:dfundamental}. 

\vspace{2mm}

\subsection{Settings}

In this subsection, we introduce the notions of the fundamental solution for a general time fractional equation and the weak scaling properties for non-negative function. Then, we list our main assumptions in this paper.

\smallskip

Let $(M, \p, m)$ be a separable locally compact Hausdorff metric measure space and $D \subset M$ be an open subset. Let $\{T^D_t, t \ge 0\}$ be a uniformly bounded strongly continuous semigroup with infinitesimal generator $(\sL^D, \sD(\sL^D))$ in some Banach space $(\bB, \lVert \cdot \rVert)$. Let $w$ be a function satisfying condition {\bf (Ker.)}. Then, we consider the following time fractional equation with Dirichlet boundary condition.
\begin{align}\label{e:dtfe}
\begin{cases}
\partial^{w}_t u(t,x) = \sL^D u(t,x),    & x \in D, \;\; t>0, \\
u(0, x) = f(x), & x \in D, \\
u(t, x) = 0,  & \text{vanishes continuously on} \; \partial D \;\; \text{for all} \;\; t>0.
\end{cases}
\end{align}
Examples of the problem \eqref{e:dtfe} can be found in \cite{CMN, MNV}. If we overlook the boundary condition, then it is established in \cite[Theorem 2.3]{Ch} that for all $f \in \sD(\sL^D)$, $u(t,x):=\E[T^D_{E_t}f(x)]$ is a unique solution to \eqref{e:dtfe} in the following sense: 

\vspace{2mm}

\noindent(i) $\sup_{t>0} \lVert u(t, \cdot) \rVert < \infty$, $x \mapsto u(t,x)$ is in $\sD(\sL^D)$ for each $t \ge 0$ with $\sup_{t \ge 0} \lVert \sL^D u(t, \cdot)\rVert < \infty$, and both $t \mapsto u(t, \cdot)$ and $t \mapsto \sL^D u(t, \cdot)$ are continuous in $(\bB, \lVert \cdot \rVert)$;

\noindent(ii) for every $t>0$, $I^w_t[u]:= \int_0^t w(t-s) (u(s,x)-f(x)) ds$ is absolutely convergent in $(\bB, \lVert \cdot \rVert)$ and
\begin{align*}
\lim_{\delta \to 0} \frac{1}{\delta} (I^w_{t + \delta}[u] - I^w_t[u]) = \sL^D u(t,x) \quad \text{in} \;\; (\bB, \lVert \cdot \rVert).
\end{align*}

\smallskip

\noindent Indeed, we will see that if $\{T^D_t, t \ge 0\}$ admits a transition density enjoying certain types of estimates, then the solution $u(t,x)$ satisfies the following boundary condition (see Corollary \ref{c:boundary} for a precise statement). 

\smallskip

\noindent (iii) if $f$ is bounded, then for all $t>0$, $x \mapsto u(t,x)$ vanishes continuously on $\partial D$.

\vspace{3mm}

As discussed in \cite{CKKW}, if the semigroup $\{T_t^D, t \ge 0 \}$ has a transition density $q(t,x,y)$ with respect to $m$ on $M$, for any function $f \in \sD(\sL^D)$,
\begin{align*}
u(t, x) & = \E_x[T_{E_t}^Df(x)] =\int_0^{\infty} T^D_r f(x) d_r \P(E_t \le r) = \int_0^{\infty} T^D
_r f(x) d_r\P(S_r \ge t) \\
& = \int_0^{\infty}\int_M f(y)q(r,x,y) m(dy) d_r \P(S_r \ge t)\\
& = \int_M f(y) \left(\int_0^{\infty} q(r,x,y) d_r \P(S_r \ge t) \right) m(dy).
\end{align*}
Therefore, it is natural to say that
\begin{align}\label{e:dfundamental}
p(t,x,y) := \int_0^{\infty} q(r,x,y) d_r \P(S_r \ge t)
\end{align}
is the fundamental solution to the equation \eqref{e:dtfe}.

\smallskip

Next, we introduce the weak scaling properties for non-negative functions.

\begin{defn}\label{d:wsp}
{\rm 
Let $f:(0, \infty) \to [0, \infty)$ be a given function and $\alpha_1, \alpha_2 \in \R$ and $c_0>0$ be given constants.

\noindent (1) We say that $f$ satisfies ${\bf LS^0}(\alpha_1, c_0)$ (resp. ${\bf LS^{\infty}}(\alpha_1, c_0)$) if there exists a constant $c_1>0$ such that
\begin{align*}
\frac{f(R)}{f(r)} \ge c_1 \left(\frac{R}{r}\right)^{\alpha_1} \quad \text{for all} \;\; r \le R \le c_0 \;\; \text{(resp. for all} \;\; c_0 \le r \le R).
\end{align*}

\noindent (2) We say that $f$ satisfies ${\bf US^0}(\alpha_2, c_0)$ (resp. ${\bf US^{\infty}}(\alpha_2, c_0)$) if there exists a constant $c_2>0$ such that
\begin{align*}
\frac{f(R)}{f(r)} \le c_2 \left(\frac{R}{r}\right)^{\alpha_2} \quad \text{for all} \;\; r \le R \le c_0 \;\; \text{(resp. for all} \;\; c_0 \le r \le R).
\end{align*}

\noindent (3) If $f$ satisfies both ${\bf LS^0}(\alpha_1,c_0)$ and ${\bf US^{0}}(\alpha_2, c_0)$ (resp. ${\bf LS^\infty}(\alpha_1,c_0)$ and ${\bf US^{\infty}}(\alpha_2, c_0)$), we say that $f$ satisfies ${\bf WS^0}(\alpha_1,\alpha_2, c_0)$ (resp. ${\bf WS^\infty}(\alpha_1,\alpha_2, c_0)$).
Moreover, if $f$ satisfies both ${\bf WS^0}(\alpha_1,\alpha_2, c_0)$ and ${\bf WS^\infty}(\alpha_1,\alpha_2, c_0)$, then we say that $f$ satisfies ${\bf WS}(\alpha_1,\alpha_2)$.

\smallskip

}
\end{defn}

Throughout this paper, we always assume that the kernel $w$ satisfies condition {\bf (Ker.)}. Here, we enumerate our main assumptions for $w$.

\vspace{2mm}

\noindent {\bf (S.Poly.)($t_s$)} There exist constants $t_s>0$ and $\delta_1 >0$ such that $w$ satisfies ${\bf LS^0}(-\delta_1, t_s)$;

\smallskip

\noindent {\bf (L.Poly.)} There exists a constant $\delta_2 > 0$ such that $w$ satisfies ${\bf LS^\infty}(-\delta_2, 1)$;

\smallskip

\noindent {\bf (Sub.)($\bbeta, \bth$)} There exist constants $c_0, \theta>0$ and $\beta \in (0,1]$ such that 
\begin{align*}
w(t) \le c_0\exp(-\theta t^\beta) \quad \text{for all} \;\; t \ge 1.
\end{align*}

\smallskip

\noindent {\bf (Trunc.)($t_f$)} There exists a constant $t_f>0$ such that 

(i) $w(t)>0$ for $0<t<t_f$ and $w(t_f) = 0$;

(ii) $w$ is bi-Lipschitz continuous on $[t_f/4,t_f]$, i.e. there exists a constant $K \ge 1$ such that
\begin{align*}
K^{-1}|t-s| \le |w(t)-w(s)| \le K|t-s|, \quad \text{for all} \;\; t_f/4 \le s \le t \le t_f;
\end{align*}

(iii) there exists a constant $\delta_3>0$ such that $w$ satisfies ${\bf LS^0}(-\delta_3, t_f/2)$.

\vspace{2mm}

\begin{remark}
{\rm (1) Condition {\bf (S.Poly.)($t_s$)} implies that the corresponding Laplace exponent $\phi$ satisfies ${\bf US^\infty}(\min\{\delta_1,1 \}, 1)$. Conversely, if $\phi$ satisfies ${\bf US^\infty}(\delta_1, 1)$ for some $\delta_1<1$, then there exists a constant $t_s>0$ such that condition {\bf (S.Poly.)($t_s$)} holds with constant $\delta_1$. Analogously, condition {\bf (L.Poly.)} implies that $\phi$ satisfies ${\bf US^0}(\min\{\delta_2,1 \}, 1)$ and if $\phi$ satisfies ${\bf US^0}(\delta_2,1)$ with $\delta_2<1$, then condition {\bf (L.Poly.)} holds. (See, Lemma \ref{l:phiHw}.) 

(2) If condition {\bf (L.Poly.)} or {\bf (Sub.)($\bbeta, \bth$)} holds, then we can replace the constant $1$ with arbitrary positive constant since $w$ is a monotone function. However, we can not replace the constant $t_s$ in condition {\bf (S.Poly.)($t_s$)} with other positive constants in general. For instance, if $w(t)=(t^{-1/2}-1)\1_{(0,1]}(t)$, then we can only take $t_s$ strictly smaller than $1$. Moreover, the constant $t_f$ in condition {\bf (Trunc.)($t_f$)} is uniquely determined by its first condition.
}
\end{remark}

\vspace{2mm}

{\bf Notations}: In this paper, we use the symbol ``$:=$'' to denote a definition, 
which is read as ``is defined to be.'' For $a,b\in \R$, we use the notations $a\wedge b:=\min\{a,b\}$ and $a\vee b:=\max\{a,b\}$. For $x \in \R$, we define $\log^+ x:= 0 \vee \log x$ and $\lfloor x \rfloor:= \max \{ n \in \mathbb{Z} : x \ge n \}$. We denote by $\partial_t$ the partial derivative with respect to the variable $t$.
  
The notation $f(x) \asymp g(x)$ means that there exist constants $c_1,c_2>0$ such that $c_1g(x)\leq f (x)\leq c_2 g(x)$ for the specified range of the variable $x$. The notation $f(x) \lesssim g_1(x) + g_2(x)h(cx)$ (resp.  $f(x) \gtrsim g_1(x) + g_2(x)h(cx)$) means that there exist constants $c_1,c_2>0$ such that
$$f(x) \le c_1\big(g_1(x) + g_2(x)h(c_2x)\big) \quad \text{(resp. } \; f(x) \ge c_1\big(g_1(x) + g_2(x)h(c_2x)\big) ),$$
for the specified range of $x$. Then, the notation $f(x) \simeq g_1(x) + g_2(x)h(cx)$ means that both $f(x) \lesssim g_1(x) + g_2(x)h(cx)$ and $f(x) \gtrsim g_1(x) + g_2(x)h(cx)$ hold for the specified range of $x$.

For a subset $D$ of some metric space $(M, \p)$, we let $\text{diam}(D):= \sup_{u,v \in D} \p(u,v)$ and $\delta_D(x):= \sup_{z \in D} \p(x,z)$ for $x \in D$. Then, for $x,y \in D$, we define 
\begin{align}\label{e:prmM}
&\delta_*(x,y):= \delta_D(x) \delta_D(y), \quad \delta_\wedge(x,y):= \delta_D(x) \wedge \delta_D(y) \quad \text{and} \quad \delta_\vee(x,y):= \delta_D(x) \vee \delta_D(y).
\end{align}

Lower case letters $c$'s without subscripts denote strictly positive
constants  whose values
are unimportant and which  may change even within a line, while values of lower case letters with subscripts
$c_i, i=0,1,2,  \dots$, are fixed in each statement and proof, 
and the labeling of these constants starts anew in each proof.

\vspace{1mm}

\subsection{Some toy models with explicit Dirichlet estimates.}\label{s:special}

Our general estimates on the fundamental solution include a term which is described in an integral form. (See, \eqref{e:defD}.) However, in many applications, we can obtain explicit forms of them. We first represent some special versions of our results which can be described explicitly.

\smallskip

 Suppose that the operator $(\sL^D, \sD(\sL^D))$ on $(D, \p, m)$ admits a heat kernel $q(t,x,y)$ with respect to the measure $m$. 
We further assume that one of the following assumptions holds for all $(t,x,y) \in (0,\infty) \times D \times D$.

\smallskip

\noindent {\bf (J1)} $\text{diam}(D)<\infty$ and there exist constants $\alpha, d>0$ and $\lambda>0$ such that
\begin{align*}
q(t,x,y) \asymp 
\begin{cases} 
\displaystyle \left(1 \wedge \frac{\delta_D(x)}{t^{1/\alpha}}\right)^{\alpha/2}\left(1 \wedge \frac{\delta_D(y)}{t^{1/\alpha}}\right)^{\alpha/2} \left( t^{-d/\alpha} \wedge \frac{t}{\p(x,y)^{d+\alpha}}\right), & \mbox{if}  \;\; 0< t\le 1; \\
e^{-\lambda t}\delta_D(x)^{\alpha/2}\delta_D(y)^{\alpha/2}, & \mbox{if} \;\; t \ge 1;
\end{cases}
\end{align*}

\noindent {\bf (J2)} There exist constants $\alpha>0$ and $d>0$ such that for all $t>0$,
\begin{align*}
q(t,x,y) \asymp 
\displaystyle \left(1 \wedge \frac{\delta_D(x)}{t^{1/\alpha}}\right)^{\alpha/2}\left(1 \wedge \frac{\delta_D(y)}{t^{1/\alpha}}\right)^{\alpha/2}\left( t^{-d/\alpha} \wedge \frac{t}{\p(x,y)^{d+\alpha}}\right);
\end{align*}

\noindent {\bf (J3)} There exist constants $\alpha>0$ and $d>0$ such that for all $t>0$,
\begin{align*}
q(t,x,y) \asymp
\displaystyle \left(1 \wedge \frac{\delta_D(x)}{t^{1/\alpha} \wedge 1}\right)^{\alpha/2}\left(1 \wedge \frac{\delta_D(y)}{t^{1/\alpha} \wedge 1}\right)^{\alpha/2}\left( t^{-d/\alpha} \wedge \frac{t}{\p(x,y)^{d+\alpha}}\right);
\end{align*}

\noindent {\bf (J4)} $\text{diam}(D)<\infty$ and there exist constants $\alpha>1, d>0$ and $\lambda>0$ such that
\begin{align*}
q(t,x,y) \asymp 
\begin{cases} 
\displaystyle \left(1 \wedge \frac{\delta_D(x)}{t^{1/\alpha}}\right)^{\alpha-1}\left(1 \wedge \frac{\delta_D(y)}{t^{1/\alpha}}\right)^{\alpha-1} \left( t^{-d/\alpha} \wedge \frac{t}{\p(x,y)^{d+\alpha}}\right), & \mbox{if} \;\; 0< t\le 1; \\
e^{-\lambda t}\delta_D(x)^{\alpha-1}\delta_D(y)^{\alpha-1}, & \mbox{if} \;\; t \ge 1;
\end{cases}
\end{align*}

\noindent {\bf (D1)} $\text{diam}(D)<\infty$ and there exist positive constants $\alpha>1, d>0$ and $\lambda>0$ such that
\begin{align*}
q(t,x,y) \simeq 
\begin{cases} 
\displaystyle \left(1 \wedge \frac{\delta_D(x)}{t^{1/\alpha}}\right)^{\alpha/2}\left(1 \wedge \frac{\delta_D(y)}{t^{1/\alpha}}\right)^{\alpha/2} t^{-d/\alpha} \exp\bigg(-c \frac{\p(x,y)^{\alpha/(\alpha-1)}}{t^{1/(\alpha-1)}}\bigg), & \mbox{if}  \;\; 0< t\le 1; \\
e^{-\lambda t}\delta_D(x)^{\alpha/2}\delta_D(y)^{\alpha/2}, & \mbox{if} \;\; t \ge 1;
\end{cases}
\end{align*}

\noindent {\bf (D2)} There exist positive constants $\alpha>1$ and $d>0$ such that for all $t>0$,
\begin{align*}
q(t,x,y) \simeq 
\displaystyle \left(1 \wedge \frac{\delta_D(x)}{t^{1/\alpha}}\right)^{\alpha/2}\left(1 \wedge \frac{\delta_D(y)}{t^{1/\alpha}}\right)^{\alpha/2} t^{-d/\alpha} \exp\bigg(-c \frac{\p(x,y)^{\alpha/(\alpha-1)}}{t^{1/(\alpha-1)}}\bigg) ;
\end{align*}

\noindent {\bf (D3)} There exist positive constants $\alpha>1$ and $d>0$ such that for all $t>0$,
\begin{align*}
q(t,x,y) \simeq 
\displaystyle \left(1 \wedge \frac{\delta_D(x)}{t^{1/\alpha} \wedge 1}\right)^{\alpha/2}\left(1 \wedge \frac{\delta_D(y)}{t^{1/\alpha} \wedge 1}\right)^{\alpha/2} t^{-d/\alpha} \exp\bigg(-c \frac{\p(x,y)^{\alpha/(\alpha-1)}}{t^{1/(\alpha-1)}}\bigg).
\end{align*}

\vspace{3mm}

An open subset $D \subset \R^d$ $(d \ge 2)$ is said to be a $C^{1,1}$ open set if there exist a localization radius $R_0>0$ and a constant $\Lambda>0$ such that for every $z \in \partial D$, there is a $C^{1,1}$ function $\Gamma:\R^{d-1} \to \R$ satisfying $\Gamma(0) = 0, \nabla \Gamma(0) = (0,...,0), \lVert \Gamma \rVert_{\infty} \le \Lambda, |\nabla \Gamma(y) - \nabla \Gamma(z)| \le \Lambda |y-z|$ and an orthonormal coordinate system $CS_z: x = (\tilde{x}, x_d) :=(x_1, ..., x_{d-1}, x_d)$ with origin at $z$ such that
\begin{align*}
D \cap B(z,R_0) = \{ x \in B(0, R_0) \;\; \text{in} \;\; CS_z  : x_d > \Gamma(\tilde{x})\}.
\end{align*}
A $C^{1,1}$ open set in $\R$ is the union of disjoint intervals so that the minimum of their lengths and the distances between them is positive.

\begin{remark}\label{r:examples}
{\rm
When $M$ is $\R^d$, $\p$ is the usual metric on $\R^d$ and $m$ is the Lebesgue measure, there are many examples of generators $(\sL^D, \sD(\sL^D))$ on $(D, \p, m)$ which admit a transition density satisfying one of the estimates among {\bf (J1)}, {\bf (J2)}, {\bf (J3)}, {\bf (J4)}, {\bf (D1)}, {\bf (D2)} and {\bf (D3)}. For instance, if $\sL^D$ is a generator of a killed symmetric $\alpha$-stable process with $0<\alpha<2$ or a censored $\alpha$-stable process with $1<\alpha<2$ and $D \subset \R^d$ is a bounded $C^{1,1}$ open set, then estimate {\bf (J1)} or {\bf (J4)} holds, respectively. (See, \cite{CKS, CKS2}.) Else if $\sL^D$ is a generator of a killed symmetric $\alpha$-stable process with $0<\alpha<2\wedge d$ and $D$ is a half space-like $C^{1,1}$ open set or exterior of a bounded $C^{1,1}$ open set, then estimate {\bf (J2)} or {\bf (J3)} holds, respectively. (See, \cite[Theorems 5.4 and 5.8]{BGR}.) Moreover, when $d \ge 3$, $\sL$ is the Dirichlet laplacian on $D$ and $D \subset \R^d$ is a bounded connected $C^{1,1}$ open set or half space-like $C^{1,1}$ open set or exterior of a bounded $C^{1,1}$ open set, then estimate {\bf (D1)} or {\bf (D2)} or {\bf (D3)} holds with $\alpha=2$, respectively. (See, \cite{CK, So, Zh, Zh2}.)
}
\end{remark}

Recall that $\delta_*$, $\delta_\wedge$  and $\delta_\vee$ are defined in \eqref{e:prmM}.
For $\alpha>0$, we define two auxiliary functions $F^\alpha_k, F^\alpha_c:\R \times (0,\infty) \times D \times D \to [0,\infty)$ as follows.

\begin{align*}
&F^\alpha_k(s,t,x,y):= \\
&  \begin{cases}
\vspace{1mm}
\1_{\{\delta_*(x,y)^{\alpha/2} \le \phi(t^{-1})^{-1} \}} \big( \p(x,y)^\alpha \vee \delta_*(x,y)^{\alpha/2}\big) \phi(t^{-1})^{-s/\alpha},    & \text{if } \;   s < 0;\\
\vspace{1mm}
\1_{\{\delta_*(x,y)^{\alpha/2} \le \phi(t^{-1})^{-1} \}} \big(\p(x,y)^\alpha \vee \delta_*(x,y)^{\alpha/2}\big)\log^+ \left(\displaystyle\frac{2\phi(t^{-1})^{-1}}{\p(x,y)^{\alpha} \vee \delta_\vee(x,y)^\alpha }\right), & \text{if } \; s= 0 ;\\
\vspace{1mm}
\1_{\{\delta_*(x,y)^{\alpha/2} \le \phi(t^{-1})^{-1} \}}\big(\p(x,y)^{\alpha-s} \vee \delta_*(x,y)^{\alpha/2} \delta_\vee(x,y)^{-s} \big),    & \text{if } \; s< \displaystyle\frac{\alpha}{2}  ;\\
\vspace{1mm}
\1_{\{\delta_*(x,y)^{\alpha/2} \le \phi(t^{-1})^{-1} \}} \bigg(\p(x,y)^{\alpha/2}
+ \delta_\wedge(x,y)^{\alpha/2}\log \left(\displaystyle\frac{\p(x,y) \vee 2\delta_\vee(x,y)}{\p(x,y) \vee \delta_\wedge(x,y)}\right)\bigg),    & \text{if } \; s= \displaystyle\frac{\alpha}{2} ;\\
\vspace{1mm}
\1_{\{\delta_*(x,y)^{\alpha/2} \le \phi(t^{-1})^{-1} \}}\big(\p(x,y)^{\alpha-s} \vee \delta_\wedge(x,y)^{\alpha-s}\big),    & \text{if } \;  \displaystyle\frac{\alpha}{2} < s < \alpha  ;\\
\vspace{1mm}
1+ \log^+ \left(\displaystyle\frac{2\phi(t^{-1})^{-1} \wedge 2\delta_\wedge(x,y)^{\alpha}}{\p(x,y)^{\alpha} }\right), & \text{if } \;\;  s=\alpha ;\\
\p(x,y)^{\alpha-s},  & \text{if } \;\; s>\alpha.
\end{cases}
\end{align*}

\begin{align*}
&F^\alpha_c(s,t,x,y):=\\
&  \begin{cases}
\vspace{1mm}
\1_{\{\delta_*(x,y)^{\alpha/2} \le \phi(t^{-1})^{-1} \}} \big( \p(x,y)^{2\alpha -2} \vee \delta_*(x,y)^{\alpha-1}\big) \phi(t^{-1})^{-(2-\alpha-s)/\alpha},   & \text{if } \; s<2-\alpha  ;\\
\vspace{1mm}
\1_{\{\delta_*(x,y)^{\alpha/2} \le \phi(t^{-1})^{-1} \}} \big( \p(x,y)^{2\alpha-2} \vee \delta_*(x,y)^{\alpha-1}\big)\log^+\left(\displaystyle\frac{2\phi(t^{-1})^{-1}}{\p(x,y)^\alpha \vee \delta_\vee(x,y)^\alpha}\right) ,   & \text{if } \; s=2-\alpha  ;\\
\vspace{1mm}
\1_{\{\delta_*(x,y)^{\alpha/2} \le \phi(t^{-1})^{-1} \}} \big(\p(x,y)^{\alpha-s}  \vee \delta(x,y)^{\alpha-1}\delta_\vee(x,y)^{2-\alpha-s} \big),   & \text{if } \; 2-\alpha<s<1  ;\\
\vspace{1mm}
\1_{\{\delta_*(x,y)^{\alpha/2} \le \phi(t^{-1})^{-1} \}} \bigg(\p(x,y)^{\alpha-1} +   \delta_\wedge(x,y)^{\alpha-1} \log \left(\displaystyle\frac{\p(x,y) \vee 2\delta_\vee(x,y)}{\p(x,y) \vee \delta_\wedge(x,y)}\right)\bigg),   & \text{if } \; s=1  ;\\
\vspace{1mm}
\1_{\{\delta_*(x,y)^{\alpha/2} \le \phi(t^{-1})^{-1} \}} \big(\p(x,y)^{\alpha-s} \vee \delta_\wedge(x,y)^{\alpha-s}\big),    & \text{if } \; 1<s<\alpha ;\\
\vspace{1mm}
1+ \log^+ \left(\displaystyle\frac{2\phi(t^{-1})^{-1} \wedge 2\delta_\wedge(x,y)^{\alpha}}{\p(x,y)^{\alpha} }\right),    & \text{if } \;   s = \alpha  ;\\
\vspace{1mm}
\p(x,y)^{\alpha-s}, & \text{if } \;  s>\alpha.
\end{cases}
\end{align*}

We also define 
\begin{align}\label{barphi}
\bar{\phi}_{\alpha}(\lambda):= \inf \{ s>0 : s^{\alpha}\phi(s)^{-1} \ge \lambda \} \qquad \text{for} \;\; \lambda \ge 0.
\end{align}

Recall that for an integral kernel $w$ satisfying condition {\bf (Ker.)}, the fundamental solution $p(t,x,y)$ of the general fractional-time equation \eqref{e:dtfe} is given by \eqref{e:dfundamental}. We first give the small time estimates for $p(t,x,y)$ under condition {\bf (S.Poly.)($t_s$)}.

\begin{thm}\label{t:specialsmall}
Assume that $w$ satisfies conditions {\bf (Ker.)} and {\bf (S.Poly.)($t_s$)}. \\
Then, the follwing estimates for $p(t,x,y)$ hold for all $(t,x,y) \in (0, t_s] \times D \times D$.

\smallskip

\noindent (i) (Near diagonal estimates) Suppose that $\phi(t^{-1})\p(x,y)^{\alpha} \le 1/(4e^2)$.

\indent (a) If one of the estimates among {\bf (J1)}, {\bf (J2)}, {\bf (J3)}, {\bf (D1)}, {\bf (D2)} and {\bf (D3)} holds, then we have
\begin{align}\label{smallon}
p(t, x, y) &\asymp \left(1 \wedge \frac{\delta_*(x,y)}{\phi(t^{-1})^{-2/\alpha}}\right)^{\alpha/2}\phi(t^{-1})^{d/\alpha} + w(t)\left(1 \wedge \displaystyle\frac{\delta_*(x,y)}{\p(x,y)^{2}}\right)^{\alpha/2}F^\alpha_k(d, t,x,y).
\end{align}

\indent (b) Otherwise, if {\bf (J4)} holds, then we have
\begin{align*}
p(t, x, y) &\asymp \left(1 \wedge \frac{\delta_*(x,y)}{\phi(t^{-1})^{-2/\alpha}}\right)^{\alpha-1}\phi(t^{-1})^{d/\alpha} + w(t)\left(1 \wedge \displaystyle\frac{\delta_*(x,y)}{\p(x,y)^{2}} \right)^{\alpha-1}F^\alpha_c(d, t,x,y).
\end{align*}

\smallskip

\noindent  (ii) (Off diagonal estimates)  Suppose that $\phi(t^{-1})\p(x,y)^{\alpha} > 1/(4e^2)$. 

\indent (a) If {\bf (J1)} or {\bf (J2)} or {\bf (J3)} holds, then we have
\begin{align}\label{smalloff}
p(t,x,y) \asymp \left(1 \wedge \frac{\delta_D(x)}{\phi(t^{-1})^{-1/\alpha}}\right)^{\alpha/2}\left(1 \wedge \frac{\delta_D(y)}{\phi(t^{-1})^{-1/\alpha}}\right)^{\alpha/2}\frac{\phi(t^{-1})^{-1}}{\p(x,y)^{d+\alpha}}.
\end{align}

\indent (b) If {\bf (J4)} holds, then we have
\begin{align*}
p(t,x,y) \asymp \left(1 \wedge \frac{\delta_D(x)}{\phi(t^{-1})^{-1/\alpha}}\right)^{\alpha-1}\left(1 \wedge \frac{\delta_D(y)}{\phi(t^{-1})^{-1/\alpha}}\right)^{\alpha-1}\frac{\phi(t^{-1})^{-1}}{\p(x,y)^{d+\alpha}}.
\end{align*}

\indent (c) Otherwise, if {\bf (D1)} or {\bf (D2)} or {\bf (D3)} holds, then we have
\begin{align}\label{smalloff2}
p(t,x,y) \simeq \left(1 \wedge \frac{\delta_D(x)}{\phi(t^{-1})^{-1/\alpha}}\right)^{\alpha/2}\left(1 \wedge \frac{\delta_D(y)}{\phi(t^{-1})^{-1/\alpha}}\right)^{\alpha/2} \phi(t^{-1})^{d/\alpha} \exp\bigg(-c t \bar{\phi}_{\alpha}\big((\frac{\p(x,y)}{t})^{\alpha}\big)\bigg),
\end{align}
where the function $\bar{\phi}_{\alpha}$ is defined as \eqref{barphi}.

\end{thm}

\vspace{2mm}

Next, under condition {\bf (L.Poly.)}, we get the large time estimates for $p(t,x,y)$. Hereinafter, we let $R_D:=\text{diam}(D)$ and $T_D:=[\phi^{-1}(4^{-1}e^{-2}R_D^{-\alpha})]^{-1}$.

\begin{thm}\label{t:speciallarge}
Assume that $w$ satisfies conditions {\bf (Ker.)} and {\bf (L.Poly.)}. Then, for every fixed $T>0$, the follwing estimates hold for all $(t,x,y) \in [T, \infty) \times D \times D$.

\smallskip

\noindent  (i) If {\bf (J1)} or {\bf (D1)} holds and $R_D< \infty$, then we have
\begin{align*}
p(t, x, y) &\asymp w(t) \left(1 \wedge \displaystyle\frac{\delta_*(x,y)}{\p(x,y)^{2}}\right)^{\alpha/2}\bigg( [1 \wedge \delta_*(x,y)^{\alpha/2}] + F^\alpha_k(d,T_D,x,y) \bigg).
\end{align*}

\noindent (ii) If {\bf (J4)} holds and $R_D< \infty$, then we have
\begin{align*}
p(t, x, y) &\asymp w(t)\left(1 \wedge \displaystyle\frac{\delta_*(x,y)}{\p(x,y)^{2}} \right)^{\alpha-1}\bigg( [1 \wedge \delta_*(x,y)^{\alpha-1} ]+ F^\alpha_c(d, T_D,x,y)\bigg).
\end{align*}

\noindent (iii) If {\bf (J2)} holds, then estimates given in \eqref{smallon} and \eqref{smalloff} hold for all $(t,x,y) \in [T, \infty) \times D \times D$.

\noindent (iv) If {\bf (D2)} holds, then estimates given in \eqref{smallon} and \eqref{smalloff2} hold for all $(t,x,y) \in [T, \infty) \times D \times D$.

\noindent (v) Assume that either of the estimates {\bf (J3)} or {\bf (D3)} holds.

\indent(a) If $\phi(t^{-1}) \p(x,y)^\alpha \le 1/(4e^2)$, then we have
\begin{align*}
p(t, x, y) &\asymp \big(1 \wedge \delta_D(x) \big)^{\alpha/2}\big(1 \wedge \delta_D(y) \big)^{\alpha/2} \bigg (\phi(t^{-1})^{d/\alpha} + w(t)G_d^\alpha (t, 1 \vee \p(x,y) ) \bigg) \\
& \quad +  \1_{\{\p(x,y) \le 1\}}w(t)\left(1 \wedge \displaystyle\frac{\delta_*(x,y)}{\p(x,y)^{2}}\right)^{\alpha/2}F^\alpha_k(d,[\phi^{-1}(4^{-1}e^{-2})]^{-1}, x, y),
\end{align*}
where the function $G^\alpha_d(t,l)$ is defined as follows:
\begin{align*}
&G^\alpha_d(t,l):=  \; \begin{cases}
0,  &\text{if } \; d< \alpha ;\\
 \log \left(\displaystyle\frac{2\phi(t^{-1})^{-1}}{l \phi(T^{-1})^{-1}}\right), & \text{if } \;  d=\alpha ;\\
l^{\alpha-d},  & \text{if } \; d>\alpha.
\end{cases}
\end{align*}

\indent (b) If $\phi(t^{-1}) \p(x,y)^\alpha > 1/(4e^2)$, then we have
\begin{align*}
&p(t,x,y) \simeq \big(1 \wedge \delta_D(x)\big)^{\alpha/2}\big(1 \wedge \delta_D(y)\big)^{\alpha/2}\times \begin{cases}
\displaystyle \frac{\phi(t^{-1})^{-1}}{\p(x,y)^{d+\alpha}}, & \text{if {\bf (J3)} holds};\\
\displaystyle \phi(t^{-1})^{d/\alpha} \exp\bigg(-c t \bar{\phi}_{\alpha}\big((\frac{\p(x,y)}{t})^{\alpha}\big)\bigg),  & \text{if {\bf (D3)} holds},
\end{cases}
\end{align*}
where the function $\bar{\phi}_{\alpha}$ is defined as \eqref{barphi}.

\end{thm}

\vspace{1mm}

We mention that under condition {\bf (L.Poly.)}, even if $D$ is bounded so that $q(t,x,y)$ decreases exponentially as $t \to \infty$, the fundamental solution $p(t,x,y)$ is a polynomial decaying function which decreases with the same order as $w$. (See, Theorem \ref{t:speciallarge}(i) and (ii).) We introduce a condition which make $p(t,x,y)$ decreases subexponentially. 

\vspace{2mm}

\noindent {\bf (Sub*.)($\bbeta, \bth$)} There exist constants $c_0>1, \theta>0$ and $\beta \in (0,1)$ such that 
\begin{align*}
c_0^{-1}\exp(-\theta t^\beta) \le w(t) \le c_0\exp(-\theta t^\beta) \quad \text{for all} \;\; t \ge 1.
\end{align*}

Under condition {\bf (Sub*.)($\bbeta, \bth$)}, we obtain estimates for $p(t,x,y)$ which have an exactly the same exponential term as $w$.

\begin{thm}\label{t:specialsub}
Assume that $w$ satisfies conditions {\bf (Ker.)} and {\bf (Sub*.)($\bbeta, \bth$)}. We further assume that {\bf (J1)} or {\bf (J4)} or {\bf (D1)} holds. Then, for every fixed $T>0$, the follwing estimates hold for all $(t,x,y) \in [T, \infty) \times D \times D$.

\smallskip

\noindent  (i) If {\bf (J1)} or {\bf (D1)} holds and $R_D< \infty$, then we have
\begin{align*}
p(t, x, y) &\asymp \exp(-\theta t^\beta) \left(1 \wedge \displaystyle\frac{\delta_*(x,y)}{\p(x,y)^{2}}\right)^{\alpha/2}\bigg( [1 \wedge \delta_*(x,y)^{\alpha/2} ]+ F^\alpha_k(d, T_R,x,y) \bigg).
\end{align*}

\smallskip

\noindent  (ii) If {\bf (J4)} holds and $R_D< \infty$, then we have
\begin{align*}
p(t, x, y) &\asymp \exp(-\theta t^\beta) \left(1 \wedge \displaystyle\frac{\delta_*(x,y)}{\p(x,y)^{2}}\right)^{\alpha-1}\bigg( [1 \wedge \delta_*(x,y)^{\alpha-1}] + F^\alpha_c(d,T_R,x,y) \bigg).
\end{align*}

\end{thm}

\vspace{2mm}

Notice that condition {\bf (Trunc.)($t_f$)} implies condition {\bf (S.Poly.)($t_s$)} with $t_s = t_f/2$. Hence, we obtain the small time estimates ($0<t \le t_f/2$) under condition {\bf (Trunc.)($t_f$)} from Theorem \ref{t:specialsmall}. Here, we give the large time behaviors of $p(t,x,y)$ under condition {\bf (Trunc.)($t_f$)}. 

\begin{thm}\label{t:specialtrunc}
Assume that $w$ satisfies conditions {\bf (Ker.)} and {\bf (Trunc.)($t_f$)}.  Then, the follwing estimates hold for all $(t,x,y) \in [t_f/2, \infty) \times D \times D$. Let $n_t:= \lfloor t/t_f \rfloor +1 \in \N$. 

\smallskip

\noindent  (i) If {\bf (J1)} or {\bf (D1)} holds and $R_D<\infty$, then we have
\begin{align*}
 &p(t, x, y)\simeq  \\
& \begin{cases}
\left(1 \wedge \displaystyle\frac{\delta_*(x,y)}{\p(x,y)^{2}}\right)^{\alpha/2}  \bigg[[\delta_*(x,y)^{\alpha/2} \wedge \phi(t^{-1})^{-1} ]+ F_k^\alpha(d-\alpha n_t ,  T_D,x,y)   \\
\qquad\qquad\qquad \qquad \;\; + \big(n_t t_f-t\big)^{n_t} F_k^\alpha(d- \alpha(n_t-1) , T_D,x,y) \bigg], \qquad\; \mbox{if} \; \; t< \displaystyle \lfloor \frac{d+\alpha}{\alpha} \rfloor t_f; \\
\delta_*(x,y)^{\alpha/2}e^{-ct} , \qquad\qquad \qquad\qquad \qquad\qquad\qquad \qquad\qquad\qquad  \qquad \;   \mbox{if} \;\; t \ge \lfloor\displaystyle\frac{d+\alpha}{\alpha} \rfloor t_f. 
\end{cases}
\end{align*}

\smallskip

\noindent  (ii) If {\bf (J4)} holds and $R_D<\infty$, then we have
\begin{align*}
&p(t, x, y)\simeq  \\
& \begin{cases}
\left(1 \wedge \displaystyle\frac{\delta_*(x,y)}{\p(x,y)^{2}}\right)^{\alpha-1}  \bigg[[\delta_*(x,y)^{\alpha/2} \wedge \phi(t^{-1})^{-1}] +F_c^\alpha(d-\alpha n_t , T_D,x,y)  \\
\qquad\qquad\qquad \qquad \;\; + \big(n_t t_f-t\big)^{n_t} F_c^\alpha(d-  \alpha(n_t-1), T_D,x,y) \bigg],   \qquad  \mbox{if} \; \; t< \lfloor \displaystyle\frac{d+2\alpha-2}{\alpha} \rfloor t_f; \\
\delta_*(x,y)^{\alpha-1}e^{-ct} , \quad \qquad \qquad \qquad\qquad \qquad\qquad\qquad \qquad\qquad   \qquad \quad \mbox{if} \;\; t \ge \lfloor \displaystyle\frac{d+2\alpha-2}{\alpha} \rfloor t_f.
\end{cases}
\end{align*}

\smallskip

\noindent  (iii) If {\bf (J2)} or {\bf (J3)} or {\bf (D2)} or {\bf (D3)} holds, then we have
\begin{align*}
&p(t, x, y)\simeq  \\
& \begin{cases}
\left(1 \wedge \displaystyle\frac{\delta_*(x,y)}{\p(x,y)^{2}}\right)^{\alpha/2}  \bigg[\delta_*(x,y)^{\alpha/2} \wedge \phi(t^{-1})^{-1} + F_k^\alpha(d-\alpha n_t, t,x,y) \\
\qquad\qquad\qquad \qquad \;\;  + \big(n_t t_f-t\big)^{n_t} F_k^\alpha(d-  \alpha(n_t-1), t,x,y) \bigg],\\
\vspace{1mm}
\qquad \qquad \qquad \qquad \qquad \qquad\qquad \quad  \mbox{if} \; \; \p(x,y)^\alpha \le \phi(t^{-1})^{-1} \;\; \text{and} \;\; t< \lfloor (d+\alpha)/\alpha \rfloor t_f; \\
q(ct, x, y) , \;\;\quad\quad\qquad\qquad \qquad \quad\quad  \mbox{if} \;\; \p(x,y)^\alpha > \phi(t^{-1})^{-1} \;\;\; \text{or}\;\;\; t \ge \lfloor (d+\alpha)/\alpha \rfloor t_f.
\end{cases}
\end{align*}

\end{thm}

\smallskip

\begin{remark}
{\rm
When $d>\alpha$, we have that $F^\alpha_k(d, t, x,y) = F^\alpha_c(d,t,x,y) = \p(x,y)^{\alpha-d}$. Thus, by Theorems \ref{t:speciallarge} and \ref{t:specialsub}, under either of the conditions {\bf (L.Poly.)} or {\bf (Sub*.)($\bbeta, \bth$)}, we have that $\lim_{y \to x}p(t,x,y) = \infty$ for all large $t$ even if $D$ is bounded. However, under condition {\bf (Trunc.)($t_f$)}, by Theorem \ref{t:specialtrunc}, we see that $p(t,x,x)<\infty$ for all $t$ large enough. Indeed, we observe that when the kernel $w$ is truncated, the singularity of $p(t,x,y)$ at $x=y$ recedes as the number $\lfloor t/t_f \rfloor$ increases.
}
\end{remark}
\subsection{General results}

In this subsection, we present our estimates for the fundamental solution in full generality. 

Throughout this paper, we always assume that $\{V(x, \cdot) : x \in D\}$ is a family of strictly positive functions satisfying the condition ${\bf WS}(d_1, d_2)$ for some $d_2 \ge d_1>0$ uniformly, that is, there exist constants $c_1, c_2>0$ such that
\begin{align*}
 c_1 \left(\frac{l_2}{l_1}\right)^{d_1} \le \frac{V(x,l_2)}{V(x,l_1)} \le c_2 \left(\frac{l_2}{l_1}\right)^{d_2} \quad \text{for all} \;\; x \in D, \;\; 0<l_1 \le l_2 < \infty.
\end{align*}
We also always assume that $\Phi:[0, \infty) \to [0, \infty)$ is a strictly increasing function such that $\Phi(0)=0$ and satisfies ${\bf WS}(\alpha_1, \alpha_2)$ for some $\alpha_2 \ge \alpha_1>0$.

For a given non-decreaing function $\Psi: (0, \infty) \to [0, \infty)$ such that $\Phi(l) \le \Psi(l)$ for all $l>0$ and satisfies ${\bf WS}(\gamma_1, \gamma_2)$ for some $\gamma_2 \ge \gamma_1>0$, we define
\begin{align*}
q^j(t,x,l;\Phi, \Psi):=& \; \frac{t}{tV(x,\Phi^{-1}(t))+\Psi(l)V(x,l)}.
\end{align*}
Besides, for a given function $\sM: (0, \infty) \times (0, \infty) \to [0, \infty)$ and a constant $a>0$, we define
\begin{align*}
q^d(a, t,x,l; \Phi, \sM):=& \; \frac{\exp\big(-a\sM(t,l)\big)}{V(x, \Phi^{-1}(t))} .
\end{align*}
We will use the functions $q^j$ and $q^d$ to describe interior estimates for $q(t,x,y)$.

On the other hand, for $\gamma \in [0, 1)$ and $(t,x,y) \in (0, \infty) \times D \times D$, we define
\begin{align*}
&a_1^{\gamma}(t,x,y):=\left(\frac{\Phi(\delta_{D}(x))}{\Phi(\delta_{D}(x))+t}\right)^{\gamma}\left(\frac{\Phi(\delta_{D}(y))}{\Phi(\delta_{D}(y))+t}\right)^{\gamma}, \\
&a_2^{\gamma}(t,x,y):=a_1^\gamma(t/(t+1),x,y).
\end{align*}
These functions will be used to describe boundary behaviors of $q(t,x,y)$.

\begin{remark}
{\rm
Observe that for any positive constants $a, b$ and $c$, it holds that $a/(b+c) \le (a/b) \wedge (a/c) \le 2a/(b+c)$.
Hence, we have that
\begin{align*}
&q^j(t,x,l;\Phi, \Psi) \asymp \frac{1}{V(x,\Phi^{-1}(t))} \wedge \frac{t}{\Psi(l)V(x,l)}, \\
&a_1^{\gamma}(t,x,y) \asymp \left(1 \wedge \frac{\Phi(\delta_{D}(x))}{t}\right)^{\gamma}\left(1 \wedge \frac{\Phi(\delta_{D}(y))}{t}\right)^{\gamma}, \\
&a_2^{\gamma}(t,x,y) \asymp \left(1 \wedge \frac{\Phi(\delta_{D}(x))}{t\wedge 1}\right)^{\gamma}\left(1 \wedge \frac{\Phi(\delta_{D}(y))}{t \wedge 1}\right)^{\gamma}.
\end{align*}
}
\end{remark}

We list our candidates for the estimates of the transition density $q(t,x,y)$.

\begin{defn}\label{HKE}
{\rm Let $\gamma \in [0,1)$, $\lambda \in [0, \infty)$ and $k \in \{1, 2\}$.

\smallskip

\noindent (1) We say that $q(t,x,y)$ enjoys the estimate ${\bf HK}_J^{\gamma, \lambda, k}(\Phi, \Psi)$ if
\begin{align*}
\quad q(t,x,y) \asymp  a_1^{\gamma}(t,x,y) q^j(t,x, \p(x,y); \Phi, \Psi) \qquad \text{for all} \;\; (t,x,y) \in (0, 1] \times D \times D,
\end{align*}
and for all $(t,x,y) \in [1, \infty) \times D \times D$,
\begin{align*}
q(t,x,y) \asymp  \begin{cases}
a_k^{\gamma}(t,x,y) q^j(t,x, \p(x,y); \Phi, \Psi),    & \text{if} \quad  \lambda = 0, \\
a_1^{\gamma}(1,x,y) e^{-\lambda t} , &  \text{if} \quad  \lambda > 0.
\end{cases}
\end{align*}

\smallskip

\noindent (2) We say that $q(t,x,y)$ enjoys the estimate ${\bf HK}^{\gamma, \lambda, k}_D(\Phi)$ if $\alpha_1>1$ where $\alpha_1$ is the lower scaling index of $\Phi$, and 
\begin{align*}
\quad q(t,x,y) \simeq  a_1^{\gamma}(t,x,y) q^d(c, t, x, \p(x,y); \Phi, \sM) \qquad \text{for all} \;\; (t,x,y) \in (0, 1] \times D \times D,
\end{align*}
and for all $(t,x,y) \in [1, \infty) \times D \times D$,
\begin{align*}
q(t,x,y) \simeq  \begin{cases}
a_k^{\gamma}(t,x,y) q^d(c, t, x, \p(x,y); \Phi, \sM),    & \text{if} \quad  \lambda = 0, \\
a_1^{\gamma}(1,x,y) e^{-\lambda t}, &  \text{if} \quad  \lambda > 0,
\end{cases}
\end{align*}
where the function $\sM(t,l)$ is a strictly positive for all $t,l>0$, non-increasing on $(0,\infty)$ for each fixed $l>0$ and determined by the following relation
\begin{align}\label{e:sM}
\frac{t}{\sM(t,l)} \asymp \Phi \left(\frac{l}{\sM(t,l)}\right) \qquad \text{for all} \;\; t,l>0.
\end{align}

\smallskip

\noindent (3) We say that $q(t,x,y)$ enjoys the estimate ${\bf HK}^{\gamma, \lambda, k}_M(\Phi, \Psi)$ if $\alpha_1>1$ where $\alpha_1$ is the lower scaling index of $\Phi$, and there are functions $q^j, q^d$ such that  
$$
q(t,x,y) = q^j(t,x,y) + q^d(t,x,y) \qquad \text{for all} \quad (t,x,y) \in (0, \infty) \times D \times D,
$$
and $q^j$ and $q^d$ enjoy the estimate ${\bf HK}_J^{\gamma, \lambda, k}(\Phi, \Psi)$ and ${\bf HK}_D^{\gamma, \lambda, k}(\Phi)$, respectively.
}
\end{defn}

\vspace{1mm}

In the rest of this subsection, we always assume that $q(t,x,y)$ enjoys one of the estimates ${\bf HK}^{\gamma, \lambda, k}_J(\Phi, \Phi), {\bf HK}^{\gamma, \lambda, k}_D(\Phi)$ and ${\bf HK}^{\gamma, \lambda, k}_M(\Phi,\Psi)$ for some $\gamma \in [0,1)$, $\lambda \ge 0$ and $k \in \{1,2\}$. If $\lambda>0$, then we further assume that $D$ is bounded so that $R_D=\text{diam}(D)<\infty$.

\vspace{2mm}

\begin{example}\label{e:example}
{\rm

(1) Examples of estimates ${\bf HK}_J^{\gamma, \lambda, k}(\Phi, \Psi)$, ${\bf HK}_D^{\gamma, \lambda, k}(\Phi)$ and ${\bf HK}^{\gamma, \lambda, k}_M(\Phi, \Psi)$ include all estimates given in subsection \ref{s:special}. For example, we see that estimate {\bf (J1)} is nothing but estimate ${\bf HK}_J^{1/2, \lambda, 1}(\Phi_\alpha, \Phi_\alpha)$ for $\lambda>0$ where $\Phi_\alpha(x):=x^\alpha$.

(2) The factor $e^{-\lambda t}a_1^\gamma(1,x,y)$ usually appears in the global estimates of the Dririchlet heat kernel when $D$ is a $C^{1,1}$ bounded open set, $a_1^{\gamma}(t,x,y)$ appears when $D$ is a half space-like $C^{1,1}$ open set and $a_2^{\gamma}(t,x,y)$ appears when $D$ is a exterior of a bounded $C^{1,1}$ open set. Various examples are given in \cite{BGT, CK, CKS3, CKS4, KM1, So, Zh}.

(3) Recently, in \cite{CKSV}, we, jointly with Renming Song and Zoran Vondraček give examples of generators whose transition density satisfies estimate ${\bf HK}_J^{\gamma, \lambda, 1}(\Phi_\alpha, \Phi_\alpha)$ for each $0<\alpha<2$ and $\gamma \in [0 \vee (\alpha-1)/\alpha, 1)$.

(4) Examples of symmetric Markov processes (including non L\'evy processes) satisfying the mixed heat kernel estimates ${\bf HK}^{\gamma,\lambda,k}_M(\Phi, \Psi)$ can be found in \cite{BKKL, BKKL2, KM1, Mi}. We will show that one of the explict expressions of the function $\sM$ is given by
\begin{align*}
\sM(t, l) := \sup_{s>0} \left\{\frac{l}{s}-\frac{t}{\Phi(s)}\right\},
\end{align*}
which appears in the exponential terms in \cite{BKKL2}. (See, Lemma \ref{l:sM}(i).)
}
\end{example}

\vspace{1mm}

We introduce some functions which will be used in near diagonal estimates for the fundamental solution. Define for $(t,x,y) \in (0, \infty) \times D \times D$, $\gamma \in [0,1)$ and $k \in \{1,2\}$,
\begin{align}\label{e:defD}
&\sI_k^\gamma(t,x,y):=\int_{\Phi(\p(x,y))}^{1/(2e^2 \phi(t^{-1}))}  \frac{a_k^{\gamma}(r,x,y)}{V(x, \Phi^{-1}(r))} dr, \nn\\
&\sJ_k^\gamma(t,x,y):= \frac{a_k^{\gamma}(1/\phi(t^{-1}),x,y)}{V\big(x, \Phi^{-1}(1/\phi(t^{-1}))\big)} + w(t) \sI_k^\gamma(t,x,y).
\end{align}

Under certain weak scaling conditions for $V$ and $\Phi$, we can calculate the integral term $\sI^\gamma_k$ explicitly. (See, Proposition \ref{p:Dgamma}.)
Now, we are ready to state the main results.

\begin{thm}\label{t:mainsmall}
Let $p(t,x,y)$ be given by \eqref{e:dfundamental}. Assume that $w$ satisfies conditions {\bf (Ker.)} and
{\bf (S.Poly.)($t_s$)}. Then the follwing estimates hold for all $(t,x,y) \in (0, t_s] \times D \times D$.

\smallskip

\noindent (i) (Near diagonal estimates) If $\Phi(\p(x,y))\phi(t^{-1}) \le 1/(4e^2)$, then we have
\begin{align*}
p(t,x,y) &\asymp \sJ_k^{\gamma}(t, x, y).
\end{align*}

\noindent (ii) (Off diagonal estimates) Suppose that $\Phi(\p(x,y))\phi(t^{-1}) > 1/(4e^2)$.

\indent (a) If $q(t,x,y)$ enjoys the estimate ${\bf HK}^{\gamma, \lambda, k}_J(\Phi, \Phi)$, then we have
\begin{align*}
p(t,x,y) \asymp  \frac{a_k^{\gamma}(1/\phi(t^{-1}),x,y)}{\phi(t^{-1})\Phi(\p(x,y))V(x,\p(x,y))}.
\end{align*}

\smallskip

\indent (b) If $q(t,x,y)$ enjoys the estimate ${\bf HK}^{\gamma, \lambda, k}_D(\Phi)$, then we have
\begin{align*}
p(t,x,y) \simeq a_k^{\gamma}(1/\phi(t^{-1}),x,y) \frac{\exp\big(-c\sN(t,\p(x,y))\big)}{V\big(x,\Phi^{-1}(1/\phi(t^{-1}))\big)} ,
\end{align*}
where $\sN(\cdot, l)$ is a strictly positive function which is determined by the following relation
\begin{align}\label{e:defN}
\frac{1}{\phi\big(\sN(t,l)/t\big)} \asymp \Phi\left(\frac{l}{\sN(t,l)}\right), \qquad t,l>0.
\end{align}

\smallskip

\indent (c) If $q(t,x,y)$ enjoys the estimate ${\bf HK}^{\gamma, \lambda, k}_M(\Phi,\Psi)$, then we have 
\begin{align*}
 p(t,x,y) \simeq a_k^{\gamma}(1/\phi(t^{-1}),x,y) \left(\frac{1}{\phi(t^{-1})\Psi(\p(x,y))V(x,\p(x,y))} + \frac{\exp\big(-c \sN(t,\p(x,y))\big)}{V\big(x,\Phi^{-1}(1/\phi(t^{-1}))\big)}\right).
\end{align*}
\end{thm}

\smallskip

Recall that $R_D=\text{diam}(D)$ and $T_D=[\phi^{-1}(4^{-1}e^{-2}R_D^{-\alpha})]^{-1}$.

\begin{thm}\label{t:mainlarge}
Let $p(t,x,y)$ be given by \eqref{e:dfundamental}. Assume that $w$ satisfies conditions {\bf (Ker.)} and
{\bf (L.Poly.)}. Then for every fixed $T>0$, the following estimates hold for all $(t,x,y) \in [T, \infty) \times D \times D$.

\smallskip
\noindent (i) If $\lambda=0$, then estimates given in Theorem \ref{t:mainsmall} hold for all $(t,x,y) \in [T, \infty) \times D \times D$.

\noindent (ii) If $\lambda>0$ and $R_D<\infty$, then we have
\begin{align*}
p(t,x,y) \asymp w(t) F_1^\gamma(T_D,x,y) = w(t) \int_{\Phi(\p(x,y))}^{2\Phi(R_D)} \frac{a_1^{\gamma}(r,x,y)}{V(x, \Phi^{-1}(r))} dr.
\end{align*}
\end{thm}

\smallskip

\begin{remark}
{\rm

(1) By Lemma \ref{l:sM}(i), one of the explict expressions of the function $\sN$ satisfying \eqref{e:defN} is given by
\begin{align*}
\sN(t,l) := \sup_{s>0} \left\{ \frac{l}{s} - t\phi^{-1}(1/\Phi(s))\right\}.
\end{align*}

(2) Theorems \ref{t:mainsmall} and \ref{t:mainlarge} recover \cite[Theorems 1.6 and 1.8]{CKKW}. Indeed, the assumptions in \cite{CKKW} can be interpreted as the kernel $w$ satisfies conditions {\bf (Ker.)}, {\bf (S.Poly.)($t_s$)} and ${\bf (L.Poly.)}$ for some $0<\delta_1, \delta_2<1$ and $q(t,x,y)$ enjoys either of the estimates ${\bf HK}^{0,0,1}_J(\Phi, \Phi)$ or ${\bf HK}^{0, 0,1}_D(\Phi)$.

(3) In off diagonal situations, that is, when $\Phi(\p(x,y)) \ge \phi(t^{-1})^{-1}$, estimates for $p(t,x,y)$ can be factorized into the boundary factors and the rest.
However, there is no such factorization on near diagonal situation in general since $\sJ^\gamma_k(t,x,y)$ can not be factorized commonly. (cf. Theorem \ref{t:specialsmall}.)
}
\end{remark}

When condition {\bf (Sub.)($\bbeta, \bth$)} holds, the bounds for fundamental solution decrease subexponentially as $t \to \infty$. Moreover, when $0<\beta<1$ and $D$ is bounded, we obtain the sharp upper bounds that decrease with exactly the same rate as the upper bound for $w$ as $t \to \infty$.

\begin{thm}\label{t:mainsub}
Let $p(t,x,y)$ be given by \eqref{e:dfundamental}. Assume that $w$ satisfies conditions {\bf (Ker.)} and
{\bf (Sub.)($\bbeta, \bth$)}. Then for every fixed $T>0$, the following estimates hold for all $(t,x,y) \in [T, \infty) \times D \times D$.

\smallskip

\noindent (i) Suppose that $\lambda = 0$.

\indent  (a) If $\Phi(\p(x,y))\phi(t^{-1}) \le 1/(4e^2)$, then there exists a constant $c>1$ such that
\begin{align*}
&c^{-1}\left(\frac{a_k^{\gamma}(t,x,y)}{V\big(x, \Phi^{-1}(t)\big)} + w(t) \int_{\Phi(\p(x,y))}^{1/(2e^2 \phi(t^{-1}))}  \frac{a_k^{\gamma}(r,x,y)}{V(x, \Phi^{-1}(r))} dr \right)\\
&\quad \le p(t,x,y) \le c\left(\frac{a_k^{\gamma}(t,x,y)}{V\big(x, \Phi^{-1}(t)\big)} + \exp \big(-\frac{\theta}{2} t^{\beta} \big) \int_{\Phi(\p(x,y))}^{1/(2e^2 \phi(t^{-1}))}  \frac{a_k^{\gamma}(r,x,y)}{V(x, \Phi^{-1}(r))} dr \right),
\end{align*}
where $\theta>0$ is the constant in condition {\bf (Sub.)($\bbeta, \bth$)}.

\indent (b) If $\Phi(\p(x,y))\phi(t^{-1}) > 1/(4e^2)$, then we have
\begin{align*}
p(t,x,y) \simeq q(ct,x,y).
\end{align*}

\noindent  (ii) Suppose that $\lambda>0$ and $R_D<\infty$. Then, there exist constants $L_1, L_2>0$ independent of $\lambda$ and $c>1$ such that in the case when $\beta \in (0,1)$, we have
\begin{align*}
&c^{-1}w(t)\int_{\Phi(\p(x,y))}^{2\Phi(R_D)} \frac{a_1^{\gamma}(r,x,y)}{V(x, \Phi^{-1}(r))} dr 
\le  p(t,x,y) \le c \exp\big(-\theta t^{\beta}\big) \int_{\Phi(\p(x,y))}^{2\Phi(R_D)} \frac{a_1^{\gamma}(r,x,y)}{V(x, \Phi^{-1}(r))} dr,
\end{align*}
and in the case when $\beta =1$, we have
\begin{align*}
&c^{-1}\left(w(t)\int_{\Phi(\p(x,y))}^{2\Phi(R_D)} \frac{a_1^{\gamma}(r,x,y)}{V(x, \Phi^{-1}(r))} dr + e^{-\lambda L_1t}\Phi(\delta_D(x))^\gamma\Phi(\delta_D(y))^\gamma\right) \\
&\quad \le  p(t,x,y) \le c \left(\exp\big(-\frac{\theta}{2} t\big) \int_{\Phi(\p(x,y))}^{2\Phi(R_D)} \frac{a_1^{\gamma}(r,x,y)}{V(x, \Phi^{-1}(r))} dr + e^{-\lambda L_2t}\Phi(\delta_D(x))^\gamma\Phi(\delta_D(y))^\gamma\right),
\end{align*}
where $\theta>0$ is the constant in condition {\bf (Sub.)($\bbeta, \bth$)}.
\end{thm}

Our last theorem gives the estimates for $p(t,x,y)$ when $w$ is finitely supported.

\begin{thm}\label{t:main2}
Let $p(t,x,y)$ be given by \eqref{e:dfundamental}. Assume that $w$ satisfies conditions {\bf (Ker.)} and {\bf (Trunc.)($t_f$)}. Then the follwing estimates hold for all $(t,x,y) \in [t_f/2, \infty) \times D \times D$. Let $n_t:= \lfloor t/t_f \rfloor +1 \in \N$. 

\smallskip

\noindent (i) Suppose that $\lambda = 0$.

\smallskip

\indent (a) If $\Phi(\p(x,y)) \le t \le \lfloor d_2/\alpha_1 + 2 \gamma \rfloor t_f$, then
\begin{align*}
p(t,x,y) &\asymp \int_{\Phi(\p(x,y))}^{2t}  \frac{r^{n_t}a_1^{\gamma}(r,x,y)}{V(x, \Phi^{-1}(r))} dr + (n_t t_f -t)^{n_t} \int_{\Phi(\p(x,y))}^{2t}  \frac{r^{n_t-1}a_1^{\gamma}(r,x,y)}{V(x, \Phi^{-1}(r))} dr,
\end{align*}

\indent (b) If $\Phi(\p(x,y)) \le t$ and $t \ge \lfloor d_2/\alpha_1 + 2 \gamma \rfloor t_f$, then
\begin{align*}
p(t,x,y) \asymp \frac{a_k^\gamma(t,x,y)}{V\big(x, \Phi^{-1}(t)\big)} \asymp  q(t,x,y).
\end{align*}

\indent  (c) If $\Phi(\p(x,y))>t$, then
\begin{align*}
p(t,x,y) \simeq q(ct,x,y).
\end{align*}

\noindent (ii) Suppose that $\lambda>0$ and $R_D<\infty$.

\smallskip

\indent  (a) If $t \le \lfloor d_2/\alpha_1 + 2 \gamma \rfloor t_f$, then
\begin{align*}
p(t,x,y) &\simeq \int_{\Phi(\p(x,y))}^{2\Phi(R_D)}  \frac{r^{n_t}a_1^{\gamma}(r,x,y)}{V(x, \Phi^{-1}(r))} dr + (n_t t_f -t)^{n_t} \int_{\Phi(\p(x,y))}^{2\Phi(R_D)}  \frac{r^{n_t-1}a_1^{\gamma}(r,x,y)}{V(x, \Phi^{-1}(r))} dr,
\end{align*}

\indent  (b) If $t \ge \lfloor d_2/\alpha_1 + 2 \gamma \rfloor t_f$, then
\begin{align*}
p(t,x,y) \simeq e^{-ct}\Phi(\delta_D(x))^\gamma\Phi(\delta_D(y))^\gamma \simeq q(t,x,y).
\end{align*}

\end{thm}

\begin{remark}
{\rm
Note that under settings of Theorem \ref{t:main2}, we can apply Theorem \ref{t:mainsmall} to obtain the estimates of $p(t,x,y)$ for all $(t,x,y) \in (0, t_f/2] \times D \times D$. Hence, we have the global estimates for $p(t,x,y)$ under those settings.
}
\end{remark}

\smallskip

As a consequence of the estimates for the fundamental solution, we have that the solution to the Dirichlet problem \eqref{e:dtfe} vanishes continuously on the boundary of $D$. Indeed, under mild conditions, the solution $u(t,x)$ vanishes exactly the same rate as a transition density $q(t,x,y)$.

\smallskip

{\bf (V.)} There exists a constant $c_V>1$ such that for all $x \in D$ and $0<l \le R_D=\text{diam}(D)$,
\begin{align*}
c_V^{-1}V(x,l) \le m\big( \{y \in D : \p(x,y) \le l\}\big) \le c_V V(x,l).
\end{align*}

\begin{cor}\label{c:boundary}
Suppose that $(D,\p, m)$ satisfies {\bf (V.)} and $w$ satisfies conditions {\bf (Ker.)}, {\bf (S.Poly.)($t_s$)} and one among {\bf (L.Poly.)}, {\bf (Sub.)($\bbeta, \bth$)} and {\bf (Trunc.)($t_f$)}. We also assume that $q(t,x,y)$ enjoys one of the estimates ${\bf HK}^{\gamma, \lambda, k}_J(\Phi, \Phi), {\bf HK}^{\gamma, \lambda, k}_D(\Phi)$ and ${\bf HK}^{\gamma, \lambda, k}_M(\Phi,\Psi)$ for some $0<\gamma<1$, $\lambda \ge 0$ and $k \in \{1,2\}$. When $\lambda>0$, we further assume that $D$ is bounded. Then, for all bounded measurable function $f$ on $D$, $u(t,x):=\E[T^D_{E_t}f(x)]$ satisfies the following boundary condition:

\smallskip

\indent \indent For any fixed $t>0$, there exists a constant $c_1>0$ such that for every $x \in D$,
\begin{align*}
|u(t,x)| \le c_1 \lVert f \rVert_{\infty} \Phi(\delta_D(x))^\gamma.
\end{align*}
\end{cor}
\proof
Since the main ideas are similar, we only give the proof for the case when $w$ satisfies {\bf (Ker.)}, {\bf (S.Poly.)($t_s$)} and {\bf (L.Poly.)} and $q(t,x,y)$ enjoys estimate ${\bf HK}^{\gamma, \lambda, k}_J(\Phi, \Phi)$ for some $\gamma \in (0,1)$, $\lambda = 0$ and $k \in \{1,2\}$. Fix $t>0$ and we let $A_t:=\Phi^{-1}\big(1/(4e^2\phi(t^{-1}))\big)$. By Theorems \ref{t:mainsmall} and \ref{t:mainlarge}, for every $x \in D$,
\begin{align*}
&|u(t,x)| = \left| \int_D p(t,x,y)f(y) m(dy) \right| \\
&\le c\lVert f \rVert_{\infty} \left( \int_{\{y \in D:\p(x,y) \le A_t\}} \sJ^\gamma_k(t,x,y) m(dy) +  \int_{\{y \in D:\p(x,y) > A_t\}} \frac{\Phi(\delta_D(x))^\gamma}{\Phi(\p(x,y))V(x,\p(x,y))} m(dy) \right) \\
&=: c\lVert f \rVert_{\infty}  (I_1 + I_2).
\end{align*}

Set $\eta =  \frac{d_1}{2\alpha_2} \wedge \frac{1-\gamma}{2}$. Since $\eta<\frac{d_1}{\alpha_2}$, by \cite[Theorem 2.2.2]{BGT}, we have that for all $x \in D$ and $0<s<t$,
\begin{align}\label{almost}
\inf_{r \in (s,t]} r^{-\eta}V(x, \Phi^{-1}(r)) \asymp s^{-\eta}V(x, \Phi^{-1}(s)).
\end{align}

Then, by Fubini's theorem, \eqref{almost}, condition {\bf (V.)} and the weak scaling properties of $V$ and $\Phi$,
\begin{align*}
\frac{I_1}{\Phi(\delta_D(x))^\gamma}& \le c \sum_{k=1}^{\infty}\int_{\{y \in D:2^{-k}A_t < \p(x,y) \le 2^{-(k-1)}A_t \}} \int_{\Phi(\p(x,y))}^{1/(2e^2 \phi(t^{-1}))} \frac{1}{r^\gamma V(x,\Phi^{-1}(r))} dr m(dy) \\
&\le c \sum_{k=1}^{\infty}\int_{\{y \in D:2^{-k}A_t < \p(x,y) \le 2^{-(k-1)}A_t \}}m(dy) \int_{\Phi(2^{-k}A_t)}^{1/(2e^2 \phi(t^{-1}))} \frac{1}{r^{\gamma+\eta} r^{-\eta}V(x,\Phi^{-1}(r))} dr \\
&\le c\sum_{k=1}^{\infty} \frac{V(x, 2^{-(k-1)}A_t)}{ \Phi(2^{-k}A_t)^{-\eta}V(x,2^{-k}A_t)}\int_{\Phi(2^{-k}A_t)}^{1/(2e^2 \phi(t^{-1}))} \frac{1}{r^{\gamma+\eta}} dr \\
&\le c \Phi(A_t)^\eta \left( \int_0^{1/(2e^2 \phi(t^{-1}))} \frac{1}{r^{\gamma+\eta}} dr\right) \sum_{k=1}^{\infty} 2^{-k \eta \alpha_1}   \le  c.
\end{align*}

Moreover, we also have that by condition {\bf (V.)} and the weak scaling properties of $V$ and $\Phi$,
\begin{align*}
\frac{I_2}{\Phi(\delta_D(x))^\gamma} &\le  c \sum_{k=1}^{\infty}\int_{\{y \in D:2^{k-1}A_t < \p(x,y) \le 2^{k}A_t\}} \frac{1}{\Phi(\p(x,y))V(x,\p(x,y))}m(dy) \\
&\le  c \sum_{k=1}^{\infty} \frac{V(x, 2^kA_t)}{\Phi(2^{k-1}A_t)V(x,2^{k-1}A_t)}\le c \sum_{k=1}^{\infty} \frac{2^{-k \alpha_1}}{\Phi(A_t)} \le  c.
\end{align*}

Therefore, we get the result.

\qed

In the end of this section, we study explicit forms of $\sJ^\gamma_k(t, x, y)$ $(0\le \gamma<1)$ under some weak scaling conditions for $V$ and $\Phi$.
Recall that $\Phi(\cdot)$ satisfies ${\bf WS}(\alpha_1, \alpha_2)$ and $V(x, \cdot)$ satisfies ${\bf WS}(d_1, d_2)$ uniformly. We define $\delta_*^\Phi(x,y):=\Phi(\delta_{D}(x))\Phi(\delta_{D}(y))$.

\begin{prop}\label{p:Dgamma}
Let $\gamma \in [0,1)$. If $\gamma = 0$, then we redefine $\delta_D(x) = \infty$ for all $x \in D$.
Then, the following estimates hold for all $(t,x,y) \in (0, \infty) \times D \times D$ satisfying $\Phi(\p(x,y))\phi(t^{-1}) \le 1/(4e^2)$.

\smallskip

\indent (a) If $d_2/\alpha_1<1-2\gamma$, then
\begin{align*}
\sJ^{\gamma}_1(t,x,y) &\asymp \left(1 \wedge \frac{\delta_*^\Phi(x,y)}{\phi(t^{-1})^{-2}}\right)^\gamma \frac{1}{V\big(x, \Phi^{-1}(1/\phi(t^{-1}))\big)}.
\end{align*}

(b) If $\alpha_1=\alpha_2$, $d_1=d_2=(1-2\gamma)\alpha_1$ and $\gamma>0$, then
\begin{align*}
\sJ^{\gamma}_1(t,x,y) &\asymp \left(1 \wedge \frac{\delta_*(x,y)^{\alpha_1}}{\phi(t^{-1})^{-2}}\right)^\gamma \phi(t^{-1})^{1-2\gamma} +  \1_{\{\delta_*(x,y)^{\alpha_1/2} \le \phi(t^{-1})^{-1}\}}w(t)\\
&\qquad \times \delta_*(x,y)^{\alpha_1 \gamma}  \log^+ \left(\frac{2\phi(t^{-1})^{-1}}{\big(\p(x,y) \vee \delta_\vee(x,y) \big)^{\alpha_1}}\right).
\end{align*}

\indent (c) If $1-2\gamma < d_1/\alpha_2 \le d_2/\alpha_1 <1 - \gamma$, then
\begin{align*}
\sJ^{\gamma}_1(t,x,y) \asymp& \left(1 \wedge \frac{\delta_*^\Phi(x,y)}{\phi(t^{-1})^{-2}}\right)^\gamma \frac{1}{V\big(x, \Phi^{-1}(1/\phi(t^{-1}))\big)} +  \1_{\{\delta_*^\Phi(x,y)^{1/2} \le \phi(t^{-1})^{-1}\}}w(t)\\
&\quad \times  \left(1 \wedge \frac{\delta_*^\Phi(x,y)}{\Phi(\p(x,y))^2}\right)^\gamma \left( \frac{\Phi(\p(x,y))}{V(x,\p(x,y))} \vee \frac{\delta_*^\Phi(x,y)^\gamma \Phi(\delta_\vee(x,y))^{1-2\gamma}}{V(x,\delta_\vee(x,y))}\right).
\end{align*}

\indent (d) If $\alpha_1 = \alpha_2$, $d_1=d_2=(1-\gamma)\alpha_1$ and $\gamma>0$, then
\begin{align*}
\sJ^{\gamma}_1(t,x,y) &\asymp \left(1 \wedge \frac{\delta_*(x,y)^{\alpha_1}}{\phi(t^{-1})^{-2}}\right)^\gamma \phi(t^{-1})^{1-\gamma}+  \1_{\{\delta_*(x,y)^{\alpha_1/2} \le \phi(t^{-1})^{-1}\}}w(t)\\
&\quad \times \left( 1 \wedge \frac{\delta_*(x,y)}{\p(x,y)^2}\right)^{\alpha_1 \gamma} \left( \p(x,y)^{\alpha_1 \gamma} + \delta_\wedge(x,y)^{\alpha_1 \gamma} \log\left(\frac{\p(x,y) \vee 2\delta_\vee(x,y)}{\p(x,y) \vee \delta_\wedge(x,y)}\right)\right).
\end{align*}

\indent (e)  If $1-\gamma < d_1/\alpha_2 \le d_2/\alpha_1 <1$, then
\begin{align*}
\sJ^{\gamma}_1(t,x,y) &\asymp \left(1 \wedge \frac{\delta_*^\Phi(x,y)}{\phi(t^{-1})^{-2}}\right)^\gamma \frac{1}{V\big(x, \Phi^{-1}(1/\phi(t^{-1}))\big)} +  \1_{\{\delta_*^\Phi(x,y)^{1/2} \le \phi(t^{-1})^{-1}\}}w(t)\\
&\quad \times  \left(1 \wedge \frac{\delta_*^\Phi(x,y)}{\Phi(\p(x,y))^2}\right)^\gamma \left( \frac{\Phi(\p(x,y))}{V(x,\p(x,y))} \vee \frac{\Phi(\delta_\wedge(x,y))}{V(x,\delta_\wedge(x,y))}\right).
\end{align*}

\indent (f) If $\alpha_1 = \alpha_2 = d_1=d_2$, then
\begin{align*}
\sJ^{\gamma}_1(t,x,y) &\asymp \left(1 \wedge \frac{\delta_*(x,y)^{\alpha_1}}{\phi(t^{-1})^{-2}}\right)^\gamma \phi(t^{-1}) +  w(t)\\
& \quad \times \left( 1 \wedge \frac{\delta_*(x,y)}{\p(x,y)^2}\right)^{\alpha_1 \gamma} \left( 1 + \log^+ \left(\frac{2\phi(t^{-1})^{-1} \wedge 2\delta_\wedge(x,y)^{\alpha_1}}{\p(x,y)^{\alpha_1}}\right)\right).
\end{align*}

\indent (g)  If $1<d_1/\alpha_2$, then
\begin{align*}
\sJ^{\gamma}_1(t,x,y) \asymp& \left(1 \wedge \frac{\delta_*^\Phi(x,y)}{\phi(t^{-1})^{-2}}\right)^\gamma \frac{1}{V\big(x, \Phi^{-1}(1/\phi(t^{-1}))\big)} +  w(t) \left(1 \wedge \frac{\delta_*^\Phi(x,y)}{\Phi(\p(x,y))^2}\right)^\gamma  \frac{\Phi(\p(x,y))}{V(x,\p(x,y))}.
\end{align*}

\end{prop}
\proof
See Appendix.
\qed

\begin{remark}\label{r:sJ2}
{\rm
We can obtain closed forms of $\sJ^\gamma_2$ from closed forms of $\sJ^\gamma_1$ and $\sJ^0_1$. Indeed, for every fixed $T>0$, we can check that for all $\gamma \in [0,1)$ and $(t,x,y) \in (0,T] \times D \times D$, it holds that $\sJ^\gamma_2(t,x,y) \asymp \sJ^\gamma_1(t,x,y)$. 

Moreover, observe that for all large $t$ such that $\Phi(1)\phi(t^{-1}) \le 1/(8e^2)$, we have
\begin{align*}
&\int_{\Phi(\p(x,y))}^{1/(2e^2\phi(t^{-1}))} \frac{a^\gamma_2(r,x,y)}{V(x, \Phi^{-1}(r))} dr\\
& \asymp  a_1^\gamma(1,x,y)\int_{2\Phi(1) \vee \Phi(\p(x,y)) }^{1/(2e^2\phi(t^{-1}))} \frac{1}{V(x, \Phi^{-1}(r))} dr + \1_{\{\p(x,y) \le 1\}}\int_{\Phi(\p(x,y))}^{2\Phi(1)} \frac{a^\gamma_1(r,x,y)}{V(x, \Phi^{-1}(r))} dr.
\end{align*}

\noindent Add an isolated point $y_0$ to $D$ and define $\p(x,y_0) = 1$ for all $x \in D$. 
By the above observation, we have that for any fixed $T>0$, the following comparison holds for all $\gamma \in [0,1)$ and $(t,x,y) \in [T,\infty) \times D \times D$:
\begin{align*}
\sJ^\gamma_2(t,x,y) \asymp  \big(1 \wedge \Phi(\delta_D(x))\big)^\gamma \big(1 \wedge \Phi(\delta_D(y))\big)^\gamma\sJ^0_1(t,x,y') + \1_{\{\p(x,y) \le 1\}}\sJ^\gamma_1([\phi^{-1}(4^{-1}e^{-2})]^{-1},x,y),
\end{align*}
where $y'=y$ if $\p(x,y) \ge 1$ and $y'=y_0$ if $\p(x,y) < 1$. (cf. Theorem \ref{t:speciallarge}(v)(a).)
}
\end{remark}

\section{Estimates for Subordinator}\label{s:subordinator}

Throughout this section, we always assume that $S$ be the subordinator whose Laplace exponent has the following representation with a function $w$ satisfying condition {\bf (Ker.)}: 
$$
\phi(\lambda) = -\log \E\big[\exp(-\lambda S_1)\big] = \int_0^{\infty} (1-e^{-\lambda s}) (-dw(s)) \quad \text{for all} \;\; \lambda \ge 0.
$$ 

Following \cite{JP}, we let 
$$H(\lambda):=\phi(\lambda)-\lambda\phi'(\lambda) \quad \text{for all} \;\; \lambda \ge 0.
$$
In \cite{JP}, Naresh C. Jain and William E. Pruitt studied asymptotic properties of lower tail probabilities of subordinators, $\P(S_r \le t)$, in terms of the function $H$. Then, in \cite{Mi}, Ante Mimica obtained esitmates for upper tail probabilities, $\P(S_r \ge t)$, in terms of the function $H$ as well. Those estimates were crucial ingredients in \cite{CKKW} to establish the estimates for the fundamental solution $p(t,x,y)$. 

In this section, we will improve the results in \cite{Mi} and obtain tail probability estimates in terms of the tail measure $w$ instead of the function $H$. This allows us to get estimates for the fundamental solution in more general situations. 

\subsection{General estimates for subordinator.}

\begin{lemma}\label{l:phiHw}
(i) For every $\lambda>0$, we have
\begin{align*}
\phi(\lambda) \asymp \lambda \int_0^{1/\lambda} w(s)ds  \quad \text{and} \quad H(\lambda) \asymp \lambda^2 \int_0^{1/\lambda} sw(s)ds.
\end{align*}
(ii) If $w$ satisfies ${\bf WS^0}(-\alpha_1, -\alpha_2, c_0)$ {\rm (resp. ${\bf WS^\infty}(-\alpha_1, -\alpha_2, c_0)$)} for some constants $\alpha_1 \ge \alpha_2 \ge 0$ and $c_0>0$, then $\phi$ satisfies ${\bf WS^{\infty}}(\alpha_2, \alpha_1 \wedge 1, 1/c_0)$ {\rm (resp. ${\bf WS^0}(\alpha_2, \alpha_1 \wedge 1, 1/c_0)$)} and $H$ satisfies ${\bf WS^{\infty}}(\alpha_2, \alpha_1 \wedge 2, 1/c_0)$. {\rm (resp. ${\bf WS^0}(\alpha_2, \alpha_1 \wedge 2, 1/c_0)$.)} 

In particular, if there exist constants $\alpha_1<2$ and $c_0>0$ such that $w$ satisfies ${\bf LS^0}(-\alpha_1, c_0)$ {\rm (resp. ${\bf LS^\infty}(-\alpha_1, c_0)$)}, then we have
\begin{align*}
w(s) \asymp H(s^{-1}), \qquad \text{for all} \;\; 0<s \le c_0. \quad \text{\rm (resp. for all} \;\; s \ge c_0.)
\end{align*}

Conversely, if $\phi$ satisfies ${\bf WS^{\infty}}(\alpha_2, \alpha_1, c_0)$ {\rm (resp. ${\bf WS^0}(\alpha_2, \alpha_1, c_0)$)} for some constants $0 \le \alpha_2 \le \alpha_1 < 1$ and $c_0>0$ or $H$ satisfies ${\bf WS^{\infty}}(\alpha_2, \alpha_1, c_0)$ {\rm (resp. ${\bf WS^0}(\alpha_2, \alpha_1, c_0)$)} for some constants $0 \le \alpha_1 \le \alpha_2 < 2$ and $c_0>0$, then there exists a constant $c_1>0$ such that $w$ satisfies ${\bf WS^0}(-\alpha_1, -\alpha_2, c_1)$. {\rm (resp. ${\bf WS^\infty}(-\alpha_1, -\alpha_2, c_1)$.)}

\end{lemma}
\proof
(i) By the integration by parts and Fubini's theorem,
\begin{align*}
\phi(\lambda) = \lambda \int_0^{\infty} \int_0^s e^{-\lambda u}du (-dw(s)) = \lambda \int_0^{1/\lambda} e^{-\lambda u} w(u) du + \lambda \int_{1/\lambda}^{\infty} e^{-\lambda u} w(u) du =: I_1 + I_2.
\end{align*}
First, we see that $I_1 \asymp \lambda \int_0^{1/\lambda} w(s)ds$. Moreover, since $w$ is non-increasing,
\begin{align*}
I_2 \le w(1/\lambda) \int_{1/\lambda}^{\infty} \lambda e^{-\lambda u} du \le \frac{1}{2}w(1/\lambda) \le \lambda \int_{1/(2\lambda)}^{1/\lambda} w(s) ds \le \lambda \int_0^{1/\lambda} w(s)ds.
\end{align*}
Hence, the first claim holds. On the other hand, note that by the definition of $H$,
\begin{align*}
H(\lambda) = -\lambda^2 (\lambda^{-1}\phi(\lambda))' = \lambda^2 \int_0^{\infty} ue^{-\lambda u}w(u)du.
\end{align*}
Then, we can deduce that $H(\lambda) \asymp \lambda^2 \int_0^{1/\lambda} sw(s)ds$ by a similar argument.

(ii) We first assume that $w$ satisfies ${\bf WS^0}(-\alpha_1, -\alpha_2, c_0)$  for some constants $\alpha_1 \ge \alpha_2 \ge 0$ and $c_0>0$. By (i),
\begin{align}\label{phiscale}
\phi(\kappa \lambda) \asymp \kappa \lambda \int_0^{1/(\kappa \lambda)} w(s) ds =\lambda \int_0^{1/\lambda} \frac{w(s/\kappa)}{w(s)} w(s)ds \qquad \text{for all} \;\; \kappa \ge 1, \; \lambda \ge 1/c_0.
\end{align}
Moreover, by the assumption, there are constants $c_2, c_3>0$ such that 
\begin{align*}
 c_2 \kappa^{\alpha_2} \le \frac{w(s/\kappa)}{w(s)} \le c_3 \kappa^{\alpha_1} \qquad \text{for all} \;\; \kappa \ge 1, \; \lambda \ge 1/c_0.
\end{align*}
Thus, we deduce that $\phi$ satisfies ${\bf WS^{\infty}}(\alpha_2, \alpha_1, 1/c_0)$ from \eqref{phiscale} and (i). Since $\phi$ always satisfy ${\bf WS}(0,1)$, we get the result for $\phi$. Moreover, by a similar argument and the fact that $H$ satisfies ${\bf WS}(0,2)$, we can deduce that $H$ satisfies ${\bf WS^{\infty}}(\alpha_2, \alpha_1 \wedge 2, 1/c_0)$. Next, we further assume that $\alpha_1<2$. Then, for all $0<s \le c_0$, we have that
\begin{align*}
H(s^{-1}) \asymp s^{-2} w(s)\int_0^s u\frac{w(u)}{w(s)}du \asymp s^{-2} w(s)\int_0^s s^{\alpha_1}u^{1-\alpha_1}du \asymp w(s).
\end{align*}

Now, suppose that $\phi$ satisfies ${\bf WS^{\infty}}(\alpha_2, \alpha_1, c_0)$ for some constants $0 \le \alpha_2 \le \alpha_1 < 1$ and $c_0>0$ or $H$ satisfies ${\bf WS^{\infty}}(\alpha_2, \alpha_1, c_0)$ for some constants $0 \le \alpha_2 \le \alpha_1 < 1$ and $c_0>0$. In either case, by \cite[Lemma 2.6 and  Proposition 2.9]{Mi}, $H$ satisfies ${\bf WS^{\infty}}(\alpha_2, \alpha_1, c_0)$ and there exists a constant $c_1>0$ such that $w(s) \asymp H(s^{-1})$ for $0<s<c_1$. Then, the result follows.

The cases when $w$ satisfies the weak scaling properties at infinity or either of $\phi$ and $H$ satisfies the weak scaling properties at the orgin can be proved by similar arguments.
\qed

\begin{lemma}\label{cauchy}
Suppose that there exist $\delta>0$ and $t_0>0$ such that $w$ satisfies ${\bf LS^0}(-\delta,t_0)$. Then, there exists a constant $c_1>0$ such that for every $t \in (0, t_0]$,
\begin{align*}
H(t^{-1})^{\delta+1} \le c_1\phi(t^{-1})^{\delta} w(t).
\end{align*}
Similarly, if there exist $\delta'>0$ and $t_0'>0$ such that $w$ satisfies ${\bf LS^\infty}(-\delta',t_0')$, then there exists a constant $c_2>0$ such that for every $t \in [t_0', \infty)$,
\begin{align*}
H(t^{-1})^{\delta+1} \le c_2\phi(t^{-1})^{\delta'} w(t).
\end{align*}
\end{lemma}
\pf
Since the proofs are similar, we only give the proof for the first assertion.
If $\delta<2$, then by Lemma \ref{l:phiHw}(ii), we have that for all $t \in (0, t_0]$,
\begin{align*}
H(t^{-1})^{\delta+1} \le cH(t^{-1})^{\delta} w(t)  \le c\phi(t^{-1})^{\delta} w(t).
\end{align*}
Now, assume that $\delta \ge 2$.
By Lemma \ref{l:phiHw}(i) and H$\ddot{\rm o}$lder's inequality, for every $t \in (0, t_0]$,
\begin{align*}
H(t^{-1}) &\le c t^{-2}\int_0^{t} sw(s)ds \le c t^{-2} \left(\int_0^t  w(s) ds\right)^{1-1/(\delta+1)}\left(\int_0^t  s^{\delta+1} w(s) ds\right)^{1/(\delta+1)}  \\
& \le ct^{-2}\big(t\phi(t^{-1})\big)^{1-1/(\delta+1)} \big(t^{\delta+2}w(t)\big)^{1/(\delta+1)}  =  c \phi(t^{-1})^{1-1/(\delta+1)}w(t)^{1/(\delta+1)}.
\end{align*}
We used Lemma \ref{l:phiHw}(i) and \cite[2.12.16]{BGT} in the third inequality. 
\qed

\begin{lemma}\label{l:decaycomp}
Suppose that there exist $\delta>0$ and $t_0>0$ such that $w$ satisfies ${\bf LS^\infty}(-\delta,t_0)$. Then, there exists a constant $c_1>0$ such that for every $t \in [t_0, \infty)$,
\begin{align*}
\phi(t^{-1})^{\delta+1} \le c_1 w(t).
\end{align*}
\end{lemma}
\proof We first assume that $\int_{1/(2t_0)}^{\infty} w(s) ds < \infty$. By Lemma \ref{l:phiHw}(i), we have that $\phi(t^{-1}) \asymp t^{-1}$ for all $t \ge t_0$. Then, by Potter's theorem, (see, \cite[Theorem 1.5.6]{BGT},) for all $t \ge t_0$,
\begin{align*}
\phi(t^{-1})^{\delta+1} \le c t^{-\delta-1} \le cw(t).
\end{align*}
Now, assume that $\int_{t_0/2}^{\infty} w(s) ds = \infty$. In this case, by Lemma \ref{l:phiHw}(i), $\phi(t^{-1}) \asymp t^{-1}\int_{t_0/2}^{t} w(s) ds$ for all $t \ge t_0$. We also have that by \cite[2.12.16]{BGT}, $w(t) \asymp t^{-\delta-1}\int_{t_0/2}^{t} s^{\delta}w(s) ds$ for all $t \ge t_0$. Then, by l'Hospital's rule and the fact that $w$ is non-increasing, we get
\begin{align*}
\limsup_{t \to \infty} \frac{w(t)}{\phi(t^{-1})^{\delta+1}} &\le c \limsup_{t \to \infty} \frac{t^{-\delta-1}\int_{t_0/2}^{t} s^{\delta}w(s) ds}{\left(t^{-1}\int_{t_0/2}^{t} w(s) ds\right)^{\delta+1}} =  c \limsup_{t \to \infty} \frac{\int_{t_0/2}^{t} s^{\delta}w(s) ds}{\left(\int_{t_0/2}^{t} w(s) ds\right)^{\delta+1}} \\
& \le c \limsup_{t \to \infty} \frac{t^{\delta}w(t)}{w(t)\left(\int_{t_0/2}^{t} w(s) ds\right)^{\delta}}  \le  c \limsup_{t \to \infty} \frac{t^{\delta}}{\big( tw(t_0/2) \big)^{\delta}} = c.
\end{align*}
\qed

For $s>0$, we define
\begin{align*}
b(s):= s \phi'(H^{-1}(1/s)).
\end{align*}

\begin{lemma}\label{l:b}
(i) $b$ is strictly increasing on $(0,\infty)$, $\lim_{s \to 0} b(s) = 0$ and $\lim_{s \to \infty} b(s) = \infty$.\\
(ii) For every $s>0$, we have that
\begin{align*}
\phi(s^{-1})^{-1} \le b^{-1}(s) \le \frac{e^2-e}{e-2} \phi(s^{-1})^{-1}.
\end{align*}
\end{lemma}
\proof
(i) Since $H$ is strictly increasing on $(0,\infty)$ and $\phi'$ is strictly decreasing on $(0,\infty)$, $b$ is strictly increasing on $(0,\infty)$. Moreover, we have that $\lim_{s \to 0} b(s) \le \phi'(H^{-1}(1)) \lim_{s \to 0} s= 0$ and $\lim_{s \to \infty} b(s) \ge \phi'(H^{-1}(1)) \lim_{s \to \infty} s= \infty$.

(ii) From the concavity of $\phi$, since $\phi^{-1}(\lambda) \le H^{-1}(\lambda)$, we have that for all $s>0$,
\begin{align*}
b(s) \le \frac{\phi^{-1}(s^{-1})}{\phi(\phi^{-1}(s^{-1}))} \frac{\phi(H^{-1}(s^{-1}))}{H^{-1}(s^{-1})}  \frac{1}{\phi^{-1}(s^{-1})} \le \frac{1}{\phi^{-1}(s^{-1})}.
\end{align*}
Therefore, we get $b^{-1}(s) \ge \phi(s^{-1})^{-1}$ since both $\phi$ and $b$ are strictly increasing.

On the other hand, we note that from the definition of $\phi$ and $H$, for every $\lambda>0$,
\begin{align*}
& \phi(\lambda) \le \lambda \int_{(0,1/\lambda]} u (-dw(u)) + w(1/\lambda), \\
&\phi'(\lambda) \ge e^{-1} \int_{(0,1/\lambda]} u (-dw(u)), \qquad    H(\lambda) \ge e^{-1}(e-2)w(1/\lambda).
\end{align*}
Let $a:=(e^2-e)/(e-2)$. Then, for all $s>0$,
\begin{align*}
b\big(a\phi(s^{-1})^{-1}\big)& = a\phi(s^{-1})^{-1}\phi'\big(H^{-1}\big(\phi(s^{-1})/a\big)\big) \\
&\ge ae^{-1}\phi(s^{-1})^{-1} \int_{(0,[H^{-1}(\phi(s^{-1})/a)]^{-1}]} u  (-dw(u)) \\
&\ge ae^{-1}\phi(s^{-1})^{-1} \bigg[ \int_{(0,s]} u  (-dw(u)) + s\int_{(s,[H^{-1}(\phi(s^{-1})/a)]^{-1}]}  (-dw(u)) \bigg] \\
& = ae^{-1}\phi(s^{-1})^{-1} \bigg[ \int_{(0,s]} u  (-dw(u)) + sw(s)-sw\big([H^{-1}(\phi(s^{-1})/a)]^{-1}\big) \bigg] \\
& \ge  ae^{-1}\phi(s^{-1})^{-1}\bigg[ s\phi(s^{-1})-e(e-2)^{-1}sH(H^{-1}\big(\phi(s^{-1})/a\big)\big) \bigg] \\
& = ae^{-1}(1-e(e-2)^{-1}a^{-1})s = s.
\end{align*}
Again, since $b$ is strictly increasing, we conclude that $b^{-1}(s) \le a\phi(s^{-1})^{-1}$.
\qed

We will use Chebyshev's inequality in tail probability estimates several times. To applying Chebyshev's inequality for subordinators, we need the following lemma.

\begin{lemma}\label{l:laplace}
Assume that $w$ is finitely supported, that is, there exists a constant $T>0$ such that $w(T)=0$. Then, for every $\lambda \in \R$, $r>0$ and $n \in \{0\} \cup \N$, we have that
\begin{align*}
\E[(S_r)^ne^{\lambda S_r}] = \frac{d^n}{d \lambda^n}\exp \left(r \int_{(0,T]} (e^{\lambda s} -1) (-dw(s)) \right).
\end{align*}

\end{lemma}
\proof
Fix $r>0$ and let $\xi(dt):=\P(S_r \in dt)$. For $z \in \mathbb{C}$, define
\begin{align*}
f(z)=\int_{[0,\infty)} e^{-zt} \xi(dt).
\end{align*}
Then, it is well known that there exists the abscissa of convergence $\sigma_0 \in [-\infty, \infty]$ such that $f(z)$ converges for $\text{Re} \; z > \sigma_0$, diverges for $\text{Re} \; z < \sigma_0$ and has a singularity at $\sigma_0$. Moreover, $f(z)$ is analytic in the half-plane $\text{Re} \; z>\sigma_0$ so that for every $n \in \N$ and $x>\sigma_0$, it holds that
\begin{align}\label{diffL}
\frac{d^n}{dx^n} f(x) = (-1)^n\int_{[0,\infty)} t^n e^{-xt} \xi(dt). 
\end{align}
(See, \cite[p.37 and p.58]{Wi} and \cite{Mi2}.) On the other hand, we also have that for $\lambda>0$,
\begin{align*}
f(\lambda) = \E\big[\exp(-\lambda S_r)\big] = \exp \big(-r \phi(\lambda)\big) = \exp \left(r \int_{(0,T]} (e^{\lambda s} -1) (-dw(s)) \right)=:g(\lambda).
\end{align*}
Since $w$ is finitely supported, the function $\lambda \mapsto g(\lambda)$ is a well-defined differentiable function on $\R$. If $\sigma_0>-\infty$, then from the uniqueness of the analytic continuation, the function $g(\lambda)$ should have a singularity at $\lambda = \sigma_0$. Since there is no such singularity, we get $\sigma_0 = -\infty$. Then, the result follows from the definition of $f$ and \eqref{diffL}.
\qed

\subsection{Tail probability estimates for subordinator.}

In this section, we study two tail probabilities $\P(S_r \ge t)$ and $\P(S_r \le t)$ under mild assumption for $w$. We first give the general lower bounds for upper tail probability $\P(S_r \ge t)$ which are established in \cite{Mi}. Note that this bounds hold for every subordinator.

\begin{lemma}\label{l:tail0}
For every $L>0$, it holds that for all $r, t>0$ satisfying $r\phi(t^{-1}) \le L$,
\begin{align*}
\P(S_r \ge t) \ge e^{-eL} rw(t).
\end{align*}
\end{lemma}
\proof 
Note that $r\phi(t^{-1}) \le L$ implies that $rw(t) \le er\phi(t^{-1}) \le eL$. Thus, by \cite[Proposition 2.5]{Mi}, for all $r, t>0$ satisfying $r\phi(t^{-1}) \le L$, we have that
\begin{align*}
\P(S_r \ge t) \ge 1 - e^{-rw(t)} \ge rw(t) e^{-rw(t)} \ge e^{-eL}rw(t).
\end{align*}
\qed

Now, we study the upper bounds for $\P(S_r \ge t)$.

\begin{prop}\label{p:utail1}
Assume that condition {\bf (S.Poly.)($t_s$)} holds. Then, there exists a constant $c_1>0$ such that for all $r,t>0$ satisfying $0<t \le t_s$ and $r\phi(t^{-1}) \le 1/(4e^2)$,
\begin{align*}
\P(S_r \ge t) \le c_1 rw(t).
\end{align*}
\end{prop}
\proof
Fix $r, t>0$ sastisfying $0<t \le t_s$ and $r\phi(t^{-1}) \le 1/(4e^2)$. Set
\begin{align*}
\mu^1:=\1_{(0,  1/H^{-1}(1/r)]}\cdot (-dw), \quad \mu^2:=\1_{(1/H^{-1}(1/r), t]} \cdot (-dw), \quad \mu^3:= \1_{(t, \infty)} \cdot (-dw).
\end{align*}
Let $S^1, S^2$ and $S^3$ be independent subordinators without drift and having L\'evy measure $\mu^1, \mu^2$ and $\mu^3$, respectively. Then, we have $S_r \le S^1_r + S^2_r + S^3_r$ and hence
\begin{align}\label{eq}
\P(S_r \ge t) &\le \P(S^1_r + S^2_r + S^3_r \ge t) \le \P(S^1_r \ge t/2) + \P(S^2_r \ge t/2) + \P(S^3_r >0).
\end{align}

First, since $S^3$ is a compounded Poisson process, $\P(S^3_r>0) = 1 - e^{-rw(t)} \le rw(t)$.

Next, we note that by Lemma \ref{l:b}(ii), $t =b(b^{-1}(t)) \ge b(\phi(t^{-1})^{-1}) \ge b(4e^2r) \ge 4e^2b(r)$. By Chebyshev's inequality and Lemma \ref{l:laplace}, we have that for every $\lambda>0$,
\begin{align*}
&\P\big(S^1_r \ge t/2\big) \le \E\big[\exp\big(-\lambda t/2 + \lambda S^1_r\big)\big] = \exp\bigg(-\frac{\lambda t}{2} + r\int_{(0,1/H^{-1}(1/r)]} (e^{\lambda s} - 1) (-dw(s))\bigg)\\
& \le \exp\bigg(-\frac{\lambda t}{2} + \lambda re^{\lambda/H^{-1}(1/r)} \int_{(0,1/H^{-1}(1/r)]} s (-dw(s)) \bigg)  \le \exp\bigg(-\frac{\lambda t}{2} + e \lambda b(r)  e^{\lambda/H^{-1}(1/r)}  \bigg). 
\end{align*}
Thus, by letting $\lambda=H^{-1}(1/r)$, we get
\begin{align*}
\P\big(S^1_r \ge t/2\big) & \le \exp \big(-2^{-1}tH^{-1}(1/r)+e^2b(r)H^{-1}(1/r)\big)\le \exp \big(-4^{-1}tH^{-1}(1/r)\big).
\end{align*}

Thirdly, let $f_0(s):=w(s) \1_{(0,t]}(s) + w(t)t^2s^{-2}\1_{(t,\infty)}(s)$ for $s>0$. Then, we see that $f_0$ is non-increasing and for every Borel set $A \subset \R$, it holds that
$$
\mu^2(A) \le w(\dist(0,A)) \1_{(0,t]}(\dist(0,A)) \le f_0(\dist(0,A)),
$$
where $\dist(0,A):=\inf\{|y|:y \in A\}$. Moreover, since $w$ satisfies ${\bf LS^0}(-\delta,t_s)$, for all $u,v>0$,
\begin{align*}
\int_u^{\infty} f_0\big(v \vee y - \frac{y}{2}\big) \mu^2(dy) \le f_0(v/2) w(u) \le c_1 f_0(v)H(1/u).
\end{align*}
Therefore, by \cite[Proposition 1 and Lemma 9]{KS}, we have that for every $x>0$ and $\rho \in (0,x/3]$,
\begin{align*}
\P(S^2_r \in [x-\rho,x+\rho]) \le c_2 r f_0(x/3).
\end{align*}
It follows that
\begin{align*}
\P(S^2_r \ge t/2) &\le \sum_{i=0}^{\infty} \P\big(S^2_r \in [2^{i-1}t,2 \cdot 2^{i-1}t]\big) \le cr \sum_{i=0}^{\infty} f_0(2^{i-2} t) \le crw(t) \sum_{i=0}^{\infty} 2^{-2i} = crw(t).
\end{align*}

Combining the above inequalities, by \eqref{eq} and Lemma \ref{cauchy}, we deduce that
\begin{align*}
\P(S_r \ge t) & \le crw(t) + \exp \bigg(-2(\delta+1) \cdot \frac{tH^{-1}(1/r)}{8(\delta+1)}\bigg) \le crw(t) + \left(\frac{1 \vee 8(\delta+1)}{1 \vee tH^{-1}(1/r)}\right)^{2\delta+2} \\
& \le crw(t) + c\left(\frac{H(t^{-1})}{H(H^{-1}(1/r))}\right)^{\delta+1} \le crw(t) + cr \phi(t^{-1})^{-\delta} H(t^{-1})^{\delta+1} \le crw(t).
\end{align*}

In the second inequality, we used the fact that $e^x \ge x$ for all $x>0$ and in the third inequality, we used the fact that $H(\lambda x) \le (1 \vee \lambda^2) H(x)$ for all $\lambda, x>0$. Also, the fourth inequality holds since $r \le L \phi(t^{-1})^{-1}$.
\qed

\smallskip

By the same argument, we also get analogous estimates for large time $t$.

\begin{prop}\label{p:utail2}
Assume that condition {\bf (L.Poly.)} holds. Then, for every $T>0$, there exists a constant $c_1>0$ such that for all $r,t>0$ satisfying $t \ge T$ and $r\phi(t^{-1}) \le 1/(4e^2)$,
\begin{align*}
\P(S_r \ge t) \le c_1 rw(t).
\end{align*}
\end{prop}
\proof Follow the proof of Proposition \ref{p:utail1}. The only difference occurs in the definition of $f_0$. In this case, we use $f_1(s):=\frac{e}{e-2}H(s^{-1})\1_{(0,T/2]}(s) + w(s) \1_{(T/2, \infty)}(s)$ instead of $f_0(s)$.
\qed

\begin{prop}\label{p:utail3}
Assume that condition {\bf (Sub.)($\bbeta, \bth$)} holds. Then, for every $T>0$, there exist constants $c_2>0$ and $L \in (0,1]$ such that 
for all $r,t>0$ satisfying $t \ge T$ and $rt^{-1} \le L$,
\begin{align*}
\P(S_r \ge t) \le c_2r \exp\big(-\frac{\theta}{2}t^\beta \big).
\end{align*}
\end{prop}
\proof
Fix $t \ge T$ and $r \in (0, Lt)$ where the constant $L \in (0,1]$ will be chosen later. Let $\widehat{S}^1$ and $\widehat{S}^2$ be independent subordinators without drift and having L\'evy measures
\begin{align*}
\widehat{\mu}^1 := \1_{(0,t]} \cdot (-dw)  \quad \text{and} \quad   \widehat{\mu}^2 := \1_{(t,\infty)} \cdot (-dw), \qquad \text{respectively.} 
\end{align*}
Then, since $S_r = \widehat{S}^1_r + \widehat{S}^2_r$, by condition {\bf (Sub.)($\bbeta, \bth$)}, we have
$$
\P(S_r \ge t) \le \P(\widehat{S}^1_r \ge t) + \P(\widehat{S}^2_r > 0) \le \P(\widehat{S}^1_r \ge t) + rw(t) \le \P(\widehat{S}^1_r \ge t) + cr\exp \big(-\theta t^\beta \big).
$$

It remains to bound $\P(\widehat{S}^1_r \ge t)$. By Chebyshev's inequality and Lemma \ref{l:laplace}, for all $\lambda > 0$,
\begin{align}\label{e:Che}
\P(\widehat{S}^1_r \ge t) &\le \E\big[t^{-1}\widehat{S}^1_r \exp\big(\lambda\widehat{S}^1_r  - \lambda t \big)\big] \nn\\
&\le t^{-1}e^{-\lambda t} r \bigg(\int_{(0,t]} se^{\lambda s} (-dw(s))\bigg) \exp \bigg(r \int_{(0,t]} (e^{\lambda s}-1)(-dw(s))\bigg).
\end{align}
Note that by the integration by parts and condition {\bf (Sub.)($\bbeta, \bth$)}, we get
\begin{align*}
&\int_{(0,t]} se^{\lambda s} (-dw(s)) \le \int_{(0,t]} w(s)e^{\lambda s}ds + \lambda \int_{(0,t]} sw(s)e^{\lambda s}ds \\
&\le 2\lambda e^\lambda \int_{(0,1]} w(s)ds + c_0 \int_{(1,t]} \exp \big(-\theta s^\beta + \lambda s \big)ds + c_0\lambda \int_{(1,t]} s \exp \big(-\theta s^\beta + \lambda s \big)ds,
\end{align*}
and
\begin{align*}
&\int_{(0,t]} (e^{\lambda s}-1) (-dw(s)) \le \lambda e^\lambda \int_{(0,1]} w(s)ds  + c_0\lambda \int_{(1,t]} \exp \big(-\theta s^\beta + \lambda s \big)ds.
\end{align*}
Take $\lambda = 2\theta t^{\beta-1}/3 \in (0, 2\theta T^{\beta-1}/3]$. Then, since $s \mapsto -2\theta s^\beta/3 + \lambda s$ is a convex function,
\begin{align*}
&\int_{(1,t]} s\exp \big(-\theta s^\beta + \lambda s \big)ds \le \sup_{s \in (1,t]} \bigg[-\frac{2\theta s^\beta}{3} + \lambda s \bigg] \cdot \int_{(1,t]} s \exp \big(-\frac{\theta s^\beta}{3} \big)ds \\
&\qquad \le  \big( -\frac{2\theta }{3} + \lambda  -\frac{2\theta t^\beta}{3} + \lambda t \big)\int_{(1,t]} s \exp \big(-\frac{\theta s^\beta}{3} \big)ds \le c.
\end{align*}
Using this observation and the fact that $\int_{(0,1]} w(s)ds < \infty$, \eqref{e:Che} implies that
\begin{align*}
\P(\widehat{S}^1_r \ge t) \le c_3 t^{-1}r \exp \big(-\frac{2\theta }{3}t^\beta + c_4rt^{\beta-1} \big), 
\end{align*}
for some constants $c_3, c_4>0$. Now, we choose $L=1 \wedge (\theta /(6c_4))$. Then, we get
\begin{align*}
\P(\widehat{S}^1_r \ge t) \le c_3 T^{-1}r \exp \big(-\frac{2\theta }{3}t^\beta + c_4Lt^{\beta} \big) \le c_2r \exp \big(-\frac{\theta}{2}t^\beta \big).
\end{align*}
\qed

When $w$ decreases subexponentially ($0<\beta<1$), we obtain small time sharp upper bounds for $\P(S_r \ge t)$ which decrease with exactly the same rate as the bounds for $w$ as $t \to \infty$. 

\begin{prop}\label{p:utail4}
Assume that condition {\bf (Sub.)($\bbeta, \bth$)} holds with constant $0<\beta<1$. Then, for every fixed $k>0$ and $T>0$, there exist constants $c_2>0$ and $L \in (0,1]$ such that 
for all $r,t>0$ satisfying $t \ge T$ and $rt^{-1} \le L$,
\begin{align*}
\P(S_r \ge t) \le c_2r \exp\big(-\theta t^{\beta}+kr\big).
\end{align*}
\end{prop}
\proof
Let $\widetilde{S}^1$ and $\widetilde{S}^2$ be independent subordinators without drift and having L\'evy measures
\begin{align*}
\widetilde{\mu}^1 := \1_{(0,t/2]} \cdot (-dw)  \quad \text{and} \quad   \widetilde{\mu}^2 := \1_{(t/2,\infty)} \cdot (-dw), \qquad \text{respectively.} 
\end{align*}
Then, since $S_r = \widetilde{S}^1_r + \widetilde{S}^2_r$, we get
\begin{align}\label{e:eq2}
\P(S_r \ge t) &= \int_{0}^{\infty} \P(\widetilde{S}^2_r \ge t-u) \P(\widetilde{S}^1_r \in du) \nn\\
&\le  \P(\widetilde{S}^2_r \ge t-T/2) + \int_{T/2}^{t-T/2} \P(\widetilde{S}^2_r \ge t-u) \P(\widetilde{S}^1_r \in du) + \P(\widetilde{S}^1_r \ge t-T/2).
\end{align}

By Chebyshev's inequality, Lemma \ref{l:laplace} and the integration by parts, for $u>0$ and $\lambda > 0 $,
\begin{align*}
\P(\widetilde{S}^1_r \ge u) &\le \E\big[u^{-2} (\widetilde{S}^1_r)^2 \exp \big(- \lambda u + \lambda \widetilde{S}^1_r\big)\big] \\
& = u^{-2} \bigg[r \int_{(0,t/2]} s^2e^{\lambda s} (-dw(s)) + \bigg(r \int_{(0,t/2]} se^{\lambda s} (-dw(s))\bigg)^2 \bigg] \\
&\quad \times \exp \bigg(-\lambda u +  r \int_{(0, t/2]} (e^{\lambda s}-1) (-dw(s))\bigg)\\
& \le u^{-2} \bigg[r \int_0^{t/2} (2+\lambda s)se^{\lambda s} w(s)ds + \bigg(r \int_0^{t/2} (1+\lambda s)e^{\lambda s} w(s) ds\bigg)^2 \bigg] \\
&\quad \times \exp \bigg(-\lambda u + \lambda r \int_0^{t/2} e^{\lambda s}w(s) ds \bigg).
\end{align*}
Take $\lambda = \theta t^{\beta-1} \in (0, \theta T^{\beta-1}]$. Then, for all $1 \le s \le t/2$, we have that $\lambda s \le \theta t^{\beta-1}(t/2)^{1-\beta}s^{\beta} \le 2^{\beta-1}\theta s^\beta$. It follows that 
\begin{align*}
\int_0^{t/2} (2+\lambda s)se^{\lambda s} w(s)ds &\le (2+\lambda)e^{\lambda} \int_0^1 w(s)ds + c_0(2+\lambda)\int_1^{t/2} s^2\exp \big(2^{\beta-1}\theta s^{\beta} - \theta s^\beta\big) ds \\
& \le c + c\int_1^{\infty} s^2 \exp\big(-\theta (1-2^{\beta-1})s^\beta \big) ds \le c_4,
\end{align*}
where the constant $c_4>0$ is independent of $t \in [T, \infty)$. By similar calculations, by taking $c_4$ larger, we may assume that
\begin{align*} 
\int_0^{t/2}(1+\lambda s)e^{\lambda s} w(s) ds \le c_4 \quad \text{and} \quad  \int_0^{t/2} e^{\lambda s}w(s) ds \le c_4.
\end{align*}
Therefore, we have that for every $u>0$,
\begin{align*}
\P(\widetilde{S}^1_r \ge u) &\le (c_4+c_4^2)u^{-2}(r+r^2) \exp \big(-\theta t^{\beta-1}u + \theta c_4 t^{\beta-1} r \big).
\end{align*}
In particular,
\begin{align*}
\P(\widetilde{S}_r^1 \ge t-T/2) & \le  ct^{-2}(r+r^2)\exp\big( -\theta t^{\beta} + \theta t^{\beta-1}T/2 + \theta c_4t^{\beta-1}r\big) \\
& \le cT^{-2}\exp \big(\theta T^{\beta}/2\big) r \exp\big( -\theta t^{\beta} + (\theta c_4t^{\beta-1}+k/2)r\big).
\end{align*}

On the other hand, note that $\widetilde{S}^2_r = \sum_{i=1}^{N(r)} D_i$ where $N(r)$ is a Poisson process with rate $w(t/2)$ and $D_i$ are i.i.d. random variables with distribution $\P(D_i > u) = w\big(u \vee (t/2)\big)/w\big(t/2\big)$. Thus, for every $0<u<t$,
\begin{align*}
\P(\widetilde{S}_r^2 \ge u) &\le \P(N(r) = 1, \; D_1 \ge u) + \P(N(r) \ge 2) \\
&\le rw(u \vee (t/2)) + 1 - e^{-rw(t/2)} - rw(t/2)e^{-rw(t/2)} \\
&\le cr \exp \big( -\theta u^{\beta} \big) + r^2w(t/2)^2 \le cr\exp \big( -\theta u^{\beta} \big) + cLrt \exp \big( -\theta 2^{1-\beta} t^\beta \big) \\
&\le cr\exp \big( -\theta u^{\beta} \big) + cLr \exp \big( -\theta t^\beta \big).
\end{align*}
It follows that
\begin{align*}
\P(\widetilde{S}_r^2 \ge t-T/2) & \le  cr \exp\big( -c_1(t-T/2)^\beta\big) \le cr \exp\big( -c_1t^{\beta}\big).
\end{align*}
The second inequality holds since $t^{\beta} - (t-T/2)^\beta \le (T/2)^\beta$.

Using the above inequalities, by \eqref{e:eq2} and the integration by parts, we obtain
\begin{align*}
&\P(S_r \ge t) \le  c r \exp\big( -\theta t^{\beta} + (\theta c_4t^{\beta-1}+\frac{k}{2})r\big) + cr \int_{T/2}^{t}\exp \big( -\theta (t-u)^{\beta} \big) \P(\widetilde{S}^1_r \in du) \\
& \le c r \exp\big( -\theta t^{\beta} + (\theta c_4t^{\beta-1}+\frac{k}{2})r\big)+ cr \int_{T/2}^{t} \P(\widetilde{S}^1_r \ge u) (t-u)^{\beta-1} \exp \big( -\theta (t-u)^{\beta} \big) du \\
& \le c r \exp\big( -\theta t^{\beta} + (\theta c_4t^{\beta-1}+\frac{k}{2})r\big) \bigg(1+ c(T/2)^{\beta-1} \int_{T/2}^{t}u^{-2}\exp \big( - f(u) \big) du\bigg),
\end{align*}
where $f(u):= \theta (t-u)^\beta + \theta t^{\beta-1}u - \theta t^\beta$. Observe that
\begin{align*}
f'(u) = - \beta \theta (t-u)^{\beta-1} + \theta t^{\beta-1} = -\theta t^{\beta-1}(t-u)^{\beta-1} (\beta t^{1-\beta} - (t-u)^{1-\beta}).
\end{align*}
Hence, $f$ is decreasing on $(0, (1-\beta^{1/(1-\beta)})t)$ and increasing on $((1-\beta^{1/(1-\beta)})t, t)$. Since $f(0) = f(t) =0$, we deduce that $f(u) \le 0$ for $u \in (0,t)$ and hence $\int_{T/2}^{t}u^{-2}\exp \big( - f(u) \big) du \le \int_{T/2}^{\infty} u^{-2} du \le c$. It follows that
\begin{align*}
\P(S_r \ge t) \le  c r \exp\big( -\theta t^{\beta} + (\theta c_4t^{\beta-1}+\frac{k}{2})r\big).
\end{align*}
Hence, if $t \ge (k/(2\theta c_4))^{-1/(1-\beta)}=:c_5$, we are done. Moreover, if $t<c_5$, then we get
\begin{align*}
\exp \big((\theta c_4t^{\beta-1}+\frac{k}{2})r \big) \le \exp \big(\theta c_4t^{\beta} + \frac{k}{2}t\big)\le c,
\end{align*}
since $r \le Lt \le t$. This completes the proof.

\qed

Here, we state the estimates on lower tail probabilities $\P(S_r \le t)$ when $r$ is large enough compare to $b^{-1}(t)$, which are established in \cite{JP}.

\begin{lemma}\label{l:Offprob}{\rm\cite[Lemma 5.2]{JP}}
For every $N>0$, there exist constants $c_1,c_2>0$ such that
\begin{align*}
c_1 \exp\big(- c_2 r H((\phi')^{-1}(t/r))\big) \le \P(S_r \le t)  \le \exp\big(- r H((\phi')^{-1}(t/r))\big),
\end{align*}
for all $r,t>0$ satisfying $ r \ge Nb^{-1}(t)$.
\end{lemma}
\proof
If $N \ge 1$, then $ r \ge Nb^{-1}(t)$ implies that $rH((\phi')^{-1}(t/r)) \ge r H((\phi')^{-1}(b(r)/r)) = 1$ and hence the result follows from \cite[Lemma 5.2]{JP}. Suppose that $N \in (0,1)$. Since $r \mapsto S_r$ is strictly increasing almost surely, we deduce that for all $r \in (Nb^{-1}(t), b^{-1}(t)]$,
\begin{align*}
\P(S_r \le t) \ge \P(S_{b^{-1}(t)} \le t) \ge c \ge c \exp\big(- c_2 r H((\phi')^{-1}(t/r))\big).
\end{align*}
\qed

\begin{cor}\label{c:meandist}
If condition {\bf (S.Poly.)($t_s$)} holds, then there exist constants $N > \varepsilon_1>0$ such that for all $t \in (0,t_s]$, it holds that
\begin{align}\label{meandist}
\P(S_{N/\phi(t^{-1})} \ge t) -\P(S_{\varepsilon_1/\phi(t^{-1})} \ge t) \ge 1/4.
\end{align}

On the other hand, if either of the conditions {\bf (L.Poly.)} or {\bf (Sub.)($\bbeta, \bth$)} holds, then for every fixed $T>0$, there exist constants $N > \varepsilon_1>0$ such that \eqref{meandist} holds for all $ t \in [T, \infty)$.
\end{cor}
\proof
By Lemma \ref{l:Offprob} and Lemma \ref{l:b}(ii), there exists a constant $N>0$ such that for all $t>0$, $\P(S_{N/\phi(t^{-1})} < t) \le 1/4$ and hence $\P(S_{N/\phi(t^{-1})} \ge t) \ge 3/4$. On the other hand, by Proposition \ref{p:utail1} (resp. Proposition \ref{p:utail2} or Proposition \ref{p:utail3}) and the facts that $\phi(t^{-1}) \asymp t^{-1}$ for all $t \ge T$ under condition {\bf (Sub.)($\bbeta, \bth$)} and $\phi(t^{-1}) \ge e^{-1}w(t)$ for all $t>0$, we can find a constnat $\varepsilon_1>0$ such that $\P(S_{\varepsilon_1/\phi(t^{-1})} \ge t) \le 1/2$ for all $t \in (0, t_s]$ (resp. for all $t \in [T, \infty)$).
\qed

By Corollary \ref{c:meandist}, we get a priori estimates for the fundamental solution $p(t,x,y)$.

\begin{cor}\label{c:prebound}
Assume that condition {\bf (S.Poly.)($t_s$)} holds. Let $p(t,x,y)$ be given by \eqref{e:dfundamental}. Then, there exist constants $N > \varepsilon_1>0$ and $c>0$ such that for all $t \in (0,t_s]$, 
\begin{align}\label{prebound}
p(t,x,y) \ge c \inf_{r \in (\varepsilon_1/\phi(t^{-1}), N/\phi(t^{-1}))} q(r,x,y).
\end{align}

On the other hand, if either of the conditions {\bf (L.Poly.)} or {\bf (Sub.)($\bbeta, \bth$)} holds, then for every fixed $T>0$, there exist constants $N > \varepsilon_1>0$ such that \eqref{prebound} holds for all $ t \in [T, \infty)$.
\end{cor}

\vspace{1.5mm}

\subsection{Estimates for truncated subordinator.}\label{s:fs}

In this subsection, we obtain tail probability estimates when the kernel $w$ is finitely supported. Throughout this subsection, we always assume that condition {\bf (Trunc.)($t_f$)} holds.
An example of such kernel is given by $w(t):= \frac{1}{\Gamma(1-\beta)} (t^{-\beta} - 1) \1_{(0,1]}(t)$ ($0<\beta<1$). Those integral kernels are used in the fractional-time derivative whose value at time $t$ depends only on the finite range of the past. (See, ~ \cite[Example 2.5]{Ch}.)

\begin{prop}\label{p:trunc}
There exists a constant $r_0 > 0$ such that for all $r \in (0,r_0]$ and $t \ge t_f/2$,
$$
\P(S_r \ge t) \simeq  \big[r + (nt_s - t)^n \big]r^n\exp \big(-ct \log t\big),
$$
where $n:=\lfloor t/t_f \rfloor + 1$.
\end{prop}
\proof
Take $r_0$ small enough so that $r\phi(r^{-1}) \le 1/(4e^2)$ and $r \le t_f/6$ for all $r \in (0,r_0]$. Since $ \lim_{ r \to 0 } r\phi(r^{-1}) =  0$, we can always find such constant $r_0$.
Then, fix $r \in (0,r_0]$ and $t \ge t_f/2$. Note that since $n= \lfloor t/t_f \rfloor + 1$, we have $\big((n-1)\vee 1/2\big)t_f \le t < nt_f$. 

\smallskip
(Lower bound)
Let $U^1$ and $U^2$ be the driftless subordinators with L\'evy measures
\begin{align*}
\nu_1:=\1_{(t/(n+1),\infty)} \cdot (-dw) \quad \text{and} \quad \nu_2:=\1_{(t/n,\infty)} \cdot (-dw), \qquad \text{respectively.}
\end{align*}
Observe that both $U^1$ and $U^2$ are compounded Poisson processes and thier jump sizes are at least bigger than $t/(n+1)$ and $t/n$, respectively. Since $S_r \ge U^1_r \ge U^2_r$, it follows that
\begin{align}\label{poisson}
&2\P(S_r \ge t) \ge \P(U^1_r \ge t) + \P(U^2_r \ge t) \nn\\
&\ge \P\big(U^1 \; \text{jumps} \; (n+1) \; \text{times before time} \; r\big) + \P\big(U^2 \; \text{jumps} \; n \; \text{times before time} \; r\big) \nn\\
&\ge \exp \big(-r w(t/(n+1)) \big) \frac{\big(rw(t/(n+1))\big)^{n+1}}{(n+1)!} +  \exp \big(-r w(t/n) \big) \frac{\big(rw(t/n)\big)^{n}}{n!}.
\end{align}
Since $s \mapsto w(s)$ is non-increasing, we have $w(t/(n+1)) \le w(t_f/4)$ and $w(t/n) \le  w(t_f/2)$.
Moreover, by condition {\bf (Trunc.)($t_f$)}(i) and (ii),
\begin{align*}
&w(t/(n+1)) \ge K^{-1} (t_f - t/(n+1)) \ge K^{-1}(n+1)^{-1}t_f, \\
&w(t/n) \ge K^{-1} (t_f - t/n) \ge K^{-1}n^{-1}(nt_f-t).
\end{align*}
Using these observations, Stirling's formula and the fact that $n \asymp t$, by \eqref{poisson}, we obtain
\begin{align*}
\P(S_r \ge t) &\ge e^{-rw(t_f/4)} \frac{t_f^{n+1}r^{n+1}}{2K^{n+1}(n+1)^{n+1}(n+1)!} +  e^{-rw(t_f/2)} \frac{(nt_f-t)^{n}r^{n}}{2K^{n}n^{n}n!} \\
& \gtrsim r^{n+1}\exp \big(-ct - 2n \log n \big) + (nt_f-t)^nr^n \exp \big(-ct - 2n \log n\big) \\
& \gtrsim \big[ r  + (nt_f-t)^n \big] r^n \exp \big(- ct  \log t). 
\end{align*}

\smallskip

(Upper bound) Let $U^3$ and $U^4$ be the driftless subordinators with L\'evy measures
\begin{align*}
\nu_3:=\1_{(0,t/9]} \cdot (-dw) \quad \text{and} \quad \nu_4:=\1_{(t/9,\infty)} \cdot (-dw), \qquad \text{respectively.}
\end{align*}
Then, we have that $S_r = U^3_r + U^4_r$ and $U^4_r = \sum_{i=1}^{P(r)} J_i$ where $P(r)$ is a Poisson process with rate $w(t_f/9)$ and $J_i$ are i.i.d. random variables with distribution
$$
F(u):=\P(J_i \ge u) = w(t_f/9)^{-1}w\big(u \vee (t_f/9)\big).
$$ Hence, we get
\begin{align*}
\P(S_r \ge t)  &= \sum_{j=0}^\infty \P\big(U^3_r + U^4_r \ge t, P(r) = j\big)\\
&\le \P\big(U^3_r \ge t\big) + \sum_{j=1}^n \P\big(U^3_r + U^4_r \ge t | P(r) = j\big) \P\big(P(r) = j\big) + \P\big(P(r) > n\big).
\end{align*}

First, by Stirling's formula, the definition of Poisson process and the fact that $n \asymp t$,
\begin{align*}
\P(P(r)>n) \le \frac{er^{n+1}}{(n+1)!} \simeq r^{n+1} \exp \big(-ct \log t \big).
\end{align*}

Secondly, by Chebyshev's inequality and Lemma \ref{l:laplace}, for all $u>0$ and $\lambda>0$,
\begin{align*}
\P(U^3_r \ge u) &\le \E\big[\exp\big(-\lambda u + \lambda U^3_r \big) \big] = \exp \big(-\lambda u + r \int_{(0,t_s/9]} (e^{\lambda s}-1) (-dw(s)) \big) \\
&\le \exp\big(-\lambda u + \lambda e^{\lambda t_f/9} r \int_{(0,t_f/9]} s (-dw(s))\big) \le \exp\big(-\lambda u + c_1 \lambda e^{\lambda t_f/9} r \big).
\end{align*}
Hence, by taking $\lambda=9t_f^{-1} \log\big(u/(9c_1r)\big)$, we have that for every $u>0$,
\begin{align}\label{smalljump}
\P(U^3_r \ge u) \le \exp \big(-8\lambda u/9 \big) = \big(9c_1r/u\big)^{8u/t_f}.
\end{align}
In particular, since $t \ge \big((n-1) \vee 1/2 \big)t_f$, we have that 
\begin{align*}
\P(U^3_r \ge t) \le \big(9c_1r/t\big)^{8t/t_f} \lesssim r^{8t/t_f} \exp \big(-ct \log t\big) \le cr^{n+1}\exp \big(-ct \log t\big).
\end{align*}
Moreover, we also have that
\begin{align*}
\sum_{j=1}^{n-2} & \; \P\big(U^3_r + U^4_r \ge t | P(r) = j\big) \P\big(P(r) = j\big) \le \sum_{j=1}^{n-2} \frac{r^jw(t_f/9)^j}{j!} \P\big(U^3_r \ge (n-1-j)t_f \big) \\
&\le \sum_{j=1}^{n-2} \frac{r^jw(t_f/9)^j}{j!} \left(\frac{9c_1r}{(n-j-1)t_f}\right)^{8(n-1-j)} \lesssim e^{ct} \sum_{j=1}^{n-2} r^{8(n-1)-7j}\frac{1}{j!(n-j-1)^{8(n-j-1)}}  \\
&\lesssim r^{n+1}\sum_{j=1}^{n-2} \exp \big( ct - c j \log j - c (n- j - 1) \log (n-j-1) \big) \\
& \lesssim r^{n+1} \exp \big(ct - c (n-1) \log (n-1) \big) \simeq r^{n+1} \exp \big( - ct \log t \big).
\end{align*}
The first inequality holds since the jump sizes of $U_r^4$ are at most $t_f$ and the third line follows form Stirling's formula. Lastly, the fourth line holds by the facts that $4(a \log a + b \log b) \ge 2(a \vee b) \log (2 (a\vee b)) \ge (a+b) \log (a+b)$ for all $a, b \ge 1$ satisfying $a \vee b \ge 2$ and that $n \asymp t$.

It remains to bound probabilities $\P\big(U^3_r + U^4_r \ge t, P(r) = j\big)$ for $j = n-1$ (when $n \ge 2$) and $j=n$. Observe that by Stirling's formula, we have
\begin{align*}
&\P\big(U^3_r + U^4_r \ge t | P(r) = n-1\big)\P\big(P(r)=n-1\big) \\
& \le \frac{r^{n-1}w(t_f/9)^{n-1}}{(n-1)!} \int_{0}^{(n-1)t_f} \P\big(U^3_r \ge t-(n-1)t_f+ u\big) d_u \P\big(\sum_{i=1}^{n-1} J_i \ge (n-1)t_f-u\big)\\
& \lesssim r^{n-1} \exp \big( - ct \log t \big)  \\
& \quad \times \bigg[\left(\int_{0}^{r} + \int_{t_f/4}^{(n-1)t_f} + \int_{r}^{t_f/4}\right) \P\big(U^3_r \ge t-(n-1)t_f+ u\big) d_u \P\big(\sum_{i=1}^{n-1} J_i \ge (n-1)t_f-u\big)\bigg]\\
& \le r^{n-1} \exp \big( - ct \log t \big) \bigg[\P(U^3_r \ge t_s/4) + \P\big(U^3_r \ge t- (n-1)t_f\big)\P\big(\sum_{i=1}^{n-1} J_i \ge (n-1)t_f-r\big)  \\
& \qquad \qquad\qquad \qquad\qquad\;  +\int_{r}^{t_f/4} \P\big(U^3_r \ge t-(n-1)t_f + u\big) d_u \P\big(\sum_{i=1}^{n-1} J_i \ge (n-1)t_f-u\big) \bigg]\\
& =: r^{n-1} \exp \big( - ct \log t \big)\big[A_1 + A_2 + A_3\big]
\end{align*}
and by the same way, we also have that
\begin{align*}
&\P\big(U^3_r + U^4_r \ge t |  P(r) = n\big)\P\big(P(r)=n\big)  \\
&\qquad \lesssim r^{n} \exp \big( - ct \log t \big) \bigg[\P(U^3_r \ge t_f/4) + \P\big(\sum_{i=1}^{n} J_i \ge nt_f-(nt_f-t+r)\big)  \\
&\qquad \qquad \qquad\qquad \qquad\quad  +\int_{nt_f-t+r}^{t_f/4} \P\big(U^3_r \ge t-nt_f + u\big) d_u \P\big(\sum_{i=1}^{n} J_i \ge nt_f-u\big) \bigg]\\
&\qquad =: r^{n} \exp \big( - ct \log t \big)\big[B_1 + B_2 + B_3\big].
\end{align*}

To bound $A_i$ and $B_i$, we claim that for every $k \in \N$ and  $u \in (0,t_f/4]$, it holds that
\begin{align}\label{tailorder}
& \P\big(\sum_{i=1}^k J_i \ge kt_f-u\big) \le \big(Kw(t_f/9)^{-1}\big)^ku^k,
\end{align}
where $K  \ge 1$ is the constant in {\bf (Trunc.)($t_f$)}(ii). Indeed, if $k=1$, then by {\bf (Trunc.)($t_f$)}(i) and (ii), we get $\P(J_1 \ge t_f-u) = F(t_f-u) = w(t_f/9)^{-1}w(t_f-u) \le Kw(t_f/9)^{-1}u$. Suppose that the claim holds for $k$. Then, by {\bf (Trunc.)($t_f$)}(i) and (ii), for all $u \in (0, t_f/4]$,
\begin{align*}
\P\big(\sum_{i=1}^{k+1} J_i \ge (k+1)t_f&-u\big) = \int_{\{\sum_{i=1}^k u_i \le u\}} F\big(t_f-u+\sum_{i=1}^k u_i\big) d_{u_{k}}F(t_f-u_{k})...d_{u_1}F(t_f-u_1) \\
& \le Kw(t_f/9)^{-1}\int_{\{\sum_{i=1}^k u_i \le u\}} \big(u-\sum_{i=1}^k u_i \big)d_{u_{k}}F(t_f-u_{k})... d_{u_1}F(t_f-u_1)\\
& \le Kw(t_f/9)^{-1}u\int_{\{\sum_{i=1}^k u_i \le u\}} d_{u_{k}}F(t_f-u_{k})... d_{u_1}F(t_f-u_1) \\
& \le Kw(t_f/9)^{-1} u\P\big(\sum_{i=1}^{k} J_i \ge kt_f-u\big)  \le \left( Kw(t_f/9)^{-1} u \right)^{k+1}.
\end{align*}
Therefore, the claim holds by induction. We consider the following two cases that when $t$ is very close to $nt_f$ and not.

\smallskip

\noindent {\it Case 1.} $(n-1/12)t_f \le t < nt_f$;

At first, by \eqref{smalljump}, we obtain $A_1 + A_2 + A_3 \le 3\P(U_r^3 \ge t_f/4) \le cr^2.$
On the other hand, by \eqref{smalljump}, \eqref{tailorder}, Proposition \ref{p:utail1}, the change of the variables and the integration by parts,
\begin{align*}
B_1+ B_2 + B_3 &\le cr^2 + c^n(nt_f-t+r)^n + cr\int_{r}^{t_f/4+t-nt_f}w(u) d_u \P\big(\sum_{i=1}^{n} J_i \ge t-u \big) \\
& \le c^nr + c^n(nt_f-t)^n + c^nr \int_{r}^{t_f/4+t-nt_f} (nt_f-t+u)^n (-dw(u))\\
& \le c^nr + c^n(nt_f-t)^n + c^n(nt_f-t)^nrw(r) +  c^nr \int_{r}^{t_f/4+t-nt_f} u^n (-dw(u))\\
& \le c^n \big(r + (nt_f-t)^n\big).
\end{align*}
In the third inequality, we used the fact that $(a+b)^k \le 2^k(a^k+b^k)$ for all $a,b>0$ and $k \in \N$ and in the fourth inequality, we used the assumption that $rw(r) \le er\phi(r^{-1}) \le 1/(4e)$. Therefore, since $n \asymp t$ so that $c^n \le ce^{ct}$, we get the result in this case.

\smallskip

\noindent {\it Case 2.} $(n-1)t_f \le t < (n-1/12)t_f$;

By \eqref{smalljump}, \eqref{tailorder}, Proposition \ref{p:utail1} and the integration by parts, we obtain
\begin{align*}
A_1 + A_2 + A_3 &\le (36c_1r/t_f)^2 + c^nr^{n-1} + cr \int_r^{t_f/4} w(u) d_u \P\big(\sum_{i=1}^{n-1} J_i \ge (n-1)t_f-u\big)  \\
& \le c^nr +  c^n r\int_r^{t_f/4} u(-dw(u)) \le c^nr.
\end{align*}
Since $B_1 + B_2 + B_3 \le 3$, $n \asymp t$ and $(nt_f-t) \asymp 1$ in this case, we finish the proof.

\qed

\begin{lemma}\label{l:trunc2}
There exists a constant $L \in (0,1)$ such that for all $t,r>0$ satisfying $t \ge t_f/2$ and $rt^{-1} \le L$, 
$$
\P(S_r \ge t) \simeq \left(\frac{r}{t}\right)^{ct} \simeq \exp \big(-ct \log \frac{t}{r}\big).
$$
\end{lemma}
\proof
Fix $r, t>0$ satisfying $t \ge t_f/2$ and $rt^{-1} \le L$ where the constant $L$ will be chosen later. Pick any $t_e \in (0, t_f)$ such that $w(t_e) \ge 1$ and let $S^*$ be the driftless subordinator with L\'evy measure $\1_{(t_e, \infty)} \cdot (-dw)$. By condition {\bf (Ker.)}, we can always find such constant $t_e$. Since $S_r \ge S^*_r$ and jump sizes of $S^*$ are at least bigger than $t_e$, by Stirling's formula, we get
\begin{align*}
\P(S_r \ge t) &\ge \P(S^*_r \ge t) \ge \P\big(S^* \; \text{jumps} \; (\lfloor t/t_e \rfloor +1) \; \text{times before time} \; r\big) \\
& = \exp \big( - rw(t_e) \big) \frac{\big(rw(t_e)\big)^{(\lfloor t/t_e \rfloor +1)}}{(\lfloor t/t_e \rfloor +1)!} \\
&\ge \exp \big( -rw(t_e) - (\lfloor t/t_e \rfloor +3/2) \log (\lfloor t/t_e \rfloor +1) + \lfloor t/t_e \rfloor + (\lfloor t/t_e \rfloor +1) \log r\big) \\
& \ge \exp \big( -c t \log \frac{t}{r} + t/(2t_e) - rw(t_e) \big) \ge \exp \big( -c t \log \frac{t}{r} + t/(2t_e) - Ltw(t_e) \big).
\end{align*}
Hence, by taking $L$ sufficiently small so that $Lw(t_e) \le 1/(2t_e)$, we get the lower bound.

On the other hand, by Chebyshev's inequality and Lemma \ref{l:laplace}, for all $\lambda > 0$, 
\begin{align*}
\P(S_r \ge t) &\le e^{-\lambda t} \E\big[e^{ \lambda S_r} \big] = \exp \big(-\lambda t + r \int_0^{t_f} (e^{\lambda u}-1) (-dw(u)) \big) \le \exp \big(-\lambda t + c_0\lambda re^{\lambda t_f}  \big),
\end{align*}
where $c_0:=\int_0^{t_f} u (-dw(u)) \in (0, \infty)$. Then, by taking $\lambda = t_f^{-1} \log \big(t/(2c_0r)\big)$, we obtain
\begin{align*}
\P(S_r \ge t) \le \exp \big( - \frac{\lambda t}{2} \big) \lesssim \exp \big( -c t \log \frac{t}{r}\big).
\end{align*}

\qed

\section{Properties of the Estimates ${\bf HK}^{\gamma, \lambda,k}_J(\Phi, \Psi)$, ${\bf HK}^{\gamma, \lambda,k}_D(\Phi)$ and ${\bf HK}^{\gamma, \lambda,k}_M(\Phi, \Psi)$}

A function $f : (0, \infty) \to \R$ is called a completely monotone fucntion if $f$ is infinitely differentiable and $(-1)^n f^{(n)}(\lambda) \ge 0$ for all $n \in \N$ and $\lambda>0$. A Bernstein function is said to be a complete Bernstein function if its L\'evy measure has a completely monotone density with respect to Lebesgue measure.

\begin{lemma}{\rm (\cite[Lemmas 3.1 and 3.2]{CKKW})}\label{l:wsp}
Assume that a family of non-negative functions $\{f(x, \cdot)\}_{x \in M}$ satisfies the weak scaling property uniformly with $(\alpha_1, \alpha_2)$ for some $0<\alpha_1 \le \alpha_2 < \infty$, that is, there are constants $c_1, c_2>0$ such that for all $x \in M$,
$$
c_1(R/r)^{\alpha_1} \le f(x, R)/f(x,r) \le c_2(R/r)^{\alpha_2}, \quad 0<r \le R < \infty.
$$
Then for any $\alpha_3>\alpha_2$, there is a family of complete Bernstein functions $\{\varphi(x, \cdot)\}_{x \in M}$ such that for all $x \in M$ and $r>0$, we have that
\begin{align*}
f(x, r) \asymp \varphi(x, r^{-\alpha_3})^{-1} \quad \text{and} \quad \partial_r \varphi(x,r) \asymp  r^{-1}\varphi(x,r).
\end{align*}
\end{lemma}

By Lemma \ref{l:wsp}, we can assume that all functions $\Phi(r), \Psi(r)$ and $V(x, r)$ are differentiable in variable $r$ and their derivatives are comparable to the function obtained by dividing $r$, i.e., $\Phi'(r) \asymp r^{-1}\Phi(r)$, $\Psi'(r) \asymp r^{-1}\Psi(r)$ and $\partial_rV(x,r) \asymp r^{-1}V(x,r)$ for all $r>0$ and $x \in M$. Indeed, for example, by Lemma \ref{l:wsp}, we have $V(x, r) \asymp \widetilde{V}(x, r):=\varphi(x, r^{-d_3})^{-1}$ for some complete Bernstein functions $\{\varphi(x, \cdot)\}_{x \in M}$ and $d_3>d_2$. Then, for all $r>0$ and $x \in M$,
\begin{align*}
r \partial_r \widetilde{V}(x,r) \asymp \frac{r^{-d_3} \varphi'(x, r^{-d_3})}{\varphi(x, r^{-d_3})^2} \asymp \widetilde{V}(x,r).
\end{align*}
Therefore, by using $\widetilde{V}$ instead of $V$, we get the desired properties.

\vspace{1mm}

Recall that for a strictly increasing function $\Phi:[0, \infty) \mapsto [0, \infty)$ which satisfies ${\bf WS}(\alpha_1, \alpha_2)$ for some $\alpha_2 \ge \alpha_1 >1$ and $\Phi(0)=0$, a function $\sM$ is determined by the relation \eqref{e:sM},
\begin{align*}
\frac{t}{\sM(t,l)} \asymp \Phi \left(\frac{l}{\sM(t,l)}\right) \qquad \text{for all} \;\; t,l>0.
\end{align*}
For example, if $\Phi(l)=l^{\alpha}$ for some $\alpha>1$, then we have $\sM(t,l) = l^{\alpha/(\alpha-1)}t^{-1/(\alpha-1)}$.

\begin{lemma}\label{l:sM}
(i) For $t, l>0$, define
\begin{align*}
\Phi_1 (t, l) := \sup_{s>0} \left\{\frac{l}{s} - \frac{t}{\Phi(s)}\right\}.
\end{align*}
Then, $\Phi_1(t,l)$ is strictly positive for all $t, l>0$, non-increasing on $(0, \infty)$ for fixed $l>0$ and satisfies \eqref{e:sM}. In other words, $\Phi_1(t,l)$ is one of the explicit forms of the function $\sM(t,l)$.\\
(ii) $\sM(\Phi(l), l) \asymp 1$ for all $l>0$. \\
(iii) There are constants $c_3, c_4>0$ such that for all $l>0$ and $0<t\le T$,
\begin{align*}
 c_3\left(\frac{T}{t}\right)^{-1/(\alpha_1-1)} \le \frac{\sM(T,l)}{\sM(t,l)} \le c_4 \left(\frac{T}{t}\right)^{-1/(\alpha_2-1)}.
\end{align*}
\end{lemma}
\proof
(i) Fix $t,l>0$ and define for $s>0$,
\begin{align*}
g(s):=\frac{l\Phi(s) - ts}{s \Phi(s)}, \quad k(s):= \frac{\Phi(s)}{s}.
\end{align*}
 We also define $k^{-1}(x):= \inf\{s: k(s) \ge x \}$ for $x>0$. Since $\Phi(s) \asymp s\Phi'(s)$ for all $s>0$, there exists a constant $c_1>0$ such that
\begin{align*}
 (s\Phi(s))^2 g'(s) = s( ts \Phi'(s) - l k(s) \Phi(s)) \ge s\Phi(s) (c_1 t - l k(s)).
\end{align*}
It follows that for $s_*:=k^{-1}(c_1t/l)$, we have $\Phi_1(t,l) = \sup_{s>0} g(s) = \sup_{s \ge s_*} g(s) \le l/s_*$. 

On the other hand, for any $a>1$, we have
\begin{align*}
\Phi_1(t,l) \ge \frac{l}{as_*} - \frac{t}{\Phi(as_*)} \ge \frac{l}{as_*} - c_2 \frac{t}{a^{\alpha_1}\Phi(s_*)} = \frac{l}{as_*} \left( 1- \frac{c_1^{-1}c_2}{a^{\alpha_1-1}} \right).
\end{align*}
Hence, by choosing $a=2 \vee (2c_1^{-1}c_2)^{1/(\alpha_1-1)}$, we get $\Phi_1(t,l) \asymp l/s_*$. Then, we conclude that
\begin{align*}
\Phi\left(\frac{l}{\Phi_1(t,l)}\right) \asymp \Phi(s_*) = s_* k(s_*) \asymp \frac{t}{\Phi_1(t,l)}.
\end{align*}
(ii), (iii) These are consequences of the relation \eqref{e:sM}.
\qed

By Lemma \ref{l:sM}(iii) and Lemma \ref{l:wsp}, we can assume that $\sM(t,l)$ is differentiable in variable $t$ for every fixed $l>0$ and there exists a constant $c_1>1$ such that for all $t, l>0$, 
\begin{align}\label{e:sMdiff}
c_1^{-1}t^{-1}\sM(t,l) \le  - \partial_t \sM(t,l) \le c_1t^{-1}\sM(t,l).
\end{align}

\smallskip

From \cite[Lemma 5.1]{CKKW}, we get the following time derivative estimates for $q(a,t,x,l;\Phi, \sM)$.

\begin{lemma}\label{l:qdiff}
For every $a>0$, there are constants $c_1, c_2>0$ such that
\begin{align*}
\left| \partial_t q^d(a,t,x,l;\Phi, \sM) \right| \le c_1 t^{-1}q^d(c_2, t, x, l;\Phi, \sM), \quad t,l>0, \;\; x \in D,
\end{align*}
Moreover, there are constants $c_3>0$ and $c_u \in (1, \infty)$ such that for all $x\in D$, 
\begin{align*}
\partial_t q^d(a,t,x,l;\Phi, \sM) \ge c_3t^{-1}q^d(a,t,x,l;\Phi, \sM) \;\;\quad \text{if} \;\; \Phi(l) \ge c_ut.
\end{align*}
\end{lemma}

\smallskip

We obtain the upper time derivative estimates for $q^j(t,x,l;\Phi, \Psi)$ and $a_k^{\gamma}(t,x,y)$.

\begin{lemma}\label{l:qdiff2}
(i) There is a constant $c_1>0$ such that for all $t,l>0$ and $x \in D$,
\begin{align*}
\left| \partial_t q^j(t,x,l;\Phi, \Psi) \right| \le c_1 t^{-1}q^j(t,x,l;\Phi, \Psi).
\end{align*} 
(ii) For all $\gamma \in [0,1)$, $t>0$, $x,y \in D$ and $j \in \{1,2\}$,
\begin{align*}
\left| \partial_t a_k^{\gamma}(t,x,y) \right| \le 2 t^{-1}a_k^{\gamma}(t,x,y).
\end{align*}
\end{lemma}
\proof
(i) Observe that
\begin{align*}
\partial_t q^j(t,x,l;\Phi, \Psi) & = \frac{\Psi(l)V(x,l) - t^2 \partial_r V(x, \Phi^{-1}(t)) \partial_t \Phi^{-1}(t)}{(tV(x,\Phi^{-1}(t))+ \Psi(l)V(x,l))^2},
\end{align*}
By
using the comparisons $\partial_r V(x, r) \asymp r^{-1}V(x,r)$ and $\partial_t \Phi^{-1}(t) \asymp t^{-1} \Phi^{-1}(t)$, we get
\begin{align*}
|\partial_t q^j(t,x,l;\Phi, \Psi)| \le \frac{ct V(x, \Phi^{-1}(t)) + \Psi(l)V(x,l) }{(tV(x,\Phi^{-1}(t))+ \Psi(l)V(x,l))^2} \le c t^{-1}q^j(t,x,l;\Phi, \Psi)
\end{align*}

(ii) From the definition of $a^\gamma_j$, we get
\begin{align*}
|\partial_t {a}_1^{\gamma}(t,x,y)| =\left( \frac{\gamma}{t+\Phi(\delta_{D}(x))} + \frac{\gamma}{t+\Phi(\delta_{D}(y))} \right){a}_1^{\gamma}(t,x,y) \le 2t^{-1}a_1^{\gamma}(t,x,y),
\end{align*}
\begin{align*}
|\partial_t {a}_2^{\gamma}(t,x,y)| &=\left( \frac{\gamma}{t+(1+t)\Phi(\delta_{D}(x))} + \frac{\gamma}{t+(1+t)\Phi(\delta_{D}(y))} \right)\frac{{a}_2^{\gamma}(t,x,y)}{1+t}\le 2t^{-1}a_2^{\gamma}(t,x,y).
\end{align*}
\qed

\section{Proof of Main Theorems}\label{s:mainproof}

In this section, we give the proof for Theorems \ref{t:mainsmall}, \ref{t:mainlarge}, \ref{t:mainsub} and \ref{t:main2}. Throughout this section, we assume that there exist $\gamma \in [0,1)$, $\lambda \ge 0$ and $k \in \{1,2\}$ such that $q(t,x,y)$ enjoys the one of the esimates  ${\bf HK}_J^{\gamma,\lambda,k}(\Phi, \Phi)$, ${\bf HK}_D^{\gamma,\lambda,k}(\Phi)$ and ${\bf HK}_M^{\gamma,\lambda,k}(\Phi, \Psi)$. Let $p(t,x,y)$ be given by \eqref{e:dfundamental}. 

\begin{prop}\label{p:ndl}{\rm \bf (Near diagonal lower bounds)} If condition {\bf (S.Poly.)($t_s$)} holds, then there exists a constant $c>0$ such that for all $(t,x,y) \in (0,t_s] \times D \times D$ satisfying $\Phi(\p(x,y)) \phi(t^{-1}) \le 1/(4e^2)$,
\begin{align}\label{ndl}
p(t,x,y) \ge c \sJ^\gamma_k(t,x,y).
\end{align}

On the other hand, if condition {\bf (L.Poly.)} holds and $\lambda =0$, then for every fixed $T >0$, \eqref{ndl} holds for all $(t,x,y) \in [T,\infty) \times D \times D$ satisfying $\Phi(\p(x,y)) \phi(t^{-1}) \le 1/(4e^2)$.
\end{prop}
\proof Since the proofs are similar, we only give the proof when condition {\bf (S.Poly.)($t_s$)} holds. Fix $(t,x,y) \in (0,t_s] \times D \times D$ satisfying $\Phi(\p(x,y)) \phi(t^{-1}) \le 1/(4e^2)$ and set $l:=\p(x,y)$. By Proposition \ref{p:utail1}, there is a constant $\varepsilon_2 \in (0,1/2]$ such that for all $t \in (0, t_s]$, we have that $\P(S_{\varepsilon_2 \Phi(l)} \ge t) \le 1/2$. Then, by the Markov property, we get
\begin{align*}
\P(S_{2\varepsilon_2\Phi(l)} \ge t) &\ge \P(S_{2\varepsilon_2\Phi(l)}-S_{\varepsilon_2\Phi(l)} \ge t \quad \text{or} \quad S_{\varepsilon_2\Phi(l)} \ge t) \\
 & \ge 1 - (1 - \P(S_{\varepsilon_2\Phi(l)}\ge t) )^2 \ge \frac{3}{2} \P(S_{\varepsilon_2\Phi(l)} \ge t).
\end{align*}
We used the inequality that $1-(1-x)^2 \ge 3x/2$ for $x \in (0,1/2]$. It follows that
\begin{align*}
\P(S_{2\varepsilon_2\Phi(l)} \ge t) -\P(S_{\varepsilon_2 \Phi(l)} \ge t) \ge \frac{1}{2} \P(S_{\varepsilon_2 \Phi(l)} \ge t).
\end{align*}
and hence by the scaling properties of $V$ and $\Phi$ and the monotonicity of $r \mapsto a_k^{\gamma}(r,x,y)$,
\begin{align}\label{e:lower2}
p(t,x,y) \ge c \int_{\varepsilon_2 \Phi(l)}^{2\varepsilon_2\Phi(l)}  \frac{a_k^{\gamma}(r,x,y)}{V(x, \Phi^{-1}(r))} d_r \P(S_r \ge t) \ge c_2  \frac{a_k^{\gamma}(\Phi(l),x,y)}{V(x,l)}\P(S_{\varepsilon_2 \Phi(l)} \ge t).
\end{align}

Besides, by the integration by parts and Lemma \ref{l:tail0},
\begin{align}\label{e:lower3}
&p(t,x,y) \ge c \int_{\varepsilon_2\Phi(l)}^{1/(2e^2\phi(t^{-1}))}  \frac{a_k^{\gamma}(r,x,y)}{V(x, \Phi^{-1}(r))} d_r \P(S_r \ge t) \nn\\
&\quad \ge  -cw(t) \int_{\varepsilon_2\Phi(l)}^{1/(2e^2\phi(t^{-1}))} rd_r\left( \frac{a_k^{\gamma}(r,x,y)}{V(x, \Phi^{-1}(r))}\right) -c_3\frac{a_k^{\gamma}(\Phi(l),x,y)}{V(x, l)}\P(S_{\varepsilon_2 \Phi(l)} \ge t)\nn \\
&\quad \ge c_4w(t) \int_{\Phi(l)}^{1/(2e^2\phi(t^{-1}))} \frac{a_k^{\gamma}(r,x,y)}{V(x,\Phi^{-1}(r))}dr  -c_3\frac{a_k^{\gamma}(\Phi(l),x,y)}{V(x,l)}\P(S_{\varepsilon_2 \Phi(l)} \ge t).
\end{align}

Finally, by Corollary \ref{c:prebound}, \eqref{e:lower2} and \eqref{e:lower3}, we deduce that
\begin{align*}
(1+c_3+c_2)p(t,x,y) &\ge c\frac{a_k^{\gamma}(1/\phi(t^{-1}),x,y)}{V\big(x,\Phi^{-1}(1/\phi(t^{-1}))\big)}+ c_2c_4w(t) \int_{\Phi(l)}^{1/(2e^2\phi(t^{-1}))} \frac{a_k^{\gamma}(r,x,y)}{V(x,\Phi^{-1}(r))}dr.
\end{align*}
\qed

In the rest of this section, we fix $(x,y) \in D \times D$ and then define $l:= \p(x,y)$ and $V(r):=V(x,r)$.

\subsection{Pure jump case}\label{s:purejump}

In this subsection, we give the proofs when $q(t,x,y)$ enjoys the estimate ${\bf HK}_J^{\gamma,\lambda,k}(\Phi, \Phi)$.

\smallskip

\noindent{\bf Proof of Theorem \ref{t:mainsmall}.}
 Fix $t \in (0,t_s]$. Since we only deal with small time $t$, we can assume that $\lambda = 0$. By \eqref{e:dfundamental} and the integration by parts, we have that for $L:=1/(4e^2)$,
\begin{align}\label{decomp}
&p(t,x,y) \asymp \int_0^{\infty}  q(r,x,y) d_r \P(S_r \ge t)\nn\\
&= \int_0^{L / \phi(t^{-1})}  q(r,x,y) d_r \P(S_r \ge t) - \int_{L / \phi(t^{-1})}^{\infty} q(r,x,y) d_r \P(S_r \le t) \nn\\
& = q\big(L / \phi(t^{-1}), x, y\big)- \int_0^{L / \phi(t^{-1})} \P(S_r \ge t) d_r q(r,x,y) +  \int_{L /\phi(t^{-1})}^{\infty} \P(S_r \le t) d_rq(r,x,y) \nn\\
& =: q\big(L / \phi(t^{-1}), x, y\big) - I_1 + I_2.
\end{align}

\vspace{2mm}

\noindent {\it Case 1.}  $\Phi(l)\phi(t^{-1}) \le 1/(4e^2)$;

By Proposition \ref{p:ndl}, it remains to prove the upper bound. We first note that
\begin{align*}
q\big(L/\phi(t^{-1}),x,y\big) \le c \frac{a_k^\gamma(1/\phi(t^{-1}),x,y)}{V\big(\Phi^{-1}(1/\phi(t^{-1}))\big)}.
\end{align*}

Next, by Proposition \ref{p:utail1}, Lemma \ref{l:qdiff2} and the definition of ${\bf HK}_J^{\gamma, \lambda, k}(\Phi,\Phi)$,
\begin{align*}
|I_1| &\le c w(t)\int_0^{L / \phi(t^{-1})} q(r,x,y) dr \\
& \le c w(t)  \int_{0}^{\Phi(l)/2}  \frac{ra_k^{\gamma}(r,x,y)}{\Phi(l)V(l)} dr+  cw(t)\int_{\Phi(l)/2}^{L / \phi(t^{-1})} \frac{a_k^{\gamma}(r,x,y) }{V(\Phi^{-1}(r))} dr =: I_{1,1} + I_{1,2}.
\end{align*}

Observe that since $\gamma<1$ and $r \mapsto r^{2\gamma} a_k^{\gamma}(r,x,y)$ is increasing, we have that
\begin{align*}
I_{1,1} & \le cw(t) \int_0^{\Phi(l)/2} \frac{r^{2\gamma}a_k^{\gamma}(r,x,y) r^{1-2\gamma}}{\Phi(l)V(l)} dr \le cw(t) \frac{\Phi(l)^{2\gamma}a_k^{\gamma}(\Phi(l),x,y)}{\Phi(l)V(l)} \int_0^{\Phi(l)/2} r^{1-2\gamma} dr  \\
&\le cw(t)a_k^{\gamma}(\Phi(l),x,y) \frac{\Phi(l)}{V(l)} \le cw(t) \int_{\Phi(l)/2}^{\Phi(l)}\frac{a_k^{\gamma}(r,x,y) }{V(\Phi^{-1}(r))} dr \le I_{1,2}.
\end{align*}
Therefore, we have that $|I_1| \le c I_{1,2} \le c \sJ^\gamma_k(t,x,y)$.

Lastly, by Lemma \ref{l:qdiff2} and the change of variables, we get
\begin{align*}
|I_2| & \le c \int_{L / \phi(t^{-1})}^{\infty}  \frac{a_k^{\gamma}(r,x,y) }{rV(\Phi^{-1}(r))} dr \le ca_k^{\gamma}(1/\phi(t^{-1}),x,y)  \int_{L}^{\infty} \frac{1}{sV\big(\Phi^{-1}(s/\phi(t^{-1}))\big)} ds \\
& \le c\frac{a_k^{\gamma}(1/\phi(t^{-1}),x,y) }{V\big(\Phi^{-1}(1/\phi(t^{-1}))\big)}\int_{L}^{\infty} s^{-1-d_1/\alpha_2} ds  = c\frac{a_k^{\gamma}(1/\phi(t^{-1}),x,y) }{V\big(\Phi^{-1}(1/\phi(t^{-1}))\big)}.
\end{align*}

Therefore, we obtain the upper bound from \eqref{decomp}.

\smallskip

\noindent {\it Case 2.}  $\Phi(l)\phi(t^{-1})>1/(4e^2)$;

In this case, we have
\begin{align*}
q\big(L/\phi(t^{-1}), x, y\big) \le c \frac{ a_k^{\gamma}(1/\phi(t^{-1}),x,y)}{\phi(t^{-1})\Phi(l)V(l)}.
\end{align*}

By Proposition \ref{p:utail1}, Lemma \ref{l:qdiff2} and the fact that $\phi(t^{-1}) \ge e^{-1}w(t)$,
\begin{align*}
|I_1| & \le c w(t)  \int_{0}^{L/\phi(t^{-1})}  \frac{ra_k^{\gamma}(r,x,y)}{\Phi(l)V(l)} dr \le cw(t) \frac{a_k^{\gamma}(1/\phi(t^{-1}),x,y)}{\phi(t^{-1})^{2\gamma}} \int_{0}^{L/\phi(t^{-1})}  \frac{r^{1-2\gamma}}{\Phi(l)V(l)} dr \\
&\le cw(t) a_k^{\gamma}(1/\phi(t^{-1}),x,y) \frac{\phi(t^{-1})^{-2}}{\Phi(l)V(l)} \le c \frac{ a_k^{\gamma}(1/\phi(t^{-1}),x,y)}{\phi(t^{-1})\Phi(l)V(l)}.
\end{align*}

Moreover, by Lemma \ref{l:qdiff2}, Lemma \ref{l:Offprob}, Lemma \ref{l:b}(ii) and the change of variables, 
\begin{align*}
|I_2| & \le c\int_{L/\phi(t^{-1})}^{b^{-1}(t)}  \frac{a_k^{\gamma}(r,x,y)}{V(l)\Phi(l)} dr+ c\int_{ b^{-1}(t)}^{\infty} \frac{a_k^{\gamma}(r,x,y)\exp\big(-r H(\phi'^{-1}(t/r))\big) }{V(l)\Phi(l)} dr \\
& \le c \frac{a_k^{\gamma}(1/\phi(t^{-1}),x,y)}{\phi(t^{-1})V(l)\Phi(l)} \left[1 + \int_{1}^{\infty} \exp\bigg(-b^{-1}(t)s H\big(\phi'^{-1}(t/(b^{-1}(t)s))\big)\bigg)ds\right] \\
& \le c \frac{a_k^{\gamma}(1/\phi(t^{-1}),x,y)}{\phi(t^{-1})V(l)\Phi(l)}\left[1+ \int_1^{\infty} \exp(-s) ds \right] \le c \frac{a_k^{\gamma}(1/\phi(t^{-1}),x,y)}{\phi(t^{-1})V(l)\Phi(l)}.
\end{align*}
In the third inequality, we used the fact that $s \mapsto H\big(\phi'^{-1}(t/(b^{-1}(t)s))\big)$ is increasing and $b^{-1}(t)H\big(\phi'^{-1}(t/(b^{-1}(t)))\big) = 1$.
This proves the upper bound.

On the other hand, by Corollary \ref{c:prebound} and the definition of ${\bf HK}_J^{\gamma, \lambda, k}(\Phi, \Phi)$, we obtain
\begin{align*}
p(t,x,y) \ge c \frac{a_k^{\gamma}(1/\phi(t^{-1}),x,y)}{\phi(t^{-1})V(l)\Phi(l)}.
\end{align*}
\qed

\noindent{\bf Proof of Theorem \ref{t:mainlarge}.}
If $\lambda = 0$, then by using Proposition \ref{p:utail2} instead of Proposition \ref{p:utail1}, the proof is essentially the same as the one for Theorem \ref{t:mainsmall}. Hence, we omit it in here. Now, assume that $\lambda>0$ and $R_D=\text{diam}(D)< \infty$. Let $T_*:=1/(4e^2 \phi(T^{-1}))$. Then, by Proposition \ref{p:utail2}, Lemma \ref{l:qdiff2} and the integration by parts,
\begin{align*}
&p(t,x,y)  \asymp \int_0^{T_*} q(r,x,y) d_r \P(S_r \ge t)  + \Phi(\delta_{D}(x))^{\gamma} \Phi(\delta_{D}(y))^{\gamma} \int_{T_*}^{\infty} e^{-\lambda r} d_r \P(S_r \ge t) \\
&\le  q(T_*,x,y) \P(S_{T_*} \ge t)  + c w(t) \int_0^{T_*} q(r,x,y) dr  + \lambda \Phi(\delta_{D}(x))^{\gamma} \Phi(\delta_{D}(y))^{\gamma} \int_{T_*}^{\infty} e^{-\lambda r} \P(S_r \ge t) dr \\
&=:  q(T_*,x,y) \P(S_{T_*} \ge t)  + J_1 + J_2 \le cw(t)\Phi(\delta_{D}(x))^{\gamma} \Phi(\delta_{D}(y))^{\gamma} + J_1 + J_2.
\end{align*}

By Proposition \ref{p:utail2} and Lemma \ref{l:decaycomp}, we obtain
\begin{align*}
J_2 & \le c\Phi(\delta_{D}(x))^{\gamma} \Phi(\delta_{D}(y))^{\gamma}\bigg(w(t)\int_{T_*}^{1/(4e^2 \phi(t^{-1}))} re^{-\lambda r}  dr + \int_{1/(4e^2 \phi(t^{-1}))}^{\infty} e^{-\lambda r} dr \bigg)\\
& \le \Phi(\delta_{D}(x))^{\gamma} \Phi(\delta_{D}(y))^{\gamma}   \big(w(t) + \exp\big(-c/\phi(t^{-1})\big)\big)\\
& \le c\Phi(\delta_{D}(x))^{\gamma} \Phi(\delta_{D}(y))^{\gamma}   \big(w(t) + \phi(t^{-1})^{\delta_2 +1}\big) \le cw(t)\Phi(\delta_{D}(x))^{\gamma} \Phi(\delta_{D}(y))^{\gamma}.
\end{align*}
In the third inequality, we used the fact that for every $\delta>0$, $e^{-1/x} \le \delta^\delta e^{-\delta}x^{\delta}$ for all $x>0$.

On the other hand, we note that 
\begin{align*}
2\int_{T_*\Phi(l)/(2\Phi(R_D))}^{T_*} \frac{a_1^{\gamma}(r,x,y) }{V(\Phi^{-1}(r))} dr &\ge  \int_{T_*\Phi(l)/(2\Phi(R_D))}^{T_*\Phi(l)/\Phi(R_D)} \frac{ra_1^{\gamma}(r,x,y)}{V(l)\Phi(l)} dr + \int_{T_*/2}^{T_*} \frac{a_1^{\gamma}(r,x,y) }{V(\Phi^{-1}(r))} dr \\
& \ge c\frac{a_1^{\gamma}(\Phi(l),x,y)\Phi(l)}{V(l)} + c \Phi(\delta_{D}(x))^{\gamma} \Phi(\delta_{D}(y))^{\gamma}\\
& \ge c\int_{0}^{T_*\Phi(l)/(2\Phi(R_D))} \frac{ra_1^{\gamma}(r,x,y)}{V(l)\Phi(l)} dr + c \Phi(\delta_{D}(x))^{\gamma} \Phi(\delta_{D}(y))^{\gamma}.
\end{align*}
Thus, by the scaling properties of $a_1^\gamma, V$ and $\Phi$, we get
\begin{align*}
J_1 &\asymp w(t) \int_{T_*\Phi(l)/(2\Phi(R_D))}^{T_*} \frac{a_1^{\gamma}(r,x,y) }{V(\Phi^{-1}(r))} dr \asymp w(t) \int_{\Phi(l)}^{2\Phi(R_D)} \frac{a_1^{\gamma}(r,x,y) }{V(\Phi^{-1}(r))} dr \\
&\ge  cw(t) \Phi(\delta_{D}(x))^{\gamma} \Phi(\delta_{D}(y))^{\gamma} \ge cJ_2.
\end{align*}
This proves the upper bound.

On the other hand, by essentially the same proof as the one for Proposition \ref{p:ndl}, we get the lower bound.  We omit the details in here. \qed

\smallskip

\noindent{\bf Proof of Theorem \ref{t:mainsub}.}
If $\lambda = 0$, then by using Proposition \ref{p:utail3} instead of Proposition \ref{p:utail1} and the fact that $\phi(t^{-1}) \asymp t^{-1}$ for all $t \ge T$ which follows from Lemma \ref{l:phiHw}, we get the desired results. Hence, we assume that $\lambda>0$ and $R_D=\text{diam}(D)< \infty$. Let $L>0$ be the minimum of the constants in Propositions \ref{p:utail3} and \ref{p:utail4}. By the integration by parts, Proposition \ref{p:utail4} with $k = \lambda/2$ and the argument given in the proof of Theorem \ref{t:mainlarge},
\begin{align*}
p(t,x,y)& \le c \int_0^{LT} q(r,x,y) d_r \P(S_r \ge t) + c\Phi(\delta_D(x))^\gamma\Phi(\delta_D(y))^\gamma \int_{LT}^{\infty} e^{-\lambda r} d_r \P(S_r \ge t)\\
&\le c\bigg[\1_{\{0<\beta<1\}}\exp\big(-c_1 t^{\beta}\big) + \1_{\{\beta=1\}}\exp\big(-\frac{c_1}{2}t\big) \bigg] \bigg(q(LT,x,y) + \int_{\Phi(l)}^{2\Phi(R_D)} \frac{a_1^\gamma(r,x,y)}{V(\Phi^{-1}(r))} dr\bigg)\\
& \quad + c \lambda \Phi(\delta_D(x))^\gamma\Phi(\delta_D(y))^\gamma \bigg[\1_{\{0<\beta<1\}}\exp\big(-c_1 t^{\beta}\big) + \1_{\{\beta=1\}}\exp\big(-\frac{c_1}{2}t\big) \bigg] \int_{LT}^{Lt} re^{-\lambda r/2} dr  \\
& \quad + c \lambda \Phi(\delta_D(x))^\gamma\Phi(\delta_D(y))^\gamma \int_{Lt}^{\infty}  e^{-\lambda r}\P(S_r \ge t) dr \\
& \le c\bigg[\1_{\{0<\beta<1\}}\exp\big(-c_1 t^{\beta}\big) + \1_{\{\beta=1\}}\exp\big(-\frac{c_1}{2}t\big) \bigg] \int_{\Phi(l)}^{2\Phi(R_D)} \frac{a_1^\gamma(r,x,y)}{V(\Phi^{-1}(r))} dr \\
& \quad + c \lambda e^{-\lambda Lt} \Phi(\delta_D(x))^\gamma\Phi(\delta_D(y))^\gamma .
\end{align*}
This proves the upper bound.

On the other hand, by the proof for Proposition \ref{p:ndl}, we can obtain that
\begin{align*}
p(t,x,y) \ge c w(t)\int_{\Phi(l)}^{2\Phi(R_D)} \frac{a_1^\gamma(r,x,y)}{V(\Phi^{-1}(r))} dr.
\end{align*}
Furthermore, by Corollary \ref{c:meandist} and the fact that $\phi(t^{-1})^{-1} \asymp t$ for all $t \ge T$, there exists a constant $L_1>0$ such that
\begin{align*}
p(t,x,y) \ge c\inf_{r \in (0, L_1t)}q(r,x,y) \ge ce^{-\lambda L_1t}\Phi(\delta_D(x))^\gamma\Phi(\delta_D(y))^\gamma.
\end{align*}
Hence, we get the lower bound.
\qed

\subsection{Diffusion case}\label{s:diffusion}

In this subsection, we provide the proof when $q(t,x,y)$ enjoys the estimate ${\bf HK}_D^{\gamma,\lambda, k}(\Phi)$. Set $k(c_0,r):= a_k^{\gamma}(r,x,y)q(c_0,r,x,l;\Phi, \sM)$ for $c_0>0$ and $r>0$ where the function $\sM$ is determined by the relation \eqref{e:sM}.

\smallskip

\noindent{\bf Proof of Theorem \ref{t:mainsmall}.}
Since we only consider small time $t$, we can assume that $\lambda = 0$. For every fixed $t \in (0,t_s]$, by the integration by parts, we have that for $L:=1/(4e^2)$, 
\begin{align*}
&p(t,x,y) \simeq \int_0^{\infty} k(c, r) d_r \P(S_r \ge t) \\
& = k(c, L/\phi(t^{-1})) - \int_0^{ L/\phi(t^{-1})} \P(S_r \ge t) d_r  k(c, r) +  \int_{ L/\phi(t^{-1})}^{\infty} \P(S_r \le t) d_r  k(c, r)\\
& =: k(c,  L/\phi(t^{-1})) - I_1 + I_2.
\end{align*}

\smallskip

\noindent {\it Case 1.}  $\Phi(l)\phi(t^{-1})\le 1/(4e^2)$;

Note that by a similar proof as the one given in section \ref{s:purejump}, we obtain
\begin{align*}
I_2 & \le ck\big(c_1,  1/(4e^2\phi(t^{-1}))\big) \le \frac{a_k^{\gamma}(1/\phi(t^{-1}),x,y)}{V\big(\Phi^{-1}(1/\phi(t^{-1}))\big)}.
\end{align*}
Hence, by Propostion \ref{p:ndl}, it remains to get upper bound for $I_1$.

By Lemma \ref{l:qdiff}, Proposition \ref{p:utail1}, the change of variables and Lemma \ref{l:sM}(iii),
\begin{align*}
|I_1| & \le cw(t)\int_0^{\Phi(l)/2} \frac{a_k^{\gamma}(r,x,y)}{V(\Phi^{-1}(r))}e^{-c_2\sM(r,l)}dr + cw(t)\int_{\Phi(l)/2}^{ L/\phi(t^{-1})} \frac{a_k^{\gamma}(r,x,y)}{V(\Phi^{-1}(r))} dr \\
& \le cw(t)\Phi(l)^{1-2\gamma}\int_0^{1/2} \frac{(\Phi(l)s)^{2\gamma}a_1^{\gamma}(\Phi(l)s,x,y)}{s^{2\gamma}V(\Phi^{-1}(\Phi(l)s))}e^{-c_2 \sM(\Phi(l)s,l)}ds + c\sJ^\gamma_k(t,x,y)\\
&\le c w(t) \frac{\Phi(l) a_1^{\gamma}(\Phi(l)/2, x, y)}{V(l)} \int_0^{1/2} s^{-d_2/\alpha_1-2\gamma} \exp\big(-c_3 s^{-1/(\alpha_2-1)}\big)ds  + c\sJ^\gamma_k(t,x,y)\\
&\le cw(t) a_k^{\gamma}(\Phi(l), x, y) \frac{\Phi(l)}{V(l)} + c\sJ^\gamma_k(t,x,y) \le  cw(t)\int_{\Phi(l)/2}^{\Phi(l)} \frac{a_k^{\gamma}(r,x,y)}{V(\Phi^{-1}(r))} dr + c\sJ^\gamma_k(t,x,y) \\
& \le c\sJ^\gamma_k(t,x,y).
\end{align*}

\smallskip

\noindent {\it Case 2.}  $\Phi(l)\phi(t^{-1}) > 1/(4e^2)$;

Define for every $a>0$ and $r>0$, 
$\displaystyle{
g(a,r):=\frac{\exp\big(-a \sM(r,l)\big)}{r^{1+2\gamma}V(\Phi^{-1}(r))}.
}$
 Then, we see that
\begin{align*}
\frac{dg(a,r)}{dr} &=\left(-ar \partial_r \sM(r,l) -(1+2\gamma) -  r\partial_r V(\Phi^{-1}(r)) \cdot V(\Phi^{-1}(r))^{-1} \right)  \frac{\exp\big(-a \sM(r,l)\big)}{r^{2+2\gamma}V(\Phi^{-1}(r))} \\
&\ge \big( ac_4 \sM(r,l) - c_5 \big)\frac{\exp\big(-a \sM(r,l)\big)}{r^{2+2\gamma}V(\Phi^{-1}(r))} ,
\end{align*}
for some positive constants $c_4$ and $c_5$ independent of $a$ and $r$. By Lemma \ref{l:sM}(ii) and (iii), for each fixed $a>0$, there exists a constant $\delta>0$ such that $g(a,r)$ is increasing on $0<r< \delta \Phi(l)$.
By Lemma \ref{l:qdiff} and the fact that $r \mapsto r^{2\gamma} a_k^\gamma(r,x,y)$ is increasing on $r>0$, we get
\begin{align*}
|I_1| & \le c \int_0^{L/\phi(t^{-1})} r^{-1}k(c_6,r) dr = c \int_0^{L/\phi(t^{-1})} r^{2\gamma}a_k^{\gamma}(r,x,y)g(c_6, r) dr \\
& \le c\phi(t^{-1})^{-(1+2\gamma)} a_k^{\gamma}(1/\phi(t^{-1}),x,y) \sup_{0 <r < \phi(t^{-1})^{-1}} g(c_6, r).
\end{align*}
Therefore, if $\phi(t^{-1})^{-1}< \delta(c_6)\Phi(l)$, then we get 
\begin{align*}
|I_1| &\le ca_k^{\gamma}(1/\phi(t^{-1}),x,y) \frac{\exp\big(- c_7 \sM(1/\phi(t^{-1}),l)\big)}{V\big(\Phi^{-1}(1/\phi(t^{-1}))\big)}.
\end{align*}
Otherwise, if $\phi(t^{-1})^{-1} \ge \delta(c_6)\Phi(l)$, then $\phi(t^{-1})^{-1} \asymp \Phi(l)$ and hence by Lemma \ref{l:sM}(iii),
\begin{align*}
|I_1| &\le c\phi(t^{-1})^{-(1+2\gamma)} a_k^{\gamma}(1/\phi(t^{-1}),x,y) \sup_{ \delta(c_6)\Phi(l) < r<  \phi(t^{-1})^{-1}} g(c_6, r)\\
&\le ca_k^{\gamma}(1/\phi(t^{-1}),x,y) \frac{\exp\big(- c_9 \sM(1/\phi(t^{-1}),l)\big)}{V\big(\Phi^{-1}(1/\phi(t^{-1}))\big)}.
\end{align*}

Next, by Lemma \ref{l:qdiff}, Lemma \ref{l:Offprob} and Lemma \ref{l:b}(ii), we have
\begin{align*}
|I_2| & \le c \int_{L/\phi(t^{-1})}^{b^{-1}(t)}  r^{-1}k(c_6,r) dr+ c\int_{b^{-1}(t)}^{\infty} r^{-1}k(c_6,r) \exp\big(-r H(\phi'^{-1}(t/r))\big)dr  =: I_{2,1} + I_{2,2}.
\end{align*}

By Lemma \ref{l:b}(ii) and Lemma \ref{l:sM}(iii), we have 
\begin{align*}
I_{2,1} \le  ca_k^{\gamma}(1/\phi(t^{-1}),x,y) \frac{\exp\big(- c_{10} \sM(1/\phi(t^{-1}),l)\big)}{V\big(\Phi^{-1}(1/\phi(t^{-1}))\big)}.
\end{align*}

To control the exponential terms in $I_{2,2}$, we consider the following two functions that $e_1(r):=rH(\phi'^{-1}(t/r))$ and $e_2(r):=\sM(r,l)$. (cf.~ \cite{CKKW}.) Note that $e_1$ is non-decreasing and $e_2$ is non-increasing. Moreover, by the definition of the function $b$, $e_1(b^{-1}(t)) = 1$ for all $t>0$ and $e_1(\infty) = \infty$ and by Lemma \ref{l:sM}(ii) and (iii), $e_2(\Phi(l)) \asymp 1$ for all $l>0$ and $e_2(\infty)=0$. Thus, by the intermediate value theorem, there are constants $a_1>0$ and $a_2>0$ independent of $t$ and $l$ such that for all $t,l>0$ with $\Phi(l)\phi(t^{-1}) > 1/(4e^2)$, there exists a unique $r^*=r^*(t,l) \in (b^{-1}(t), a_1 \Phi(l))$ such that $e_1(r^*) = a_2 e_2(r^*)$. Now, we have
\begin{align*}
 a_k^{\gamma}(1/\phi(t^{-1}),x,y)^{-1}I_{2,2} &\le c\int_{b^{-1}(t)}^{r^*}\frac{\exp\big(-c_6 \sM(r,l)\big)}{rV(\Phi^{-1}(r))} dr  + c\int_{r^*}^{\infty}\frac{\exp\big(-r H(\phi'^{-1}(t/r))\big)}{rV(\Phi^{-1}(r))} dr \\
& =: I_{2,2,1} + I_{2,2,2}.
\end{align*}
By the change of variables and Lemma \ref{l:sM}(ii) and (iii), we get
\begin{align*}
I_{2,2,1} & = c\int_{b^{-1}(t)/r^*}^1 \frac{\exp\big(-c_6 \sM(r^*s,l)\big)}{sV(\Phi^{-1}(r^*s))} ds \\
& \le c \frac{ \exp\big(-\frac{c_6}{2}\sM(r^*,l)\big)}{V(\Phi^{-1}(r^*))}\int_0^1 s^{-1-d_2/\alpha_1} \exp \big(-cs^{-1/(\alpha_2-1)}\big) ds \le c \frac{ \exp\big(-\frac{c_6}{2}\sM(r^*,l)\big)}{V\big(\Phi^{-1}(1/\phi(t^{-1}))\big)}.
\end{align*}
Also, by the change of variables, we have
\begin{align*}
I_{2,2,2} & = c\int_{r^*/b^{-1}(t)}^{\infty} \frac{\exp \big( - b^{-1}(t)s H\big( \phi'^{-1}(t/(b^{-1}(t)s)) \big) \big)}{sV(\Phi^{-1}(b^{-1}(t)s))} ds \\
& \le c\frac{\exp\big(-e_1(r^*)\big)}{V(\Phi^{-1}(b^{-1}(t)))} \int_{1}^{\infty} s^{-1-d_1/\alpha_2} ds \le c \frac{ \exp\big(-a_2\sM(r^*,l)\big)}{V\big(\Phi^{-1}(1/\phi(t^{-1}))\big)}.
\end{align*}

To determine the function $\sM(r^*,l)$, we note that by \eqref{e:sM}, $e_1(r^*) \asymp e_2(r^*)$ implies that
\begin{align*}
\frac{r^*}{r^*H(\phi'^{-1}(t/r^*))} \asymp \Phi\left(\frac{l}{r^*H(\phi'^{-1}(t/r^*))}\right).
\end{align*}
Let $s^*=1/H(\phi'^{-1}(t/r^*))$. Then, $b(s^*)/s^* = \phi'(H^{-1}(1/s^*)) = t/r^*$. Therefore, by Lemma \ref{l:b}(ii), the function $\sN(t,l):=a_2\sM(r^*,l) = e_1(r^*)=r^*/s^*$ is determined by the relation
\begin{align*}
\frac{1}{\phi\big(\sN(t,l)/t\big)} \asymp b^{-1}\left(\frac{t}{\sN(t,l)}\right) = s^* \asymp \Phi\left(\frac{l}{\sN(t,l)}\right).
\end{align*}
Since $\sM(b^{-1}(t),l) \ge ce_1(r^*)$, we finish the proof for the upper bound.

\smallskip

Now, we prove the lower bound. By Lemma \ref{l:qdiff} and the integration by parts, we have
\begin{align}\label{e:offl}
&p(t,x,y) \ge -c \int_{Nb^{-1}(t)}^{\infty} k(c_8, r) d_r\P(S_r\le t) \nn \\
&\;\ge c \int_{Nb^{-1}(t)}^{\Phi(l)/c_u} r^{-1}\P(S_r \le t) k(c_8, r) dr - c \int_{\Phi(l)/c_u}^{\infty} r^{-1} \P(S_r \le t) k(c_9,r) dr \nn \\
&\; \ge c_{10} \int_{Nb^{-1}(t)}^{\infty} r^{-1}\P(S_r \le t) k(c_8, r) dr - c_{11} \int_{\Phi(l)/c_u}^{\infty} r^{-1} \P(S_r \le t) k(c_9, r) dr  :=J_1-J_2,
\end{align}
where $N:=(e-2)/(8c_ue^2(e^2-e))$. Note that by Lemma \ref{l:b}(ii), we have that $Nb^{-1}(t) \le 1/(8e^2c_u\phi(t^{-1})) \le \Phi(l)/(2c_u)$. By taking $c_u$ large enough, we may assume that $N \in (0,1/2)$. 
Then, by Lemma \ref{l:Offprob} and Lemma \ref{l:b}(ii),
\begin{align*}
& J_1 \ge c a_k^{\gamma}(1/\phi(t^{-1}),x,y) \phi(t^{-1})^{-2\gamma}\int_{Nb^{-1}(t)}^{\infty} \frac{\exp \big(-c_{12} rH (\phi'^{-1}(t/r)) -c_8\sM(r,l) \big)}{r^{1+2\gamma}V(\Phi^{-1}(r))}dr.
\end{align*}
Let $e_3(r)=c_{12}rH(\phi'^{-1}(t/r))$ and $e_4(r)=c_8\sM(r,l)$ for $r>0$. By the same argument as in the proof for the upper bounds, there are constants $a_3, a_4>0$ independent of $t$ and $l$ such that for all $t,l>0$ with $\Phi(l)\phi(t^{-1}) > 1/(4e^2)$, there exists a unique $r_*=r_*(t,l) \in (b^{-1}(t), a_3\Phi(l))$ such that $e_3(r_*)=a_4e_4(r_*).$ Moreover, from the monotonicity,
$$
e_3(r)<a_4e_4(r) \quad \text{for} \;\; r \in (b^{-1}(t),r_*) \quad \text{and} \quad e_3(r)>a_4e_4(r) \quad \text{for} \;\; r>r_*.
$$
Therefore,  by the change of variables, Lemma \ref{l:b}(ii) and the weak scaling properties,
\begin{align*}
& \phi(t^{-1})^{-2\gamma} \int_{Nb^{-1}(t)}^{\infty} \frac{\exp \big(-e_3(r) - e_4(r)\big)}{r^{1+2\gamma}V(\Phi^{-1}(r))}dr  =   \int_{N}^{\infty} \frac{\exp\big(-e_3(b^{-1}(t)s)-e_4(b^{-1}(t)s)\big)}{s^{1+2\gamma}V\big(\Phi^{-1}(b^{-1}(t)s)\big)}ds \\
& \;\; \ge \frac{c}{V\big(\Phi^{-1}(\phi(t^{-1})^{-1})\big)}\int_{r_*/(2b^{-1}(t))}^{r_*/b^{-1}(t)} s^{-1-d_2/\alpha_1-2\gamma}\exp \big(-(1+a_4)e_4(b^{-1}(t)s)\big) ds \\
& \;\; \ge \frac{c}{V\big(\Phi^{-1}(\phi(t^{-1})^{-1})\big)}  (r_*/b^{-1}(t))^{-d_2/\alpha_1-2\gamma}\exp \big(-c_{12}e_4(r_*)\big) \\
& \;\; \ge c \frac{\exp \big(-2c_{12}e_4(r_*)\big)}{V\big(\Phi^{-1}(\phi(t^{-1})^{-1})\big)}  (r_*/b^{-1}(t))^{-d_2/\alpha_1-2\gamma}\exp \big(\frac{c_{12}}{a_4}e_3(r_*)\big) \\
& \;\; \ge c \frac{\exp \big(-2c_{12}e_4(r_*)\big)}{V\big(\Phi^{-1}(\phi(t^{-1})^{-1})\big)}  (r_*/b^{-1}(t))^{-d_2/\alpha_1-2\gamma}\exp \big(\frac{c_{12}r_*}{a_4b^{-1}(t)}e_3(b^{-1}(t))\big)  \ge c \frac{\exp \big(-2c_{12}e_4(r_*)\big)}{V\big(\Phi^{-1}(\phi(t^{-1})^{-1})\big)}.
\end{align*}
In the last inequality, we used the fact that $e_3(b^{-1}(t)) = c_{12}$ and that for every $p>0$, there exists a constant $c(p)>0$ such that $e^x \ge c(p)x^p$ for all $x>0$. It follows that
\begin{align*}
J_1 \ge c_{13}a_k^{\gamma}(1/\phi(t^{-1}),x,y) \frac{ \exp \big( - c_{14} e_3(r_*)\big)}{V\big(\Phi^{-1}(\phi(t^{-1})^{-1})\big)},
\end{align*}
for some constants $c_{13}, c_{14}>0$.

On the other hand, by Lemma \ref{l:Offprob}, we have that
\begin{align*}
J_2 & \le c a_k^{\gamma}(\Phi(l),x,y)\int_{\Phi(l)}^{\infty} \frac{\exp \big(-rH(\phi'^{-1}(t/r))\big)}{rV(\Phi^{-1}(r))}dr \le ca_k^{\gamma}(1/\phi(t^{-1}),x,y)\frac{\exp \big(-c_{15} e_3(\Phi(l))\big)}{V(l)}.
\end{align*}
Since $e_3(Ar) \ge Ae_3(r)$ for all $r>0$ and $A \ge 1$, from \eqref{e:offl}, we deduce that there exists a constant $A>0$ such that $\Phi(l)>A r_* $ implies that
\begin{align*}
p(t,x,y) \ge ca_k^{\gamma}(1/\phi(t^{-1}),x,y)\frac{\exp \big( - c_{14} e_3(r_*)\big)}{V\big(\Phi^{-1}(\phi(t^{-1})^{-1})\big)},
\end{align*}
which yields the result. Otherwise, if $\Phi(l) \le Ar_*$, then by Lemma \ref{l:sM}(ii) and (iii),
$$
e_3(b^{-1}(t)) =c_{12} \ge c e_4(\Phi(l)) \ge c e_4(r_*) = ca_4^{-1} e_3(r_*) \ge ca_4^{-1} (r_*/b^{-1}(t))e_3(b^{-1}(t)).
$$
It follows that $b^{-1}(t) \asymp r_* \asymp \Phi(l)$ in this case. Since by Corollary \ref{c:prebound}, we have that
\begin{align*}
p(t,x,y) \ge ca_k^{\gamma}(1/\phi(t^{-1}),x,y)\frac{1}{V\big(\Phi^{-1}(\phi(t^{-1})^{-1})\big)},
\end{align*}
we still get the result in this case.
\qed

\noindent{\bf Proof of Theorems \ref{t:mainlarge} and \ref{t:mainsub}.} Observe that both ${\bf HK}^{\gamma, \lambda, k}_J(\Phi,\Phi)$ and ${\bf HK}_D^{\gamma, \lambda, k}(\Phi)$ give the same estimates for $q(t,x,y)$ on near diagonal situation, that is, when $t \ge c\p(x,y)$ for some constant $c>0$. Using this fact, we deduce the result by the same argument given in section \ref{s:purejump}.
\qed

\subsection{Mixed type case}

In this subsection, we give the proof when $q(t,x,y)$ enjoys the estimate ${\bf HK}_M^{\gamma,\lambda,k}(\Phi, \Psi)$. Since the ideas for proofs are similar, we only provide the proof of Theorem \ref{t:mainsmall}. This completes the proof for Theorems \ref{t:mainsmall}, \ref{t:mainlarge} and \ref{t:mainsub}.

\smallskip

\noindent{\bf Proof of Theorem \ref{t:mainsmall}.}
 Fix $t \in (0,t_s]$. Define for $r>0$ and $c_0>0$,
\begin{align*}
m_1(r):= a_k^{\gamma}(r,x,y)q^j(r,x,l;\Phi, \Psi), \qquad m_2(c_0, r):= a_k^{\gamma}(r,x,y)q^d(c_0,r,x,l;\Phi, \sM).
\end{align*}
We also define for $t>0$ and $c_0>0$,
\begin{align*}
p_1(t):= \int_0^{\infty} m_1(r) d_r \P(S_r \ge t), \quad p_2(c_0,t):= \int_0^{\infty} m_2(c_0, r) d_r \P(S_r \ge t).
\end{align*}
Then, from the definition, we get
\begin{align}\label{e:mixedsoln}
p(t,x,y) \simeq p_1(t) + p_2(c, t).
\end{align}

\smallskip

\noindent {\it Case 1.}  $\Phi(l)\phi(t^{-1}) \le 1/(4e^2)$;

By the proof given in section \ref{s:diffusion}, for each fixed $c_0>0$, $p_2(c_0, t) \asymp \sJ^\gamma_1( t, x, y)$.
On the other hand, since $\Psi(l) \ge \Phi(l)$ for all $l>0$, by the proof given in section \ref{s:purejump}, $p_1(t) \le c\int_0^{\infty} a_1^{\gamma}(r,x,y)q^j(r,x,l;\Phi,\Phi) d_r\P(S_r \ge t) \le c \sJ^\gamma_1( t, x, y)$. Therefore, \eqref{e:mixedsoln} yields the result.

\smallskip

\noindent {\it Case 2.}  $\Phi(l)\phi(t^{-1}) > 1/(4e^2)$;

By the proof given in section \ref{s:diffusion}, we get 
\begin{align*}
&p_2(c, t) \simeq a_k^{\gamma}(1/\phi(t^{-1}),x,y)\frac{\exp\big(-c \sN(t, \p(x,y))\big)}{V\big( \Phi^{-1}(1/\phi(t^{-1}))\big)}.
\end{align*}

On the other hand, by Lemma \ref{l:qdiff2}, the integration by parts, Proposition \ref{p:utail1} and Lemma \ref{l:Offprob},
\begin{align*}
p_1(t) &= m_1\big(1/(4e^2 \phi(t^{-1}))\big) - \int_0^{1/(4e^2 \phi(t^{-1}))} \P(S_r \ge t) d_r m_1(r) + \int_{1/(4e^2 \phi(t^{-1}))}^{\infty} \P(S_r \le t) d_r m_1(r) \\
& \le c  \frac{a_k^{\gamma}(1/\phi(t^{-1}),x,y)}{\phi(t^{-1})V(l)\Psi(l)}  + c w(t) \int_0^{1/(4e^2 \phi(t^{-1}))} \frac{r^{2\gamma}a_k^{\gamma}(r,x,y)}{V(l)\Psi(l)}r^{1-2\gamma} dr \\
&\quad+ c  a_k^{\gamma}(1/\phi(t^{-1}),x,y) \int_{1/(4e^2 \phi(t^{-1}))}^{\infty} \frac{\exp \big(-r H(\phi'^{-1}(t/r))\big) }{V(l)\Psi(l)} dr  \le c\frac{ a_k^{\gamma}(1/\phi(t^{-1}),x,y)}{\phi(t^{-1})V(l)\Psi(l)}.
\end{align*}
We also have that by Corollary \ref{c:prebound},
\begin{align*}
p_1(t,x,y) \ge c \frac{a_k^{\gamma}(1/\phi(t^{-1}),x,y)}{\phi(t^{-1})V(l)\Psi(l)}.
\end{align*}
Hence, we get the result from \eqref{e:mixedsoln}.
\qed

\subsection{Truncated kernel}
In this subsection, we give the proof for Theorem \ref{t:main2}. Throughout this subsection, we further assume that condition {\bf (Trunc.)($t_f$)} holds.

\begin{prop}\label{p:pretrunc}
There are comparison constants independent of $x$ and $y$ such that for all $t  \ge (\lfloor d_2/\alpha_1 + 2 \gamma \rfloor \vee 1/2)t_f$, it holds that
\begin{align*}
p(t,x,y) \simeq q(ct,x,y).
\end{align*}
\end{prop}
\proof 
Note that by Lemma \ref{l:phiHw}(i), $\phi(t^{-1}) \asymp t^{-1}$ for all $t \ge t_f$. Thus, by Corollary \ref{c:prebound}, we obtain the lower bound.
Since condition {\bf (Trunc.)($t_f$)} implies condition {\bf (Sub.)($1, 1$)}, by Theorem \ref{t:mainsub}, there exists a constant $a>0$ such that if $\lambda =0$ and $a \Phi(\p(x,y)) \ge t$, then $p(t,x,y) \simeq q(ct,x,y)$. Moreover, if $\lambda>0$, then since $D$ is bounded, by taking $a$ small enough, we can assume that there is no $x,y \in D$ such that $a \Phi(\p(x,y)) \ge t$. Hence, it remains to prove the upper bound when $a \Phi(\p(x,y)) < t$. Assume that $a \Phi(\p(x,y)) < t$.

 Let $r_0$ and $L$ be the constants in Propositions \ref{p:trunc} and \ref{l:trunc2}, respectively. Using the same arguments as in the ones given in the proof of Theorem \ref{t:mainsmall},
\begin{align*}
p(t,x,y) &\asymp \int_{0}^{Lt} q(r,x,y) d_r \P(S_r \ge t) - \int_{Lt}^{\infty} q(r,x,y) d_r \P(S_r \le t) \\
& \le cq(Lt,x,y) + c\int_{aL\Phi(l)/2}^{Lt} r^{-1}q(r,x,y)\P(S_r \ge t) dr.
\end{align*}
 
\noindent {\it Case 1.} $\lambda = 0$;

\smallskip

If $aL\Phi(l)/2 \ge r_0$, then by Lemma \ref{l:trunc2} and the fact that $r \mapsto r^{2\gamma} a^\gamma_k(r,x,y)$ is increasing, 
\begin{align*}
&\int_{aL\Phi(l)/2}^{Lt} r^{-1}q(r,x,y)\P(S_r \ge t) dr \le c\int_{r_0}^{Lt} r^{-1-2\gamma}\frac{r^{2\gamma}a^\gamma_k (r,x,y)}{V(\Phi^{-1}(r))}\left(\frac{r}{t}\right)^{ct} dr \\
&\quad  \le c t^{2\gamma} a^\gamma_k(t,x,y) L^{ct} \int_{r_0}^{Lt} dr \le c a^\gamma_k(t,x,y) e^{-ct} \le c\frac{a^\gamma_k(t,x,y)}{V(\Phi^{-1}(t))} \asymp q(Lt,x,y).
\end{align*}

Otherwise, if $aL\Phi(l)/2 < r_0$, then by Propositions \ref{p:trunc} and \ref{l:trunc2} and the weak scaling properties of $V$ and $\Phi$,
\begin{align*}
&\int_{aL\Phi(l)/2}^{Lt} r^{-1}q(r,x,y)\P(S_r \ge t) dr \\
&\le c\exp \big(-ct \log t \big)\int_{aL\Phi(l)/2}^{r_0} \frac{r^{\lfloor t/t_f \rfloor + 2\gamma}a^\gamma_k (r,x,y)}{r^{2\gamma}V(\Phi^{-1}(r))}dr+ c\int_{r_0}^{Lt} r^{-1-2\gamma}\frac{r^{2\gamma}a^\gamma_k (r,x,y)}{V(\Phi^{-1}(r))}\left(\frac{r}{t}\right)^{ct} dr \\
&\le c\frac{a^\gamma_k(t,x,y)}{V(\Phi^{-1}(t))} \left(1 + \int_{aL\Phi(l)/2}^{r_0} r^{\lfloor t/t_f \rfloor - 2\gamma - d_2/\alpha_1}dr \right) \le c\frac{a^\gamma_k(t,x,y)}{V(\Phi^{-1}(t))} \asymp q(Lt,x,y).
\end{align*} 
In the last inequality, we used the assumption that $t/t_f \ge \lfloor d_2/\alpha_1 + 2 \gamma \rfloor$.

\smallskip

\noindent {\it Case 2.} $\lambda > 0 $;

\smallskip

If $aL\Phi(l)/2 \ge r_0$, then by Lemma \ref{l:trunc2},
\begin{align*}
&\int_{aL\Phi(l)/2}^{Lt} r^{-1}q(r,x,y)\P(S_r \ge t) dr \le c\Phi(\delta_D(x))^\gamma \Phi(\delta_D(y))^\gamma \int_{r_0}^{Lt} r^{-1}e^{-\lambda r}\left(\frac{r}{t}\right)^{ct} dr \\
&\quad  \le c\Phi(\delta_D(x))^\gamma \Phi(\delta_D(y))^\gamma L^{ct} \int_{r_0}^{Lt} dr \le c \Phi(\delta_D(x))^\gamma \Phi(\delta_D(y))^\gamma e^{-ct} \simeq q(ct,x,y).
\end{align*}

Otherwise, if $aL\Phi(l)/2 < r_0$, then by Propositions \ref{p:trunc} and \ref{l:trunc2} and the above calculation,
\begin{align*}
&\int_{aL\Phi(l)/2}^{Lt} r^{-1}q(r,x,y)\P(S_r \ge t) dr \\
&\le ce^{-ct\log t}\int_{aL\Phi(l)/2}^{r_0} \frac{r^{\lfloor t/t_f \rfloor + 2\gamma}a^\gamma_k (r,x,y)}{r^{2\gamma}V(\Phi^{-1}(r))}dr+ c\Phi(\delta_D(x))^\gamma \Phi(\delta_D(y))^\gamma e^{-ct}\\
&\le c \Phi(\delta_D(x))^\gamma \Phi(\delta_D(y))^\gamma e^{-ct \log t} \int_{aL\Phi(l)/2}^{r_0} r^{\lfloor t/t_f \rfloor - 2\gamma - d_2/\alpha_1}dr +c\Phi(\delta_D(x))^\gamma \Phi(\delta_D(y))^\gamma e^{-ct} \\
& \le c \Phi(\delta_D(x))^\gamma \Phi(\delta_D(y))^\gamma e^{-ct} \simeq q(ct,x,y).
\end{align*}

\qed

\noindent{\bf Proof of Theorem \ref{t:main2}.} By Proposition \ref{p:pretrunc} and the second paragraph in its proof, it remains to consider the case when $\Phi(l) \le t \le\lfloor d_2/\alpha_1 + 2 \gamma \rfloor t_f$. Then, by using Proposition \ref{p:trunc} instead of Proposition \ref{p:utail1}, we get the result by the same argument as in the proof for Theorem \ref{t:mainsmall}. We omit in here.
\qed

\vspace{2mm}

\noindent{\bf Proof of Theorems \ref{t:specialsmall}, \ref{t:speciallarge}, \ref{t:specialsub} and \ref{t:specialtrunc}.} Let $\Phi_\alpha(x):=x^\alpha$. Then, we can check that {\bf (J1)} equals to ${\bf HK}_J^{1/2,\lambda,1}(\Phi_\alpha, \Phi_\alpha)$, {\bf (J2)} equals to ${\bf HK}_J^{1/2,0,1}(\Phi_\alpha, \Phi_\alpha)$, {\bf (J3)} equals to ${\bf HK}_J^{1/2,0,2}(\Phi_\alpha, \Phi_\alpha)$, {\bf (J4)} equals to ${\bf HK}_J^{(\alpha-1)/\alpha,\lambda,1}(\Phi_\alpha, \Phi_\alpha)$, {\bf (D1)} equals to ${\bf HK}_D^{1/2,\lambda,1}(\Phi_\alpha)$, {\bf (D2)} equals to ${\bf HK}_D^{1/2,0,1}(\Phi_\alpha)$ and {\bf (D3)} equals to ${\bf HK}_D^{1/2,\lambda,2}(\Phi_\alpha)$ where the underlying function $V(x,r):=r^d$ for all $x \in D$ and $r>0$. Hence, we can apply Theorems \ref{t:mainsmall}, \ref{t:mainlarge}, \ref{t:mainsub} and \ref{t:main2}. Combining these results with Proposition \ref{p:Dgamma} and Remark \ref{r:sJ2}, we get the result.
\qed

\section{Appendix}

In this section, we give the sketch of proof of Proposition \ref{p:Dgamma}. Fix $t>0$ and $x,y \in D$ satisfying $\Phi(\p(x,y)) \phi(t^{-1}) \le 1/(4e^2)$ and set $V(r):=V(x,r)$ and $l:=\p(x,y)$ as before.

\begin{lemma}\label{l:int}
Fix $p \in \R$. For $0<A<B/2$, define
\begin{align*}
S_p(A,B):=\int_A^B \frac{1}{r^pV(\Phi^{-1}(r))}dr.
\end{align*}
Then, the followings are true.

(i) There exists a constant $c>0$ independent of $A$ and $B$ such that
\begin{align*}
S_p(A,B) \ge c \big(A^{1-p}V(\Phi^{-1}(A))^{-1}  + B^{1-p}V(\Phi^{-1}(B))^{-1} \big).
\end{align*}

(ii) If $d_1>\alpha_2(1-p)$, then $\displaystyle S_p(A,B) \asymp A^{1-p}V(\Phi^{-1}(A))^{-1}$.

(iii) If $d_2<\alpha_1(1-p)$, then $\displaystyle S_p(A,B) \asymp B^{1-p}V(\Phi^{-1}(B))^{-1}$.

(iv) If $d_1 = d_2 = (1-p)\alpha_1 = (1-p)\alpha_2$, then $\displaystyle S_p(A,B) \asymp \log(B/A)$.
\end{lemma}
\proof (i) By the monotonicities and the weak scaling properties of $V$ and $\Phi$, 
\begin{align*}
2S_p(A,B) &\ge \int_A^{2A} \frac{1}{r^pV(\Phi^{-1}(r))}dr + \int_{B/2}^B \frac{1}{r^pV(\Phi^{-1}(r))}dr \\
&\ge \frac{A^{1-p}}{2^pV(\Phi^{-1}(2A))} + \frac{B^{1-p} }{2V(\Phi^{-1}(B))} \ge c \big(A^{1-p}V(\Phi^{-1}(A))^{-1}  + B^{1-p}V(\Phi^{-1}(B))^{-1} \big).
\end{align*}

(ii), (iii) See \cite[2.12.16]{BGT}.

(iv) In this case, since the assumptions imply that $V(r) \asymp r^{d_1}$ and $\Phi^{-1}(r) \asymp r^{1/\alpha_1}$ for all $r>0$, we get $S_p(A,B)\asymp \int_A^B r^{-p-d_1/\alpha_1} dr = \int_A^B r^{-1}dr = \log(B/A)$.
\qed

\vspace{1mm}

Recall that $\delta_*^\Phi(x,y)=\Phi(\delta_{D}(x))\Phi(\delta_{D}(y))$. Without loss of generality, by symmetry, we can assume that $\delta_D(x) \le \delta_D(y)$. We first claim that if $\Phi(l)\phi(t^{-1}) \le 1/(4e^2)$, then
\begin{align*}
\big(\phi(t^{-1})^{-1}+\Phi(\delta_D(x))\big)\big(\phi(t^{-1})^{-1}+\Phi(\delta_D(y))\big) \asymp \phi(t^{-1})^{-2} + \delta_*^\Phi(x,y).
\end{align*}
Indeed, it is clear that $(RHS) \le (LHS)$ and we also have that
\begin{align*}
(LHS) &\le \phi(t^{-1})^{-2} + \delta_*^\Phi(x,y) + 2\phi(t^{-1})^{-1} \Phi(\delta_D(x)+l) \\
&\le \phi(t^{-1})^{-2} + \delta_*^\Phi(x,y) + 2\phi(t^{-1})^{-1}\big(\Phi(2\delta_D(x)) + \Phi(2l)\big) \\
& \le c\phi(t^{-1})^{-2} + \delta_*^\Phi(x,y) +c\Phi(\delta_D(x))^2 +c\Phi(l)^2 \le c(RHS).
\end{align*}
In the third line, we used the fact that $2ab \le a^2 + b^2$ for $a, b \in \R$, the weak scaling properties of $\Phi$ and the assumption that $\phi(t^{-1})^{-1} \ge 4e^2 \Phi(l)$. Thus, if $\Phi(l)\phi(t^{-1}) \le 1/(4e^2)$, then
\begin{align*}
a_1^{\gamma}(1/\phi(t^{-1}), x, y)& =\left( \frac{\delta_*^\Phi(x,y)}{\big(\phi(t^{-1})^{-1}+\Phi(\delta_D(x))\big)\big(\phi(t^{-1})^{-1}+\Phi(\delta_D(y))\big)} \right)^{\gamma}\\
& \asymp \left(\frac{\delta_*^\Phi(x,y)}{\phi(t^{-1})^{-2} + \delta_*^\Phi(x,y)}\right)^\gamma \asymp \left(1 \wedge \frac{\delta_*^\Phi(x,y)}{\phi(t^{-1})^{-2 }}\right)^\gamma  \asymp \left(1 \wedge \frac{\delta_*^\Phi(x,y)^\gamma }{\phi(t^{-1})^{-2\gamma }}\right). 
\end{align*}

Now, We consider the following three scenarios.

\smallskip

(Sc.1) $\Phi(\delta_D(x)) \le 4\Phi(l)$.

(Sc.2) $4\Phi(l) < \Phi(\delta_D(x))$ and $\Phi(\delta_D(y)) \le 1/(4e^2\phi(t^{-1}))$.

(Sc.3)  $4\Phi(l) < \Phi(\delta_D(x))$ and $\Phi(\delta_D(y)) > 1/(4e^2\phi(t^{-1}))$.

\smallskip

If (Sc.1) is true, then we have
\begin{align*}
\sI_1^{\gamma}(t,x,y) \asymp \delta_*^\Phi(x,y)^\gamma  S_{2\gamma}\big(\Phi(l),1/(2e^2\phi(t^{-1}))\big).
\end{align*}

Else if (Sc.2) is true, then we have
\begin{align*}
\sI_1^{\gamma}(t,x,y) &\asymp S_0\big(\Phi(l), \Phi(\delta_D(x))/2\big) + \Phi(\delta_D(x))^{\gamma}S_\gamma \big(\Phi(\delta_D(x))/2, \Phi(\delta_D(y))\big)\\
&\quad + \delta_*^\Phi(x,y)^\gamma S_{2\gamma}\big(\Phi(\delta_D(y)), 1/(2e^2\phi(t^{-1}))\big).
\end{align*}

Otherwise, if (Sc.3) is true, then we get
\begin{align*}
F_1^{\gamma}(t,x,y) \asymp F_1^0(t,x,y) \asymp S_0 \big(\Phi(l),1/(2e^2\phi(t^{-1}))\big).
\end{align*}

Hence, by applying Lemma \ref{l:int} with $p=0, \gamma$ and $2\gamma$, we obtain the following estimates.
\vspace{1mm}

\noindent (a) Suppose that $d_2/\alpha_1<1-2\gamma$. Then,
\begin{align*}
\sI^\gamma_1(t,x,y) \asymp \begin{cases}
\delta_*^\Phi(x,y)^\gamma \phi(t^{-1})^{2\gamma-1}V\big(\Phi^{-1}(1/\phi(t^{-1}))\big)^{-1},    &\;\; \text{if} \; \text{(Sc.1) is true};\\
\delta_*^\Phi(x,y)^\gamma \phi(t^{-1})^{2\gamma-1}V\big(\Phi^{-1}(1/\phi(t^{-1}))\big)^{-1}, &\;\; \text{if} \; \text{(Sc.2) is true};\\
\phi(t^{-1})^{-1}V\big(\Phi^{-1}(1/\phi(t^{-1}))\big)^{-1},  &\;\; \text{if} \; \text{(Sc.3) is true}.
\end{cases}
\end{align*}

\noindent (b) Suppose that $\alpha_1=\alpha_2$, $d_1=d_2=(1-2\gamma)\alpha_1$ and $\gamma>0$. Then, $V(r) \asymp r^{d_1}, \Phi(r) \asymp r^{\alpha_1}$ and
\begin{align*}
\sI^\gamma_1(t,x,y) \asymp \begin{cases}
\delta_*^\Phi(x,y)^\gamma \displaystyle\log \left(\frac{1}{\Phi(l)\phi(t^{-1})}\right),    &\;\; \text{if} \; \text{(Sc.1) is true};\\
\delta_*^\Phi(x,y)^\gamma \displaystyle\log \left(\frac{1}{\Phi(\delta_D(y))\phi(t^{-1})}\right), &\;\; \text{if} \; \text{(Sc.2) is true};\\
\phi(t^{-1})^{-1}V\big(\Phi^{-1}(1/\phi(t^{-1}))\big)^{-1},  &\;\; \text{if} \; \text{(Sc.3) is true}.
\end{cases}
\end{align*}

\noindent (c) Suppose that $1-2\gamma<d_1/\alpha_2 \le d_2/\alpha_1<1-\gamma$. Then,
\begin{align*}
\sI^\gamma_1(t,x,y) \asymp \begin{cases}
\delta_*^\Phi(x,y)^\gamma \Phi(l)^{1-2\gamma}V(l)^{-1},    &\;\; \text{if} \; \text{(Sc.1) is true};\\
\delta_*^\Phi(x,y)^\gamma \Phi(\delta_D(y))^{1-2\gamma} V(\delta_D(y))^{-1}, &\;\; \text{if} \; \text{(Sc.2) is true};\\
\phi(t^{-1})^{-1}V\big(\Phi^{-1}(1/\phi(t^{-1}))\big)^{-1},  &\;\; \text{if} \; \text{(Sc.3) is true}.
\end{cases}
\end{align*}

\noindent (d) Suppose that $\alpha_1=\alpha_2$, $d_1=d_2=(1-\gamma)\alpha_1$ and $\gamma>0$. Then, $V(r) \asymp r^{d_1}, \Phi(r) \asymp r^{\alpha_1}$ and
\begin{align*}
\sI^\gamma_1(t,x,y) \asymp \begin{cases}
\delta_*^\Phi(x,y)^\gamma \Phi(l)^{1-2\gamma}V(l)^{-1},    &\;\; \text{if} \; \text{(Sc.1) is true};\\
\Phi(\delta_D(x))^\gamma \displaystyle\log \left(\frac{2\Phi(\delta_D(y))}{\Phi(\delta_D(x))}\right), &\;\; \text{if} \; \text{(Sc.2) is true};\\
\phi(t^{-1})^{-1}V\big(\Phi^{-1}(1/\phi(t^{-1}))\big)^{-1},  &\;\; \text{if} \; \text{(Sc.3) is true}.
\end{cases}
\end{align*}

\noindent (e) Suppose that $1-\gamma<d_1/\alpha_2 \le d_2/\alpha_1<1$. Then,
\begin{align*}
\sI^\gamma_1(t,x,y) \asymp \begin{cases}
\delta_*^\Phi(x,y)^\gamma \Phi(l)^{1-2\gamma}V(l)^{-1},    &\;\; \text{if} \; \text{(Sc.1) is true};\\
\Phi(\delta_D(x))V(\delta_D(x))^{-1}, &\;\; \text{if} \; \text{(Sc.2) is true};\\
\phi(t^{-1})^{-1}V\big(\Phi^{-1}(1/\phi(t^{-1}))\big)^{-1},  &\;\; \text{if} \; \text{(Sc.3) is true}.
\end{cases}
\end{align*}

\noindent (f) Suppose that $d_1=d_2=\alpha_1 =\alpha_2$. Then, $V(r) \asymp r^{d_1}, \Phi(r) \asymp r^{\alpha_1}$ and
\begin{align*}
\sI^\gamma_1(t,x,y) \asymp \begin{cases}
\delta_*^\Phi(x,y)^\gamma \Phi(l)^{-2\gamma},    &\;\; \text{if} \; \text{(Sc.1) is true};\\
\displaystyle\log \left(\frac{\Phi(\delta_D(x))}{\Phi(l)}\right), &\;\; \text{if} \; \text{(Sc.2) is true};\\
\displaystyle\log \left(\frac{1}{\Phi(l)\phi(t^{-1})}\right),  &\;\; \text{if} \; \text{(Sc.3) is true}.
\end{cases}
\end{align*}

\noindent (g) Suppose that $1<d_1/\alpha_2$. Then,
\begin{align*}
\sI^\gamma_1(t,x,y) \asymp \begin{cases}
\delta_*^\Phi(x,y)^\gamma \Phi(l)^{1-2\gamma}V(l)^{-1},    &\;\; \text{if} \; \text{(Sc.1) is true};\\
\Phi(l)V(l)^{-1}, &\;\; \text{if} \; \text{(Sc.2) is true};\\
\Phi(l)V(l)^{-1},  &\;\; \text{if} \; \text{(Sc.3) is true}.
\end{cases}
\end{align*}

Together with the fact that $\phi(t^{-1}) \ge t^{-1}\int_0^t e^{-s/t}w(s)ds \ge e^{-1}w(t)$, we get the result.

\section{Examples}

\begin{example}\label{example1}
{\rm

(cf. \cite[Example 2.5(ii)]{Ch}.) Let $0<\alpha \le 2$, $0<\beta < 1$ and $\delta>0$. Then, we consider the fundamental solution of the following Cauchy problem.

\begin{align}\label{e:trunc}
\frac{d}{dt} \int_{(t- \delta) \vee 0}^t \big[(t-s)^{-\beta} - \delta^{-\beta} \big] \big(u(s,x)-f(x)\big) ds = \Delta^{\alpha/2} u(t,x),  \qquad  & x \in \R^d, \;\; t>0, \nn\\
u(0, x) = f(x), \qquad & x \in \R^d.
\end{align}

In this case, we see that $w(s)=w_{\delta}(s) = (s^{-\beta} - \delta^{-\beta}) \1_{(0, \delta]}(s)$ and hence conditions {\bf (Ker.)} and {\bf (Trunc.)($\delta$)} hold. Moreover, it is well known that for the function $\Phi_{\alpha}(x) = x^\alpha$, the heat kernel $q(t,x,y)$ corresponding to the generator $\Delta^{\alpha/2}$ enjoys estimate ${\bf HK}^{0,0,0}_J(\Phi_{\alpha}, \Phi_{\alpha})$ if $0< \alpha <2$ and estimate ${\bf HK}^{0,0,0}_D(\Phi_{\alpha})$ if $\alpha = 2$. By Theorems \ref{t:mainsmall} and \ref{t:main2}, we obtain the global estimates for the fundamental solution $p(t,x,y)$ of the equation \eqref{e:trunc}.

\smallskip

(i) For every $t \in (0,\delta/2]$ and $x,y \in \R^d$, we have
\begin{align*}
p(t,x,y) \simeq \begin{cases}
t^{-\beta d/\alpha}, \quad & \mbox{if} \;\; |x-y| \le t^{\beta/\alpha} \;\; \text{and} \;\; d<\alpha; \\
t^{-\beta } \log\big( 2t^{\beta/\alpha}/|x-y|\big), \quad &\mbox{if} \;\; |x-y| \le t^{\beta/\alpha} \;\; \text{and} \;\; d=\alpha; \\
t^{-\beta}/ |x-y|^{d-\alpha}, \quad  &\mbox{if} \;\; |x-y| \le t^{\beta/\alpha} \;\; \text{and} \;\; d>\alpha; \\
t^{\beta}/ |x-y|^{d + \alpha }, \quad &\mbox{if} \;\; |x-y| > t^{\beta/\alpha} \;\; \text{and} \;\; 0< \alpha < 2; \\
t^{-\beta d/\alpha} \exp \big(-c|x-y|^{2/(2-\beta)} t^{-\beta/(2-\beta)}\big), \quad& \mbox{if} \;\; |x-y| > t^{\beta/\alpha} \; \; \text{and} \;\; \alpha =2. \\
\end{cases}
\end{align*}

(ii) Fix any $t \in [\delta/2, \infty)$ and $x,y \in \R^d$. Let $n_t = \lfloor t/\delta \rfloor + 1$. Then, we have
\begin{align*}
p(t,x,y) \simeq \begin{cases}
\big[|x-y|^\alpha t^{-1}+ (n_t \delta - t)^{n_t} \big]t^{-n_t}/|x-y|^{d-\alpha n_t} , & \;\; \mbox{if} \;\; |x-y|^\alpha \le t \;\; \text{and} \\
& \quad   \quad  \; \delta/2 \le t< \lfloor (d-\alpha)/\alpha \rfloor \delta; \\
t^{-d/\alpha}+ (n_t \delta - t)^{n_t} t^{-n_t}/|x-y|^{d-\alpha n_t} , & \;\; \mbox{if} \;\; d/\alpha \notin \N,\;\; |x-y|^\alpha \le t \;\; \text{and} \\
& \quad   \quad  \lfloor (d-\alpha)/\alpha \rfloor \delta \le t< \lfloor d/\alpha \rfloor \delta; \\
t^{-d/\alpha}+ (\displaystyle\frac{d\delta}{\alpha t}- 1)^{d/\alpha}\log \big( 2t /|x-y|^\alpha \big), & \;\; \mbox{if} \;\; d/\alpha \in \N,\;\; |x-y|^\alpha \le t \;\; \text{and} \\
& \quad   \quad  (d-\alpha)\delta/\alpha   \le t< d\delta/\alpha; \\
t^{-d/\alpha}, &   \;\;\mbox{if} \;\; |x-y|^\alpha \le t \;\; \text{and} \; \lfloor d/\alpha \rfloor \delta \le t; \\
t/|x-y|^{d+\alpha}, &  \;\;\mbox{if} \;\;  |x-y|^\alpha > t \;\; \text{and} \;\; 0< \alpha < 2, \\
t^{- d/\alpha} \exp \big(-c|x-y|^{2} t^{-1}\big), &  \;\;\mbox{if} \;\; |x-y|^\alpha >t \;\; \text{and} \;\; \alpha = 2.
\end{cases}
\end{align*}

\vspace{2mm}

In particular, for every $t>0$ and $x \in \R^d$, $p(t,x,x) < \infty$ if and only if $t \ge \lfloor d / \alpha \rfloor \delta$.

}
\end{example}

\smallskip

\begin{example}\label{example2}
{\rm

Let $d \ge 1$, $0< \alpha \le 2$ and $D \subset \R^d$ be a bounded $C^{1,1}$ open set. When $\alpha = 2$, we further assume that $D$ is connected. Let $\kappa:(0,1) \to [0, \infty)$ be a measurable function with $\int_0^1 \kappa(\beta) d\beta < \infty$. In this example, we consider the following fractional-time equation.
\begin{align}\label{e:mixCaputo}
&\int_0^1 \partial_t^\beta u(t,x) \kappa (\beta) d \beta = \Delta^{\alpha/2} u(t,x),  \qquad   x \in D, \;\; t>0, \nn\\
&u(t, x) = 0, \quad x \in \R^d \setminus D, \;\; t>0,   \qquad u(0, x) = f(x), \quad  x \in D,
\end{align}
where $\partial^\beta_t$ is the Caputo derivative of order $\beta$. Such distributed-order fractional equations were studied in \cite{MS}. In this case, we can check that the integral kernel for fractional-time derivative is given by $w(s) = \int_0^1 s^{-\beta} \frac{\kappa(\beta)}{\Gamma(1-\beta)} d \beta$ and hence condition {\bf (Ker.)} holds. Moreover, since $2w(2s) \ge w(s)$ for all $s>0$, both conditions {\bf (S.Poly.)($t_s$)} and {\bf (L.Poly.)} hold.

By \cite{CKS}, \cite{Zh} and \cite{So2}, we see that the transition density $q(t,x,y)$ corresponding to the generator $\Delta^{\alpha/2} |_D$ satisfies the estimates ${\bf HK}_J^{1/2, \lambda, 1}(\Phi_\alpha, \Phi_\alpha)$ if $0<\alpha<2$ and ${\bf HK}_D^{1/2, \lambda, 1}(\Phi_\alpha)$ if $\alpha =2$ where $-\lambda<0$ is the largest eigenvalue of the generator $\Delta^{\alpha/2} |_D$. Then, by Theorems \ref{t:specialsmall} and \ref{t:speciallarge}, we obtain the sharp estimates for the fundamental solution $p(t,x,y)$ of the equation \eqref{e:mixCaputo}. Note that $\phi(t) = \int_0^1 t^{\beta} \kappa(\beta)d \beta$ in this case.

\smallskip

(i) For every $t \in (0,1]$ and $x,y \in D$ satisfying $|x-y|^\alpha \int_0^1 t^{-\beta} \kappa(\beta)d \beta \le 1/(4e^2)$, we have
\begin{align*}
p(t,x,y)& \asymp \bigg(1 \wedge \frac{\delta_D(x)^{\alpha/2}\delta_D(y)^{\alpha/2}}{\big(\int_0^1 t^{-\beta} \kappa(\beta)d \beta\big)^{-1} }\bigg)\bigg(\int_0^1 t^{-\beta} \kappa(\beta)d \beta\bigg)^{d/\alpha} \\
& \quad+ \bigg(1 \wedge \frac{\delta_D(x)^{\alpha/2}\delta_D(y)^{\alpha/2}}{|x-y|^\alpha}\bigg) F^\alpha_k(d,t,x,y) \int_0^1 t^{-\beta} \frac{\kappa(\beta)}{\Gamma(1-\beta)} d \beta,
\end{align*}
where $F^\alpha_k$ is the function in section \ref{s:special}.

(ii) For every $t \in (0,1]$ and $x,y \in D$ satisfying $|x-y|^\alpha \int_0^1 t^{-\beta} \kappa(\beta)d \beta > 1/(4e^2)$, we have
\begin{align*}
p(t,x,y) \asymp \bigg(1 \wedge \frac{\delta_D(x)^{\alpha/2}}{\big(\int_0^1 t^{-\beta} \kappa(\beta)d \beta\big)^{-1/2} }\bigg) \bigg(1 \wedge \frac{\delta_D(y)^{\alpha/2}}{\big(\int_0^1 t^{-\beta} \kappa(\beta)d \beta\big)^{-1/2} }\bigg) \frac{\big(\int_0^1 t^{-\beta} \kappa(\beta)d \beta\big)^{-1}}{|x-y|^{d+\alpha}},
\end{align*}
if $0<\alpha<2$ and
\begin{align*}
p(t,x,y)& \simeq  \bigg(1 \wedge \frac{\delta_D(x)}{\big(\int_0^1 t^{-\beta} \kappa(\beta)d \beta\big)^{-1/2} }\bigg) \bigg(1 \wedge \frac{\delta_D(y)}{\big(\int_0^1 t^{-\beta} \kappa(\beta)d \beta\big)^{-1/2} }\bigg) \\
&\quad \times \bigg(\int_0^1 t^{-\beta} \kappa(\beta)d \beta\bigg)^{d/2} \exp \bigg(-ct\bar{\kappa}(t, |x-y|)\bigg),
\end{align*}
if $\alpha =2$ where the function $\bar{\kappa}(t,l):= \sup\{s>0 : \int_0^1 s^{\beta-2} \kappa(\beta) d \beta > t^2l^{-2} \}$.

(iii) For every $t \in [1, \infty)$ and $x,y \in D$, we have
\begin{align*}
p(t, x, y) &\asymp \left(\int_0^1 t^{-\beta} \frac{\kappa(\beta)}{\Gamma(1-\beta)} d \beta \right) \left(1 \wedge \displaystyle\frac{\delta_D(x)\delta_D(y)}{\p(x,y)^{2}}\right)^{\alpha/2}\bigg( 1 \wedge \delta_D(x)^{\alpha/2}\delta_D(y)^{\alpha/2} + F^\alpha_k(d,T_D,x,y) \bigg),
\end{align*}
where $T_D=[\phi^{-1}(4^{-1}e^{-2}R_D^{-\alpha})]^{-1}$ and $F^\alpha_k$ is the function in section \ref{s:special}.

}
\end{example}

Following \cite{CKS12}, for a function $f$ on $\R^d$, we define for $1<\alpha<2$ and $r>0$,
\begin{align*}
M^\alpha_f:= \sup_{x \in \R^d} \int_{|y-x|<r} \frac{|f(y)|}{|x-y|^{d+1-\alpha}} dy.
\end{align*}
Then, a function $f$ on $\R^d$ is said to belong to the Kato class $\mathbb{K}^{\alpha-1}$ if $\lim_{r \to 0+} M^\alpha_f(r) = 0$.

\begin{example}
{\rm Let $d \ge 2$, $1< \alpha \le 2$ and $D \subset \R^d$ be a bounded $C^{1,1}$ open set.
 In \cite{CKS12}, the authors studied the stability of Dirichlet heat kernel estimates under gradient perturbation. More precisely, for every $b \in \mathbb{K}^{\alpha-1}$, an operator $\big(\Delta^{\alpha/2} + b \cdot \nabla \big) |_D$ satisfies the estimates ${\bf HK}_J^{1/2, \lambda^b, 1}(\Phi_\alpha, \Phi_\alpha)$ for some constant $\lambda^b>0$. Notice that the estimates in Example \ref{example2} is independent of $-\lambda<0$. Therefore, we can deduce that the results in Example \ref{example2} still works not only with the operator $\Delta^{\alpha/2}$ but also operators $\Delta^{\alpha/2} + b \cdot \nabla$ for $b \in \mathbb{K}^{\alpha-1}$.
}
\end{example}

\begin{example}\label{example3}
{\rm
Since our theorem covers when $q(t,x,y)$ enjoys a mixed type estimates, we also obtain the estimates for fundamental solution with repect to the operators $\Delta + \Delta^{\beta/2}$ for $0<\beta<2$ in Examples \ref{example1} and \ref{example2}. Indeed, these are nothing but sum of two results for $\alpha = 2$ and $\alpha = \beta$.
}
\end{example}

\end{document}